\documentclass[12pt,leqno]{article}

\usepackage{amsmath,amssymb,amsfonts,theorem,latexsym}

\newtheorem{thm}{Theorem}
\newtheorem{theo}{Theorem}[section]
\newtheorem{lem}[theo]{Lemma}
\newtheorem{cor}[theo]{Corollary}
\newtheorem{prop}[theo]{Proposition}
\newtheorem{defi}[theo]{Definition}
\newtheorem{claim}[theo]{Claim}

{\theorembodyfont{\rmfamily} \newtheorem{rem}[theo]{Remark}}
{\theorembodyfont{\rmfamily} \newtheorem{ex}[theo]{Example}}
{\theorembodyfont{\rmfamily} \newtheorem{contra}[theo]{Counter Example}}
{\theorembodyfont{\rmfamily} \newtheorem{con}[theo]{Conjecture}}

\usepackage[all]{xy}

\newcommand{\Appendix}[1]{%
  \refstepcounter{section}%
  \addtocontents{toc}{\protect\setcounter{tocdepth}{1}}
  \addcontentsline{toc}{section}%
    {\bfseries\appendixname~\thesection\ #1}%
    {\medskip\noindent \Large\bfseries\appendixname\ \thesection\ #1}%
\sectionmark{#1}\smallskip\noindent
\renewcommand{\theequation}{{\bf 
{{\thesection}}.{\arabic{equation}}}}
}

\DeclareFontFamily{U}{rsf}{}
\DeclareFontShape{U}{rsf}{m}{n}{
  <5> <6> rsfs5 <7> <8> <9> rsfs7 <10->  rsfs10}{}
\DeclareMathAlphabet{\mathscr}{U}{rsf}{m}{n}
\newcommand{\mycal}[1]{\mathscr{#1}}

\newlength{\unten}
\newlength{\gesamt}
\newcommand{\dhoehe}[1]{%
\protect\settoheight{\gesamt}{\(\displaystyle#1\)}%
\protect\settodepth{\unten}{\(\displaystyle#1\)}%
\protect\addtolength{\gesamt}{\unten}}
\newcommand{\thoehe}[1]{%
\protect\settoheight{\gesamt}{\(\textstyle #1\)}%
\protect\settodepth{\unten}{\(\textstyle #1\)}%
\protect\addtolength{\gesamt}{\unten}}
\newcommand{\shoehe}[1]{%
\protect\settoheight{\gesamt}{\(\scriptstyle #1\)}%
\protect\settodepth{\unten}{\(\scriptstyle #1\)}%
\protect\addtolength{\gesamt}{\unten}}
\newcommand{\sshoehe}[1]{%
\protect\settoheight{\gesamt}{\(\scriptscriptstyle #1\)}%
\protect\settodepth{\unten}{\(\scriptscriptstyle #1\)}%
\protect\addtolength{\gesamt}{\unten}}
\newcommand{\leftidx}[3]{%
\mathchoice%
{\dhoehe{#2}\rule[-\unten]{0mm}{\gesamt}#1#2#3}%
{\thoehe{#2}\rule[-\unten]{0mm}{\gesamt}#1#2#3}%
{\shoehe{#2}\rule[-\unten]{0mm}{\gesamt}#1#2#3}%
{\sshoehe{#2}\rule[-\unten]{0mm}{\gesamt}#1#2#3}}

\DeclareFontFamily{U}{cyr}{}
\DeclareFontShape{U}{cyr}{m}{n}{
  <5> wncyr5 <6> wncyr6 <7> wncyr7 <8> wncyr8 <9> wncyr9 <10->
wncyr10}{}
\DeclareMathAlphabet{\mathcyr}{U}{cyr}{m}{n}

\input cyracc.def
\input epsf

\newcommand{\op}[1]{\operatorname{#1}}
\newcommand{\TSh}{{\mathcyr{\cyracc Sh}}}
\newcommand{\grX}{\mycal{X}}
\newcommand{\tXc}[2]{{\leftidx{_{#1}}{\widetilde{X}}{_{#2}}}}
\newcommand{\Xc}[2]{{\leftidx{_{#1}}{X}{_{#2}}}}
\newcommand{\Xl}[2]{{\leftidx{_{#1}}{\mycal{L}}{_{#2}}}}
\newcommand{\Xe}[2]{{\leftidx{_{#1}}{\mycal{E}}{_{#2}}}}
\newcommand{\ul}[2]{{\leftidx{_{#1}}{LU}{_{#2}}}}
\newcommand{\ue}[2]{{\leftidx{_{#1}}{EU}{_{#2}}}}
\newcommand{\rl}[2]{{\leftidx{_{#1}}{LR}{_{#2}}}}
\newcommand{\re}[2]{{\leftidx{_{#1}}{ER}{_{#2}}}}
\newcommand{\ot}{{1\!-\!2, \, m\cdot 3}}
\newcommand{\unalpha}{\underline{\alpha}}
\newcommand{\FM}{\boldsymbol{F}\boldsymbol{M}}
\newcommand{\MS}{\boldsymbol{M}\boldsymbol{S}}
\newcommand{\ggl}{\boldsymbol{\gamma}_{L}}
\newcommand{\gge}{\boldsymbol{\gamma}_{E}}

\newcommand{\gp}[6]{
\xymatrix@R+10pt@C+15pt{
#1 \ar[d] \ar@<.5ex>[drr]^-{s} \ar@<-.5ex>[drr]_-{t} & & & \\
{#2}\times_{#3} {#2} \ar@<.5ex>[rr]^-{#4} \ar@<-.5ex>[rr]_-{#5} & & #2
\ar[r]^-{#6} & #3
}
}

\newcommand{\et}{et}
\newcommand{\bgamma}{\boldsymbol{\gamma}}

\newcommand{\bs}{\boldsymbol{s}}
\newcommand{\bt}{\boldsymbol{t}}
\newcommand{\bi}{\boldsymbol{i}}
\newcommand{\be}{\boldsymbol{e}}
\newcommand{\bm}{\boldsymbol{m}}
\newcommand{\Tr}{\boldsymbol{T}}

\newcommand{\sw}{\boldsymbol{s}\boldsymbol{w}}
\newcommand{\bF}{\boldsymbol{F}}
\newcommand{\bS}{\boldsymbol{S}}
\newcommand{\bG}{\boldsymbol{G}}
\newcommand{\mcup}{\boldsymbol{\doublecup}}

\newcommand{\bp}{\boldsymbol{p}}

\newcommand{\insubsection}[1]{
\refstepcounter{subsection}
\bigskip\noindent {\normalfont\large\bfseries
\arabic{section}.\arabic{subsection} \quad #1} 
\bigskip\noindent
}

\newcommand\calX{{\mycal{X}}}
\newcommand\calY{{\mycal{Y}}}
\newcommand\calP{{\mycal{P}}}
\newcommand\GM{{\mathbb G}_{m}}
\newcommand\REMARK[1]{{\begin{rem} #1 \end{rem}}}
\newcommand\CCC{{\mathbb C}}
\newcommand\ZZZ{{\mathbb Z}}
\newcommand\calO{{\mathcal O}}
\newcommand\DEFINITION[1]{{\begin{defi} #1 \end{defi}}}
\newcommand\PROPOSITION[1]{{\begin{prop} #1 \end{prop}}}
\newcommand\PROOF[1]{{{\bf Proof.} #1 \ \hfill $\Box$}} 
\newcommand\EXAMPLE[1]{{\begin{ex} #1 \end{ex}}}

\numberwithin{equation}{section}

\begin{document}
\title{Torus fibrations, gerbes, and duality}
\author{Ron Donagi\thanks{Partially supported by NSF Grants
DMS-0104354 and FRG-DMS-0139799.} \and Tony Pantev\thanks{Partially
supported by NSF grants DMS-0099715 and FRG-DMS-0139799, and an
A.P.Sloan Research Fellowship.}}
\date{ }
\maketitle

\tableofcontents

\newpage

\section{Introduction} \label{sec-introduction}

\insubsection{Duality for elliptic fibrations}\label{baby}

\noindent
In this paper we are concerned with categories of sheaves on
varieties fibered by genus one curves.  For an elliptic fibration on
$X$, by which we always mean a genus one fibration $\pi : X \to B$
admitting a holomorphic section $\sigma : B \to X$, there is by now a
a well understood theory of the Fourier-Mukai transform
\cite{mukai,bbhrmp,bridgeland-elliptic,bridgeland-maciocia}. The basic
result is:

\medskip

\

\noindent
{\bf Theorem \cite{bridgeland-maciocia}} {\em Let $\xymatrix@1{X
\ar[r]^-{\pi} & B \ar@/^0.5pc/[l]^-{\sigma}}$ be an elliptic fibration
with smooth total space. Then the integral transform (Fourier-Mukai
transform)
\[
\xymatrix@R-25pt{
\FM : D^{b}(X) \ar[r] & D^{b}(X) \\
F \ar[r] & Rp_{2*}(Lp_{1}^{*}F\otimes \mycal{P}),
}
\]
induced by the Poincare sheaf $\mycal{P} \to X\times_{B} X$, is an
auto-equivalence of the bounded derived category $D^{b}(X)$ of
coherent sheaves on $X$.
}

\

\medskip

\noindent
An important feature of $\FM$ is that it transforms geometric
objects in an interesting way:

\[
\xymatrix@R-15pt{
{}\save[]+<0cm,0.7cm>*\txt{$\left\{
\begin{minipage}[c]{5in}
{
{\em Bundle data}: vector bundles on $X$, semistable of degree zero on
the generic fiber of $\pi$.
}
\end{minipage}
\right\}
$}
\ar@{<->}[dd]^-{\FM} \restore  \\
\\
{}\save[]-<0cm,1cm>*\txt{$\left\{
\begin{minipage}[c]{5in}
{
{\em Spectral data:} sheaves on $X$ with the numerics
of a line bundle on a `spectral' divisor $C \subset X$,
finite over $B$.
}
\end{minipage}
\right\}
$} \restore
}
\]
\

This spectral construction was used to study general compactifications
of heterotic string theory and their moduli, and especially the duality
with F-theory \cite{fmw,bjps,donagi-principal,aspinwall-donagi,fmw-pb,
donagi-heteroticF,fmw-vb}.
It was used also to construct special bundles
on elliptic Calabi-Yau manifolds which lead to more-or-less realistic
compactified theories \cite{smb,sm}.

\insubsection{Partial duality for genus one fibrations}

\noindent
In many applications (see e.g. \cite{sm,smb}) one is also interested in
constructing bundles on genus one fibrations $\pi : Y \to B$ which do
not necessarily admit a section. From the viewpoint of the spectral
construction one expects that vector bundles on $Y$ should correspond
to spectral data supported on a divisor $C \subset X :=
\overline{\op{Pic}^{0}}(Y/B)$, where $\overline{\op{Pic}^{0}}(Y/B)$ is
the compactified relative Jacobian of $\pi : Y \to B$. However it is
unrealistic to expect that this spectral data should again be a sheaf on
$C$. One problem is that $X$ is not a fine moduli space of objects on
$Y$ and so we do not have a Poincare sheaf on $Y\times_{B} X$.
Another problem is that in the passage from $Y$ to $X$ we seem to be losing
information. Indeed, there can be many different $Y$'s with the same
Jacobian $X$, and  there is no obvious way to recover $Y$ from
spectral data consisting of a divisor on $X$ and an ordinary sheaf
supported on this divisor.

The resolution of this problem is suggested by string theory. Namely,
in the transition from $Y$ to $X$ one should add an extra piece of
data corresponding to the physicists' $B$-field. Mathematically the
holomorphic version of this data
is encoded in an ${\mathcal O}^{\times}$-gerbe on $X$. A detailed
discussion of ${\mathcal O}^{\times}$-gerbes and their geometric
properties can be found in
\cite{giraud,brylinski,breen-2g,hitchin-slag} and in our
section~\ref{subsec-Brauer}. A baby version of our result is that:

\begin{itemize}
\item $Y$ determines a non-trivial ${\mathcal
O}^{\times}_{X}$-gerbe $\leftidx{_{Y}}{X}{}$ on $X$. \\
\item There is a gerby Fourier-Mukai transform exchanging
\[
\xymatrix@R-15pt{
{}\save[]+<0cm,0.7cm>*\txt{$\left\{
\begin{minipage}[c]{5in}
{
{\em Bundle data}: vector bundles on $Y$, semistable of degree zero on
the generic fiber of $Y \to B$.
}
\end{minipage}
\right\}
$}
\ar@{<->}[dd]^-{\FM} \restore  \\
\\
{}\save[]-<0cm,1cm>*\txt{$\left\{
\begin{minipage}[c]{5in}
{ {\em Spectral data:} sheaves on the gerbe $\leftidx{_{Y}}{X}{}$ with the
numerics of a line bundle on a `spectral' divisor
$(\leftidx{_{Y}}{X}{})\times_{X} C$, with $C \subset X$ finite over
$B$.  }
\end{minipage}
\right\}
$} \restore
}
\]
\end{itemize}
\

This statement is asymmetric - we only consider vector bundles on the
variety $Y$, and the spectral data appears only for the gerbe
$\leftidx{_{Y}}{X}{}$. The symmetry can be restored by extending this
gerby spectral construction to a full Fourier-Mukai equivalence of
derived categories.

One peculiarity of ${\mathcal O}^{\times}$-gerbes is that the
categories of coherent sheaves on them, and hence also the derived
categories, admit an orthogonal decomposition by subcategories labeled
by the characters of ${\mathcal O}^{\times}$, i.e. by the integers
(see section~\ref{subsec-Brauer}). For any $k \in {\mathbb Z}$, we
will write $D^{b}_{k}(\leftidx{_{Y}}{X}{}) \subset D^{b}(_{Y}X)$ for
the $k$-th summand and we will call $D^{b}_{k}(\leftidx{_{Y}}{X}{})$
the derived category of weight $k$ sheaves on
$\leftidx{_{Y}}{X}{}$. It may be helpful to note here that when
$\leftidx{_{Y}}{X}{}$ comes from an Azumaya algebra, the weight $k$
corresponds to the central character of the action of this algebra.

Related partial dualities were considered previously in
\cite{caldararu-cy,caldararu-k3} in the context of Fourier-Mukai
transforms and in \cite{donagi-gaitsgory} in the context of the
spectral construction. A detailed analysis of the corresponding moduli
spaces and the duality transformation between them was recently
carried out, for the particular case of Hopf-like surfaces, in
\cite{ruxandra1,ruxandra2}.

\newpage

\insubsection{Main results: duality for genus one fibrations}

\noindent
For any elliptic fibration $\xymatrix@1{X \ar[r]^-{\pi} & B
\ar@/^0.5pc/[l]^-{\sigma}}$, the twisted versions of $X$ are
parameterized by the analytic Tate-Shafarevich group $\TSh_{an}(X)$
(see section~\ref{subsec-TSh-generalities}). For a given $\beta \in
\TSh_{an}(X)$, let $\pi_{\beta} : X_{\beta} \to B$ denote the
corresponding genus one fibration. On the other hand for any analytic
space $X$, the analytic ${\mathcal O}_{X}^{\times}$-gerbes on $X$ are
parameterized by the analytic Brauer group $Br'_{an}(X) =
H^{2}_{an}({\mathcal O}_{X}^{\times})$. For any $\alpha \in
Br'_{an}(X)$ we denote the corresponding gerbe by
$\leftidx{_{\alpha}}{X}{}$
and the bounded derived category of
coherent sheaves on $_{\alpha}X$ by $D^{b}(\leftidx{_{\alpha}}{X}{})$.
The latter decomposes naturally as the orthogonal direct sum of
"pure weight" subcategories
$D^{b}_{k}(\leftidx{_{Y}}{X}{})$ indexed by characters of
${\mathcal O}_{X}^{\times}$, i.e. by integers $k \in {\mathbb Z}$.
For all $\alpha \in Br'_{an}(X)$ and
all $k \in {\mathbb Z}$ there is a canonical equivalence
$D^{b}_{k}(\leftidx{_\alpha}{X}{}) \cong
D^{b}_{1}(\leftidx{_{k\alpha}}{X}{})$. This follows immediately by
comparing representations of $_{\alpha}X$ of pure weight $k$ with
representations of $_{k\alpha}X$ of pure weight $1$ or by the
appropriate idenitfications with the category of twisted sheaves on $X$ (see
section~\ref{ssubsec-gerbes}).

Consider first the case when $X$ is a surface. We will see in
section~\ref{subsec-T}  that to any $\alpha$, $\beta \in \TSh_{an}(X)$
we can associate an ${\mathcal O}^{\times}$-gerbe $\Xc{\alpha}{\beta}$
over $X_{\beta}$. The notation $\Xc{\alpha}{\beta}$ generalizes our
previous usage: when $\alpha = 0$ we get the trivial gerbe
$\Xc{0}{\beta}$ on $X_{\beta}$, and when $\beta = 0$ we get the gerbe
$\leftidx{_{\alpha}}{X}{}$ on $X = X_{0}$.
Our main result in the case of surfaces, proved in section~\ref{sec-K3}, is:

\begin{thm} \label{Main-K3} Let
\[
\xymatrix@1{X
\ar[r]^-{\pi} & B \ar@/^0.5pc/[l]^-{\sigma}} ={\mathbb P}^{1}
\]
be a non-isotrivial elliptic fibration on a smooth complex surface $X$. Assume
that $\pi$ has $I_{1}$ fibers at worst.  Let $\alpha, \beta \in
\TSh_{an}(X)$ be two elements such that $\beta$ is torsion. Then there
is an equivalence
\[
\FM : D^{b}_{1}(\Xc{\alpha}{\beta}) \to D^{b}_{-1}(\Xc{\beta}{\alpha})
\]
of the derived category of weight $1$ coherent sheaves on the gerbe
$\Xc{\alpha}{\beta}$ and the derived category of weight $(-1)$ coherent
sheaves on the gerbe $\Xc{\beta}{\alpha}$. Equivalently $\FM$ can be
thought of as an equivalence of the derived categories
$D^{b}_{1}(\Xc{\alpha}{\beta})$ and $D^{b}_{1}(\Xc{-\beta}{\alpha})$.
\end{thm}

\

\noindent
We will see in
section~\ref{subsec-T} that for a higher dimensional $X$, the gerbes
$\Xc{\alpha}{\beta}$ may not make sense for arbitrary choices of
$\alpha, \beta \in \TSh_{an}(X)$. However in section~\ref{subsec-T} we
show that there exists a natural pairing $\langle \bullet,
\bullet\rangle : \TSh_{an}(X)\otimes_{\mathbb Z} \TSh_{an}(X) \to
H^{3}_{an}({\mathcal O}_{B}^{\times})$ and that $\Xc{\alpha}{\beta}$
can be defined whenever $\alpha$, $\beta$ are {\it complementary}, i.e.
$\langle \alpha, \beta \rangle = 0$.
The natural generalization of Theorem~\ref{Main-K3} is the following
(see Conjecture~\ref{con-duality}):

\

\medskip

\noindent
{\bf Main Conjecture} {\em For any complementary pair $\alpha,
\beta \in \TSh_{an}(X)$, there exists an equivalence
\[
D^{b}_{1}(\leftidx{_{\alpha}}{X}{}_{\beta}) \cong D^{b}_{-1}(_{\beta}X_{\alpha})
\]
of the bounded derived categories of sheaves of weights $\pm 1$ on
$\Xc{\alpha}{\beta}$ and $\Xc{\beta}{\alpha}$ respectively.
}

\

\medskip

We cannot prove this in full generality, mostly due to our inability to
handle the general singular fibers. We are able to settle the conjecture
in the non-singular case, under the somewhat more restrictive condition
that $\alpha$ is $m$-divisible and $\beta$ is
$m$-torsion for some integer $m$. This of course implies that
${\alpha},{\beta}$ are complementary, so  $\Xc{\alpha}{\beta}$ is well
defined. For a smooth $\pi$ our main result is:

\begin{thm} \label{Main-smooth} Let

\[
\xymatrix@1{X
\ar[r]^-{\pi} & B \ar@/^0.5pc/[l]^-{\sigma}}
\]
be a smooth elliptic fibration on an algebraic variety $X$ over a
smooth algebraic base $B$. Assume that $Br'_{an}(B) = 0$. Fix a
positive integer $m$ and let $\alpha, \beta \in \TSh_{an}(X)$ be two
elements such that $\alpha$ is $m$-divisible and $\beta$ is
$m$-torsion. Then there is an equivalence
\[
\FM : D^{b}_{1}(\Xc{\alpha}{\beta}) \to D^{b}_{-1}(\Xc{\beta}{\alpha})
\]
of the derived category of weight $1$ coherent sheaves on the gerbe
$\Xc{\alpha}{\beta}$ with the derived category of weight $(-1)$
coherent sheaves on the gerbe $\Xc{\beta}{\alpha}$. Equivalently $\FM$
can be thought of as an equivalence of $D^{b}_{1}(\Xc{\alpha}{\beta})$
and $D^{b}_{1}(\Xc{-\beta}{\alpha})$.
\end{thm}

\

\medskip

\noindent
In fact, the proof of Theorem B is quite a bit easier than that of
Theorem A, so we give it first, in section~\ref{sec-smooth}.
The proof is based on the construction of two explicit presentations:
the lifting presentation of $\Xc{\alpha}{\beta}$, in section
\ref{ssubsec-lifting}, and the extension presentation of
$\Xc{\beta}{\alpha}$, in section \ref{ssubsec-extension}, together with the
construction of an explicit Fourier-Mukai duality between them, in
section \ref{sec-duality}.

Our two main theorems and their proofs have fairly straightforward
analogues asserting the equivalence of the derived categories of
quasi-coherent sheaves and, in the algebraic case, asserting the
equivalence of appropriate categories of algebraic coherent sheaves.
Indeed, if the class $\alpha$ happens to be torsion as well, then the
spaces $X_{\alpha}$ and $X_{\beta}$ are algebraic and the gerbes
$\Xc{\alpha}{\beta}$ and $\Xc{\beta}{\alpha}$ are algebraic stacks in
the sense of Artin. We will see in the proofs of the two theorems that
the gerby Fourier-Mukai transform in this case will correspond to a
kernel object which is algebraic and so will give rise to an
equivalence of the derived categories of weight one algebraic coherent
sheaves.

\insubsection{Duality for commutative group stacks}

\noindent
As was pointed out by Arinkin, our Theorem~\ref{Main-smooth} (but
not Theorem~\ref{Main-K3}) fits very
naturally in the context of commutative group stacks (cgs). The
${\mathcal O}_X^\times$-gerbe $_{\alpha}X$ is a family of cgs over $B$
which is an
extension:
\[
0 \to B\GM \to \Xc{\alpha}{0} \to X \to 0 
\]
of $X$ by the classifying stack of $\GM$. The torsor $X_{\beta}$, on the
other hand, is not a cgs over $B$; but it does determine one, namely the
extension:
\[
0 \to X \to \tXc{}{\beta} \to {\mathbb Z} \to 0
\]
of ${\mathbb Z}$ by $X$, where $X_{\beta}$ is recovered as the inverse
image of $1 \in {\mathbb Z}$. Similarly, the gerbe
$\Xc{\alpha}{\beta}$ which we construct, using either the lifting
presentation (when $\beta$ is $m$-torsion) or the extension presentation
(when $\alpha$ is $m$-torsion), determines a cgs
${\mycal X} = \tXc{\alpha}{\beta}$
which has a two-step filtration, with sub $B\GM$, middle subquotient $X$, and
quotient ${\mathbb Z}$. This can be considered either as an ${\mathcal
O}_X^\times$-gerbe over $X_{\beta}$, or dually as a torsor over
$_{\alpha}X$. In particular, the derived category
$D^b(\Xc{\alpha}{\beta})$ is graded by ${\mathbb Z}\times{\mathbb Z}$.

Now quite generally, such a cgs ${\mycal X}$ has a dual cgs ${\mycal
X}^\vee$ which has a similar two-step filtration with the roles of the
sub and the quotient interchanged. There is a Poincare sheaf $\mycal P$
which is a biextension of ${\mycal X}^\vee \times {\mycal X}$ by $\GM$,
and it induces a Fourier-Mukai transform which is an equivalence of
categories $D^b({\mycal X}^\vee) \simeq D^b({\mycal X})$ interchanging the
two ${\mathbb Z}$ gradings. Our Theorem~\ref{Main-smooth} can
therefore be interpreted as saying that the cgs dual to
$\tXc{\alpha}{\beta}$ is $\tXc{\beta}{\alpha}$; the previous version
is recovered by restricting the equivalence to the piece of bidegree
$(1,-1)$. This is explained in some more detail in the
Appendix~\ref{appendix-motives} which D. Arinkin kindly wrote for us.

This duality picture has a straighforward extension to more general
cgs ${\mycal X}$ over $B$: these are again endowed with a two step
filtration $W_{-2}{\mycal X} \subset W_{-1}{\mycal X} \subset
W_{0}{\mycal X} = {\mycal X}$, where $\op{gr}_{-2} =
W_{-2}{\mycal X} = BT$ for some affine torus bundle $T \to B$, the
middle subquotient $\op{gr}_{-1}$ is some abelian scheme $A \to B$ and
the last quotient $\op{gr}_{0}$ is some bundle $\Lambda \to B$ of
finite rank free abelian groups over $B$. The dual cgs ${\mycal
X}^{\vee}$ is again the stack of homomorphisms
$\underline{\op{Hom}}_{\op{cgs}}({\mycal X},B\GM)$ and one expects
that in good cases the duality gives rise to an equivalence of the
appropriate categories of representations. 

Arinkin calls the cgs described in the previous paragraph `Beilinson's
one motives' since they were considered by Beilinson (unpublished) in
the context of the theory of mixed motives. The cgs ${\mycal X}$ are
formally very similar to the classical one motives studied by Deligne
in \cite{deligne-hodge3}. The one motives of \cite{deligne-hodge3} can
be viewed either as certain mixed Hodge structures of level $\leq 1$
or as cgs ${\mycal M}$ defined over ${\mathbb C}$. As a commutative
group stack, every Deligne's one motive ${\mycal M}$ is equipped with
a two step filtration $W_{-2}{\mycal M} \subset W_{-1}{\mycal M}
\subset W_{0}{\mycal M} = {\mycal M}$, for which 
$\op{gr}_{-2} = T$ for some affine torus $T$, $\op{gr}_{-1} = A$ is a
polarized abelian variety, and $\op{gr}_{0} = B\Lambda$ for some free
abelian group $\Lambda$ of finite rank. If we now look at families of
Deligne's one motives defined over some base $B$ we arrive at cgs over
$B$ which are of essentially the same shape as the Beilinson's one
motives, but with the stackiness appearing at a different subquotient
of the filtration. Furthermore, as explained in \cite{deligne-hodge3},
the dual of a Deligne's one motive ${\mycal M}$ is the cgs ${\mycal
M}^{\vee} := \underline{\op{Hom}}_{\op{cgs}}({\mycal M},B\GM)$, which is
again a one motive of the same type with $\op{gr}_{-2}{\mycal
M}^{\vee} = \op{Hom}(\Lambda,\GM)$, $\op{gr}_{-1}{\mycal M}^{\vee} =
\widehat{A}$ (the dual abelian variety to $A$), and
$\op{gr}_{0}{\mycal M}^{\vee} = B\op{Hom}(T,\GM)$.

In fact, we can view $\underline{\op{Hom}}_{\op{cgs}}(\bullet,B\GM)$
as a transformation acting on commutative group stacks,
which preserves the two natural families of Deligne's and Beilinson's
one motives and induces a duality on each of these families. Moreover,
since in both cases the duality is realized in terms of suitable
biextensions of ${\mycal X}\times {\mycal X}^{\vee}$ and ${\mycal
M}\times {\mycal M}^{\vee}$, one expects that the duality of cgs will
give rise to an equivalence of the corresponding categories of
representations of cgs. For the specific stacks
$\tXc{\alpha}{\beta}$ this is precisely the content of our
Theorem~\ref{Main-smooth}.

\insubsection{The non-commutative aspect}

\noindent
Results having the general shape of Theorem~\ref{Main-K3} were
anticipated in the physics literature. In fact, Ganor-Mihailov-Saulina
have conjectured in \cite{gms} that when $Y$ is a genus one
fibered $K3$ surface, there should exist a non-commutative deformation
$_{Y}X$ of $X = \overline{\op{Pic}^{0}}(Y/B)$ and a categorical
equivalence between instantons on $_{Y}X$ and spectral data on
$Y$. This is a special case of Theorem~\ref{Main-K3}.

This statement admits an intriguing interpretation in terms of
non-commutative geometry, a topic currently of high interest to
physicists \cite{ns,kko}.  According to the general yoga of
deformation quantization (see e.g. \cite{kontsevich}), any symplectic
(or Poisson) structure on $X$ is the first term in a non-commutative
deformation of its structure sheaf. In a suitable algebro-geometric
context, e.g. on a $K3$, a symplectic structure $\theta$ has three
incarnations: as a real $2$-form $\theta_{\mathbb R}$, a holomorphic
$2$-form $\theta^{2,0}$, or an antiholomorphic $2$-form
$\theta^{0,2}$.  Then $\theta_{\mathbb R}$ determines a
"non-commutative four-manifold", and $\theta^{2,0}$ determines a
"non-commutative $K3$". The third incarnation, $\theta^{0,2}$, gives
both the element $X_{\theta}$ in the Tate-Shafarevich group $\TSh(X)$
and the ${\mathcal O}_X^\times$-gerbe $_{\theta}X$. In this
sense, Theorem~\ref{Main-K3} can be viewed as an affirmative answer
and a generalization of the \cite{gms} conjecture.

\insubsection{Modified T-duality and the SYZ conjecture}

\noindent
The celebrated work of Strominger, Yau and Zaslow \cite{SYZ}
interprets mirror symmetry of Calabi-Yau spaces in terms of special
Lagrangian (SLAG) torus fibrations. If a CY manifold $X$ (with ``large
complex struture") has mirror $X'$, \cite{SYZ} conjecture the
existence of fibrations $\pi: X \to B$ and $\pi': X' \to B$ whose
generic fibers are SLAG tori dual to each other: each parameterizes
$U(1)$ flat connections on the other.  In particular, each of these
fibrations admits a SLAG zero-section, corresponding to the trivial
connection on the dual fibers. The analogy with the
theorem of \cite{bridgeland-maciocia}
is clear: the SLAG torus fibration on the Calabi-Yau
threefold replaces the elliptic fibration on the surface, and mirror
symmetry (interchanging D-branes of type B with D-branes of type A)
replaces the Fourier-Mukai transform (which interchanges vector bundles
with spectral data).

Our work suggests that the SYZ conjecture should be extended to a SLAG
analogue of Theorem~\ref{Main-K3} or of the Main Conjecture, in which
the physical B-fields $\alpha \in H^2(X, {\mathbb R}/{\mathbb Z})$
play the role of our gerbes. This extension leads to an integrable
system structure on the moduli space underlying mirror symmetry. We
give an informal discussion of these matters in section \ref{sec-SYZ},
and we hope to return to them in future work.

\insubsection{Modularity}

\noindent
As often happens in physics, the
Fourier-Mukai functor $\FM : D^{b}_{1}(\Xc{\alpha}{\beta}) \to
D^{b}_{1}(\Xc{-\beta}{\alpha})$ is just one particular element of a
whole collection of dualities. For simplicity consider only the case
of a projective elliptic surface $\pi : X \to B$. In this case,
our Fourier-Mukai
duality works for any pair of elements $(\alpha,\beta) \in
\TSh(X)\times \TSh(X)$ in the algebraic Tate-Shafarevich group. Thus
the Fourier-Mukai functor corresponds to the action of the matrix
\[
\begin{pmatrix}
0 & 1 \\ -1 & 0
\end{pmatrix} \in \op{SL}(2,{\mathbb Z})
\]
on the Cartesian square $\TSh(X)^{\times 2}$ of the abelian group
$\TSh(X)$. Moreover, one can show (see e.g. section~\ref{subsec-T})
that for surfaces the natural map $T_{\beta} : \TSh(X) \to
Br'(X_{\beta})$, used to define our gerbes, has kernel generated by
the element $\beta \in \TSh(X)$. In particular, $T_{\beta}(\alpha +
\beta) = T_{\beta}(\alpha)$ and so the gerbes $\Xc{\alpha +
\beta}{\beta}$ and $\Xc{\alpha}{\beta}$ are isomorphic. A choice of
such an isomorphism gives rise to an equivalence
$D^{b}_{1}(\Xc{\alpha +\beta}{\beta}) \cong
D^{b}_{1}(\Xc{\alpha}{\beta})$ which corresponds to the action of the
matrix
\[
\begin{pmatrix}
1 & 1 \\ 0 & 1
\end{pmatrix} \in \op{SL}(2,{\mathbb Z})
\]
on $\TSh(X)^{\times 2}$. Since these two matrices generate
$\op{SL}(2,{\mathbb Z})$, it will be very interesting to investigate
which braid group extension of $\op{SL}(2,{\mathbb Z})$ acts on
$\coprod_{\alpha, \beta \in \TSh(X)} D^{b}_{1}(\Xc{\alpha}{\beta})$,
lifting the action on $\TSh(X)^{\times 2}$. In the case when the
Mordell-Weil group of $X$ is trivial we expect this extension to be
central and to be related to the extensions appearing in
\cite{polishchuk-biextensions,polishchuk-weil}, \cite{orlov-abelian}
and \cite{seidel-thomas}. We do not discuss this question here but
hope to return to it in a future work.

\insubsection{Twisted sheaves}

\noindent
Another context in which Theorem~\ref{Main-K3} turns out to be
relevant is the theory of twisted sheaves on a complex space which
admits a genus one fibration. Recall \cite{caldararu-thesis} that for
any ${\mathcal O}^{\times}$-valued \v{C}ech 2-cocycle $\alpha$ on a
complex space $X$, one can consider the abelian category of
$\alpha$-twisted sheaves on $X$ and its derived category
$D^{b}(X,\alpha)$. By definition, an $\alpha$-twisted sheaf on $X$ is
a collection of coherent sheaves defined over open sets in $X$,
together with gluing data on overlaps which satisfy the
$\alpha$-twisted cocycle condition on triple overlaps (see
\cite{caldararu-thesis} or our section~\ref{ssubsec-gerbes} for
details). Refininig the open covering or changing the cocycle by a
coboundary results in an equivalent category of twisted
sheaves. Twisted sheaves on Calabi-Yau manifolds, and in particular on
genus one fibered Calabi-Yau manifolds, were recently studied by
A.C\u{a}ld\u{a}raru \cite{caldararu-thesis,caldararu-cy}. In particular, he
observed \cite{caldararu-k3} that in the case of a K3 surface, the
derived category of $\alpha$-twisted sheaves possesses certain natural
Fourier-Mukai partners.  The starting point of his analysis is the
observation that if $X$ is a smooth projective K3 surface, then every
element $\alpha \in H^{2}_{\et}(X,{\mathcal O}^{\times})$ can be
interpreted as a homomorphism $\alpha : \Tr_{X} \to {\mathbb
Q}/{\mathbb Z}$, where $\Tr_{X}$ denotes the transcendental lattice of
$X$ (see \cite{caldararu-k3} or section~\ref{ssubsec-gerbes}). This
interpretation suggests the following:
\

\bigskip

\noindent
{\bf C\u{a}ld\u{a}raru's Conjecture} \  Let $X$ and $Y$ be two projective K3
surfaces and let $\alpha \in H^{2}_{\et}(X,{\mathcal O}^{\times})$ and
$\beta \in H^{2}_{\et}(Y,{\mathcal O}^{\times})$. Then the derived
categories $D^{b}(X,\alpha)$ and $D^{b}(Y,\beta)$ are equivalent as
triangulated categories iff the lattices
$\ker(\alpha) \subset \Tr_{X}$ and $\ker(\beta) \subset \Tr_{Y}$ are
Hodge isometric.

\

\bigskip

\noindent
When both
$\alpha$ and $\beta$ are zero, the conjecture is true in view of a
theorem of D.Orlov \cite{orlov-k3} asserting that two smooth
projective K3 surfaces have equivalent derived categories iff their
transcendental lattices are Hodge isometric. This has been extended by
C\u{a}ld\u{a}raru, who used Mukai's quasi-universal sheaves for
non-fine moduli spaces of sheaves on K3 surfaces to deduce that the
conjecture holds whenever one of the classes, say $\beta$, is trivial.
The algebraic case of our Theorem A proves C\u{a}ld\u{a}raru's
Conjecture in a series of new cases, with both $\alpha$ and $\beta$ non-
zero.

\

\bigskip

\noindent
Indeed, if $\pi : X \to B$ is an elliptic K3 surface and if
$\alpha, \beta \in \TSh(X)$ are two elements in the algebraic
Tate-Shafarevich group, then the natural identification $\TSh(X) =
H^{2}_{\et}(X,{\mathcal O}^{\times})$ coming from the Leray spectral
sequence allows us to view both $\alpha$ and $\beta$ as homomorphisms
$\Tr_{X} \to {\mathbb Q}/{\mathbb Z}$. Using this interpretation one
checks immediately that the transcendental lattices of the K3 surfaces
$X$, $X_{\alpha}$ and $X_{\beta}$ satisfy $\Tr_{X_{\alpha}} =
\ker(\alpha) \subset \Tr_{X}$ and $\Tr_{X_{\beta}} = \ker(\beta)
\subset \Tr_{X}$, where it is understood that all the equalities are
Hodge isometries. Let $T_{\beta}(\alpha) \in Br'_{an}(X_{\beta}) =
H^{2}_{an}(X_{\beta},{\mathcal O}^{\times})$ denote the class of the
gerbe $\Xc{\alpha}{\beta}$.  Assuming that $\alpha$ and $\beta$ are in general
position in $H^{2}_{\et}(X,{\mathcal O}^{\times})$, i.e. that the cyclic
subgroups generated by $\alpha$ and $\beta$ intersect only at zero,
we have natural identifications of Hodge lattices:
\[
\begin{split}
\ker\left[\Tr_{X_{\alpha}} \xrightarrow{T_{\alpha}(\beta)} {\mathbb
Q}/{\mathbb  Z}\right] & = \ker(\alpha)\cap \ker(\beta) \subset
\Tr_{X} \\
\ker\left[\Tr_{X_{\beta}} \xrightarrow{T_{\beta}(\alpha)} {\mathbb
Q}/{\mathbb  Z}\right] & = \ker(\alpha)\cap \ker(\beta) \subset
\Tr_{X}.
\end{split}
\]
In other words - the hypothesis of Theorem~\ref{Main-K3} implies the
hypothesis of C\u{a}ld\u{a}raru's conjecture. Combined with the remark
that $D^{b}(X_{\alpha},T_{\alpha}(\beta)) \cong
D^{b}_{1}(\Xc{\alpha}{\beta})$ and $D^{b}(X_{\beta},-T_{\beta}(\alpha))
\cong D^{b}_{-1}(\Xc{\beta}{\alpha})$, this shows that
Theorem~\ref{Main-K3} implies an interesting new case of
C\u{a}ld\u{a}raru's conjecture (see Corollary~\ref{cor-andrei} for a
slightly more general statement). Note that the condition (required in the
statement of Theorem~\ref{Main-K3}) that $\beta$ should be torsion, is
vacuous in this case, since for a smooth projective surface $X$, both
the cohomological Brauer group $H^{2}_{\et}(X,{\mathcal O}^{\times})$
and the Tate-Shafarevich group $\TSh(X)$ are torsion groups.

\

\bigskip

The paper is organized as follows. In Section~\ref{sec-TSh} we recall
some standard facts about the geometry of ${\mathcal O}^{\times}$
gerbes and genus one fibrations. We also derive the compatibility
condition between two Tate-Shafarevich classes and state a general
conjecture on the equivalence of derived categories for gerbes over
genus one fibrations. In Section~\ref{sec-smooth} we introduce the
main characters appearing in the proofs of the two theorems stated
above. Working in the setup of Theorem~\ref{Main-smooth}, we define
two geometric presentations - the lifting and extension presentations
- for the gerbes $\Xc{\alpha}{\beta}$ and
$\Xc{\beta}{\alpha}$. Furthermore, we construct an integral transform
between the corresponding atlases. We show that this integral
transform sends descent data to descent data and gives rise to an
equivalence of the derived categories of the gerbes, thus proving
Theorem~\ref{Main-smooth}. Section~\ref{sec-K3} deals with the case of
surfaces. We show how, in the case of a surface, one can extend the
lifting and extension presentations across the singular fibers and
produce a Fourier-Mukai transform between the corresponding gerbes. We
again check that this transform is an equivalence, which proves
Theorem~\ref{Main-K3}. Finally, in Section~\ref{sec-SYZ} we discuss the analogy
between algebraic gerbes over  genus one fibrations and flat gerbes
over SLAG 3-torus fibrations on Calabi-Yau threefolds.  We describe a
conjectural picture which amends the Strominger-Yau-Zaslow version of
mirror symmetry to incorporate non-trivial B-fields on both sides of
the mirror correspondence.

\

\bigskip

\noindent
{\bf Acknowledgments:} We would like to thank D.Arinkin,
A.C\u{a}ld\u{a}raru, M.Gross, A.Kresch and D.Orlov for insightful
discussions on the subject of this paper. We would also like to thank
KITP and MSRI for providing a stimulating research environment during
certain stages of the preparation of this work.

\section{The Brauer group and the Tate-Shafarevich group} \label{sec-TSh}

We need some basic facts relating elements of the
Brauer group to elements of the Tate-Shafarevich group of an elliptic
fibration. We discuss ${\mathcal O}^{\times}$-gerbes  and the Brauer
groups which classify them in section \ref{subsec-Brauer}, then genus-1
fibrations and the Tate-Shafarevich group which classifies them, in
section \ref{subsec-TSh-generalities}. For an elliptic fibration there
is a simple, direct relation between these two groups. The extension to
genus-1 fibrations though is more delicate, and is defined only when a
certain alternating pairing vanishes. This is discussed in section
\ref{subsec-T}.

\subsection{Brauer groups and ${\mathcal O}^{\times}$-gerbes}
\label{subsec-Brauer}

In this section we review the notions of ${\mathcal O}^{\times}$-gerbe
and presentation, and discuss the relationship between ${\mathcal
O}^{\times}$-gerbes and  elements in the Brauer group.

\subsubsection{$\mycal{H}$-gerbes} \label{ssubsec-gerbes}

Let $\mycal{H}$ be a sheaf of abelian groups on
a topological space (or a site) $X$.  The case of main interest for us
is when $(X,{\mathcal O}_{X})$ is a ringed space and $\mycal{H} =
{\mathcal O}^{\times}_{X}$ is the sheaf of invertible elements in the
structure sheaf. In fact, most of the time we will have $\mycal{H} =
{\mathcal O}_{X}^{\times}$ in either the etale or the analytic
topologies on a complex scheme (or an algebraic or analytic space)
$X$. In section \ref{sec-SYZ} we will be interested also in the case
when $\mycal{H}$ is the sheaf of germs of smooth maps from a
$C^{\infty}$ manifold $X$ to the circle $S^{1}$.

An $\mycal{H}$-gerbe on $X$ is a global structure on $X$ which
``locally looks like the quotient of $X$ by the trivial action of
$\mycal{H}$''. More precisely ``the quotient of $X$ by the trivial
action of $\mycal{H}$'' is the classifying object $B\mycal{H}$. For 
example, in case $\mycal{H}$ is the sheaf of holomorphic maps
from $X$ to a fixed group $H$, $B\mycal{H}$ is the sheaf of
sections of $X \times BH$ over $X$, where $BH$ is the classifying
space of $H$. In the general case, $B\mycal{H}$
can be interpreted either as a topological space over $X$ (defined up
to homotopy), or as a stack in groupoids over $X$ (see \cite[\S
3]{laumon-stacks} for the definition). We adopt the second approach
and treat $B\mycal{H}$ as a stack (=`sheaf of categories'): over any
open set $V$, the objects of $B\mycal{H}(V)$ are the
$\mycal{H}$-torsors on $V$ and the morphisms are the isomorphisms of
torsors. In particular, the automorphisms of the trivial torsor
$\boldsymbol{1}_{V}$ are given by elements in $\mycal{H}(V)$. Note
that $B\mycal{H}$ is in fact a commutative group stack over $X$ with a
group structure given by convolution of
$\mycal{H}$-torsors. Explicitly, for any two $\mycal{H}$-torsors $A'$
and $A''$ over $V$ the convolution $A'\otimes A''$ is defined as the
$\mycal{H}$-torsor $(A'\times A'')/\ker(\boldsymbol{m}_{\mycal{H}})$,
where $\boldsymbol{m}_{\mycal{H}} : \mycal{H}\times \mycal{H} \to
\mycal{H}$ is the multiplication map.

\begin{defi} \label{def-H-gerbe} An $\mycal{H}$-gerbe on $X$ is a
$B\mycal{H}$ torsor, i.e. a stack of groupoids ${_{\alpha}{X}}$ over $X$,
which is equipped with a principal homogeneous action of $B\mycal{H}$.
\end{defi}

\begin{rem} \label{rem-banded} $\bullet$ Explicitly, a stack of
groupoids $\leftidx{_{\alpha}}{X}{} \to X$ is an $\mycal{H}$-gerbe if
for any open $V \subset X$ and any object $s$ of $\leftidx{_{\alpha}}{X}{}(V)$
we have chosen isomorphisms $\mycal{H}(V) \cong
\op{Aut}_{\leftidx{_{\alpha}}{X}{}(V)}(s)$, compatible with pullbacks.

\medskip

\noindent
$\bullet$ In the literature \cite{giraud}, \cite{breen-bit,breen-2g},
one encounters a more general notion of an $\mycal{H}$-gerbe, namely -
a stack $\mycal{T}$ of groupoids on $X$, which is locally isomorphic
to $B\mycal{H}$.  These more general gerbes  are classified by the
first cohomology of $X$ with coefficients in the $1$-truncated
simplicial abelian group $\mycal{H} \to \op{Aut}(\mycal{H})$
\cite{breen-bit}.  They are intimately related to
the forms of $\mycal{H}$, i.e. to sheaves of groups on $X$ which are
only locally isomorphic to $\mycal{H}$. This relatishionship is based
on the identification $\op{Out}(\mycal{H}) =
\underline{\op{Aut}}_{X}(B\mycal{H})$: to any $\mycal{T} \to X$, which is an
$\mycal{H}$-gerbe in this more general sense, one naturally associates an
$\op{Out}(\mycal{H})$-torsor $\op{band}(\mycal{T}) :=
\underline{\op{Isom}}_{X}(\mycal{T},B\mycal{H})$ - the band of the
gerbe $\mycal{T}$ \cite{giraud}, \cite{breen-bit}. A gerbe $\mycal{T}$
is said to be banded by $\mycal{H}$ if it is equipped with a
trivialization of the torsor $\op{band}(\mycal{T})$. When $\mycal{H}$
is abelian, this condition is equivalent to requiring that for any
open $V$ and any $s \in \mycal{T}$ we have chosen isomorphisms
$\mycal{H}(V) \cong \op{Aut}_{\mycal{T}}(s)$ in a way compatible with
pullbacks. In other words, the more restrictive notion of an
$\mycal{H}$-gerbe that we have adopted in this paper is the same as
the standard notion of an $\mycal{H}$-banded gerbe (at least for an
abelian $\mycal{H}$). We will casually ignore this distinction and
will call all our gerbes simply $\mycal{H}$-gerbes.
\end{rem}

\

\bigskip

In case $\mycal{H} = {\mathcal O}^{\times}_{X}$ (in the relevant
topology), the classifying stack $B{\mathcal O}_{X}^{\times}$ is the
sheaf of Picard categories $\mycal{P}ic(X)$: for an open $U$, the
objects of $\mycal{P}ic(X)(U)$ are by definition the line bundles on
$U$, and for two objects $L, M \in \op{ob}(\mycal{P}ic(X)(U))$ the set
$\op{Hom}_{\mycal{P}ic(X)}(L,M)$ is defined to be $\op{Isom}(L,M)$. An
${\mathcal O}_{X}^{\times}$ gerbe $\leftidx{_{\alpha}}{X}{}$ assigns to each open
$U$ a $\mycal{P}ic(X)(U)$-torsor, denoted $\mycal{P}ic_{\alpha}(U)$,
with a compatibility of the assignments to different $U$'s. We can
thus think of a section of an ${\mathcal O}_{X}^{\times}$-gerbe as a
twisting of the notion of a line bundle on $X$. More generally
the sections in an $\mycal{H}$-gerbe are twistings of the notion of an
$\mycal{H}$-torsor on $X$: simply replace in the previous discussion
each appearance of $\mycal{P}ic$ with $\mycal{T}ors^{\mycal{H}}$ - the
group of $\mycal{H}$-torsors. This interpretation suggests that the
group classifying $\mycal{H}$-gerbes should be $H^{2}(X,\mycal{H})$.
When $\mycal{H}$ is abelian this statement can be made precise via the
standard cohomological machinery \cite[IV.2]{milne-book} or
\cite{giraud}, \cite{breen-bit} (but keep in mind that our
$\mycal{H}$- gerbes are the $\mycal{H}$-banded gerbes of loc. cit.).

In more down to earth terms the interpretation of the elements in
$H^{2}(X,\mycal{H})$ as equivalence classes of gerbes can be seen as
follows. Assume that we are in the good situation when the cohomology
of $\mycal{H}$ can be computed in \v{C}ech terms.
Let $\{\alpha_{ijk} \}$ be an $\mycal{H}$-valued  \v{C}ech
$2$-cocycle  w.r.t. an open cover $\{ U_{i} \}$ of $X$. An object $L$
of $\mycal{T}ors_{\alpha}^{\mycal{H}}$ is defined to be an assignment
of an $\mycal{H}(U_{i})$-torsor $L(U_{i})$ to each $U_{i}$, together
with transition functions
\[
g_{ij} : L(U_{i})\otimes_{\mycal{H}(U_{i})} \mycal{H}(U_{ij})
\widetilde{\to} L(U_{j})\otimes_{\mycal{H}(U_{j})} \mycal{H}(U_{ij})
\]
satisfying the twisted cocycle condition:
\[
g_{ij}\circ g_{jk} \circ g_{ki} = \alpha_{ijk}
\]
on triple intersections.
A morphism between two $\alpha$-twisted $\mycal{H}$-torsors $L'$
and $L''$ is given by a compatible collection of isomorphisms
$L'(U_{i}) \widetilde{\to} L''(U_{i})$.

Similarly we define the category
$\mycal{T}ors_{\alpha}^{\mycal{H}}(U)$ for any open $U$. The resulting
sheaf of categories (=stack) $\mycal{T}ors_{\alpha}^{\mycal{H}}$ on
$X$ is by definition a torsor over $B\mycal{H} =
\mycal{T}ors_{\boldsymbol{1}}^{\mycal{H}}$, i.e. an $\mycal{H}$-gerbe,
which we denote by $\leftidx{_{\alpha}}{X}{}$. Clearly two cocylces which
represent the same cohomology class in $\check{H}^{2}(X,\mycal{H})$ define
isomorphic gerbes. Conversely, if sheaf cohomology on $X$ can be
computed in \v{C}ech terms, any $\mycal{H}$-gerbe arises this way from
some $\alpha$ w.r.t. a sufficiently refined cover \cite{giraud},
\cite{breen-bit,breen-2g}.

\smallskip

\noindent
{\bf Notation:} $\bullet$ Given an $\mycal{H}$-gerbe $\mycal{T}$
over $X$ we write $[\mycal{T}] \in H^{2}(X,\mycal{H})$ for the
element that classifies it.

\smallskip

\noindent
$\bullet$ The base space $X$ for an $\mycal{H}$-gerbe $\mycal{T} \to
X$ is called the {\em coarse moduli space} of $\mycal{T}$. This
terminology reflects the fact that $X$ represents the sheaf of sets
$\pi_{0}(\mycal{T})$, i.e. the sheaf of isomorphism classes of sections
in $\mycal{T}$.

\

\bigskip

\noindent
{\bf Basic construction:} Starting with an algebraic (or analytic) space
$X$ and a short exact sequence of sheaves of groups
\[
1 \to \mycal{H} \to \mycal{G} \to \mycal{K} \to 1,
\]
with $\mycal{H}$-abelian, we get a coboundary map $\delta :
H^{1}(X,\mycal{K}) \to H^{2}(X,\mycal{H})$. This admits the following
lift on the level of torsors and gerbes: a $\mycal{K}$-torsor
$\mycal{C}$ with class $[\mycal{C}] \in H^{1}(X,\mycal{K})$ determines
an $\mycal{H}$-gerbe $\delta(\mycal{C})$ with class
$\delta([\mycal{C}]) \in H^{2}(X,\mycal{H})$. Explicitly, for an open
$U$, $\delta(\mycal{C})(U)$ is the category of pairs $(\mycal{D},\iota)$
where $\mycal{D}$ is a $\mycal{G}$-torsor on $U$ and $\iota :
\mycal{D}\times_{\mycal{G}} \mycal{K} \to \mycal{C}$ is an isomorphism
of $\mycal{K}$-torsors on $U$.

\

\medskip

A familiar special case involves the
sequence
\[
1 \to {\mathcal O}_{X}^{\times} \to GL_{n}({\mathcal O}_{X}) \to
{\mathbb P}GL_{n}({\mathcal O}_{X}) \to 1.
\]
It says that every projective bundle on $X$ gives rise to an
${\mathcal O}_{X}^{\times}$-gerbe which is trivial if and only if the
projective bundle is a projectivization of a vector bundle.

\medskip

\noindent
If $(X,{\mathcal O}_{X})$ is a nice ringed space for which cohomology
can be computed in \v{C}ech terms, then
the choice of $\alpha \in H^{2}(X,{\mathcal O}_{X}^{\times})$
 gives rise to the notion of
$\alpha$-twisted sheaves on $X$.   More precisely, let ${\mathfrak U} =
\{U_{i}\}$ be an open cover of $X$ (in the topology under
consideration) and let
\[
\unalpha = \{ \alpha_{ijk} \} \in \check{C}^{2}({\mathfrak
U},{\mathcal O}_{X}^{\times})
\]
be a $2$-cocycle representing $\alpha \in H^{2}(X,{\mathcal
O}_{X})$. One defines an $\unalpha$-twisted sheaf on $X$ as a
collection $\{ F_{i} \}$ of sheaves $F_{i} \to U_{i}$ of ${\mathcal
O}_{X}$-modules, together with a collection of gluing isomorphisms
\[
\varphi_{ij} : F_{j|U_{ij}} \stackrel{\cong}{\to} F_{i|U_{ij}}
\]
satisfying $\varphi_{ii} = \op{id}$, $\varphi_{ij} =
\varphi_{ji}^{-1}$, and
$
\varphi_{ij}\circ \varphi_{jk} \circ \varphi_{jk} : F_{i|U_{ijk}}
 \to F_{i|U_{ijk}}
$
is given by multiplication by $\alpha_{ijk}$. Given two $\unalpha$
twisted sheaves $\boldsymbol{F} = \{ F_{i}, \varphi_{ij} \}$ and
$\boldsymbol{G} = \{
G_{i}, \gamma_{ij} \}$ we define a
homomorphism $\boldsymbol{f} : \boldsymbol{F} \to \boldsymbol{G}$ to be a
collection $\boldsymbol{f} = \{ f_{i} \}$ of sheaf morphisms $f_{i}  :
F_{i} \to G_{i}$ satisfying $f_{i}\circ \varphi_{ij} =
\gamma_{ij}\circ f_{j}$. Composition is defined in an obvious way and
so we obtain a category of $\unalpha$-twisted sheaves, which depends
both on the cover ${\mathfrak U}$ and on the cocycle $\unalpha$. It
can be checked \cite[Section~1.2]{caldararu-thesis} that the operations of
passing to a refinement ${\mathfrak U}'$  of ${\mathfrak U}$ and of
replacing $\unalpha$ by a cohomologous cocycle $\unalpha'$, give rise
to an equivalent category of $\unalpha'$-twisted sheaves. Thus for any
$\alpha \in H^{2}(X,{\mathcal O}_{X}^{\times})$ we get
category $({\mathcal O}_{X},\alpha)\op{-mod}$ of $\alpha$-twisted
sheaves on $X$ (defined only up to a non-canonical equivalence).
An $\alpha$-sheaf $\boldsymbol{F}$ on $X$ is called quasi-coherent
(respectively coherent) if each $F_{i}$ is quasi-coherent
(respectively coherent). We will write $\op{QCoh}(X,\alpha)$ and
$\op{Coh}(X,\alpha)$ for the categories of quasi-coherent and coherent
$\alpha$-twisted sheaves. Note that $({\mathcal
O}_{X},\alpha)\op{-mod}$, $\op{QCoh}(X,\alpha)$ and
$\op{Coh}(X,\alpha)$ are all abelian categories.

More intrinsically the $\alpha$-twisted sheaves on $X$ can be
interpreted as weight one sheaves on $\leftidx{_{\alpha}}{X}{}$,
where a sheaf on $\leftidx{_{\alpha}}{X}{}$ is understood as a
representation of the sheaf of groupoids $\leftidx{_{\alpha}}{X}{} \to
X$. To spell what this means, let us denote by $\mycal{Q}\op{Coh}_{X}$
the stack of quasicoherent sheaves on the space $X$. Let $\mycal{X}
\to X$ be any fibered category over $X$.  Recall (see
e.g. \cite[Section~3.3]{deligne-tannaka} or
\cite[Definition~13.3.3]{laumon-stacks}) that a {\em representation}
of $\mycal{X}$ is a morphism $F : \mycal{X} \to \mycal{Q}\op{Coh}_{X}$
of fibered categories defined over $X$. Explicitly this means that for
any algebraic (or analytic) space $T \to X$ we are given a finctor
$F_{T} : \mycal{X}(T) \to \op{QCoh}(T)$ so that $F_{T}$ is compatible
with base changes.

In particular, if $\leftidx{_{\alpha}}{X}{} \to X$ is an ${\mathcal
O}_{X}^{\times}$-gerbe, then a representation of
$\leftidx{_{\alpha}}{X}{}$ is a $X$-functor $F :
\leftidx{_{\alpha}}{X}{} \to \mycal{Q}\op{Coh}_{X}$. Given an integer
$n$ we say that $F$ is a {\em pure
$\leftidx{_{\alpha}}{X}{}$-representation of weight $n$} if for any
open $U \subset X$ and any section $L \in \leftidx{_{\alpha}}{X}{}(U)$
the natural sheaf homomorphism $\underline{\op{Aut}}_{U}(L) \to
\underline{\op{Aut}}_{U}(F(L))$ induced by $F$ factors as
\[
\xymatrix{
\underline{\op{Aut}}_{U}(L) \ar[r] & \underline{\op{Aut}}_{U}(F(L)) \\
{\mathcal  O}_{U}^{\times} \ar@{=}[u] \ar[r]^-{(\bullet)^{n}} &
{\mathcal O}_{U}^{\times} \ar@{^{(}->}[u]
}
\]
where in the bottom row the map is the raising into power $n$. It is
instructive to point out that when we are dealing with the trivial
gerbe $_{0}X$ on $X$, a representation of $_{0}X$ is nothing but a
quasicoherent sheaf $F$ on $X$ equipped with a direct sum
decomposition $F = \oplus_{n \in {\mathbb Z}} F_{n}$ into
quasicoherent sheaves $F_{n}$ so that a locally defined function $f
\in {\mathcal O}_{X}$ acts on $F$ as multiplication by $f^{n}$ on
$F_{n}$.  The reader can check as an exercise that the category of
representations of $\leftidx{_{\alpha}}{X}{}$ of pure weight one is
equivalent to the category $\op{QCoh}(X,\alpha)$ and that the category
of representations of $\leftidx{_{\alpha}}{X}{}$ of pure weight $n$ is
equivalent to the category $\op{QCoh}(X,n\alpha)$.

\subsubsection{Geometric gerbes and their presentations}
\label{ssubsec-presentations}

In this section we recall a more geometric approach to ${\mycal
H}$-gerbes which involves gluing of certain good local models. This
exploits the standard idea that various geometric objects can be
conveniently presented in terms of an atlas modulo certain gluing
relations on it. For example, for a a manifold $X$, an atlas $U$ can
be taken to be the disjoint union $U = \coprod _{i} U_{i}$ of
coordinate charts, and the gluing can be specified by the closed subset
of relations
\[
R := U\times_{X} U \subset U\times U,
\]
which comes together with two maps $\bs, \bt : R \to U$ (corresponding to
the two projections of $U\times U$ onto $U$) each of which is a local
diffeomorphism. Analogously, presentations can be used to define
schemes, algebraic spaces and analytic spaces.

Formally, a {\em presentation by objects in a (fibered) category
$\mycal{C}$ (or a groupoid in $\mycal{C}$)} consists of the following
data:
\begin{align*}
\tag{atlas} \text{ an object $U$ of $\mycal{C}$} \\
\tag{relations} \text{ an object $R$ of $\mycal{C}$} \\
\tag{source-target maps} \xymatrix@1{R\ar@<.5ex>[r]^-{\bs}
\ar@<-.5ex>[r]_-{\bt} & U} \\
\tag{composition map} \xymatrix@1{R\times_{U} R \ar[r]^-{\bm} & R} \\
\tag{inversion map} \xymatrix@1{R \ar[r]^-{\bi} & R} \\
\tag{identity map} \xymatrix@1{U \ar[r]^-{\be} & R}.
\end{align*}
These data are subject to the obvious analogues of the group axioms,
applied to the maps $\bm$, $\bi$ and $\be$.

Note that any morphism $\gamma : U \to X$ in $\mycal{C}$ determines a
presentation $({\mathfrak R},U, {\mathfrak m},{\mathfrak i},
{\mathfrak e})$ in $\mycal{C}$, where: ${\mathfrak R} := U\times_{X}
U$; the maps ${\mathfrak s}$, ${\mathfrak t}$ are the two projections;
the composition map ${\mathfrak m}$ sends $(a,b)\times (b,c)$ to
$(a,c)$, ${\mathfrak i}$ sends $(a,b)$ to $(b,a)$; and ${\mathfrak e}$
is the diagonal map. In this situation we identify $X$ with the
quotient $U/{\mathfrak R}$ and we will say that
$\xymatrix@1{{\mathfrak R}\, \ar@<.5ex>[r]^-{{\mathfrak s}}
\ar@<-.5ex>[r]_-{{\mathfrak t}} & U}$ is generated by $\gamma$.

Let now $\gamma : U \to X$ be a morphism of complex schemes and let
\[
\xymatrix@1{{\mathfrak R}\, \ar@<.5ex>[r]^-{{\mathfrak s}}
\ar@<-.5ex>[r]_-{{\mathfrak t}} & U \ar[r]^-{\gamma} & X}
\]
be the presentation of $X$ generated by $\gamma$. Let ${\mathfrak
p}_{1}, {\mathfrak p}_{2},
{\mathfrak m} : {\mathfrak R}\times_{U} {\mathfrak R} \to {\mathfrak R}$
denote the two projections and the multiplication map respectively.
Let $H$ be an abelian group scheme over $X$, $\mycal{H} \to X$ its
sheaf of sections, ${\mathfrak H}$ its pullback to ${\mathfrak R}$ via
$\gamma\circ {\mathfrak s} = \gamma\circ {\mathfrak t}$, and let $\pi
: R \to {\mathfrak R}$ be an ${\mathfrak H}$-torsor over ${\mathfrak
R}$. In order for
\[
\xymatrix@1{R\, \ar@<.5ex>[r]^-{{\bs}}
\ar@<-.5ex>[r]_-{{\bt}} & U }, \quad \text{with }
\bs := {\mathfrak s}\circ \pi, \bt :=
{\mathfrak t}\circ \pi,
\]
to be  a groupoid  we need a {\it biextension isomorphism}
\[
{\mathfrak p}_{1}^{*}R\otimes {\mathfrak p}_{2}^{*} R \to {\mathfrak m}^{*} R
\]
of torsors over ${\mathfrak R}\times_{U} {\mathfrak R}$, as well as a
lifting
\[
\be : U \to R
\]
of ${\mathfrak e}$, i.e. a trivialization of ${\mathfrak e}^{*}R$ on
$U$. Given these data we obtain a new presentation
\[
(\xymatrix@1{R\, \ar@<.5ex>[r]^-{{\bs}}
\ar@<-.5ex>[r]_-{{\bt}} & U }, \bm,\bi,\be).
\]
We claim that this presentation determines an $\mycal{H}$-gerbe
$[U/R]$ on $X$, which we can interpret as the (stacky) quotient of $U$
by $R$.

Indeed, for any open $V$ in $X$, define $[U/R]'(V)$ to be the category
of pairs $(T,j)$, where $T$ is an $H$-torsor over $U_{|V}$ and
\[
j : {\mathfrak t}^{*}T \stackrel{\cong}{\to} {\mathfrak s}^{*}T\otimes R
\]
is an isomorphism of ${\mathfrak H}$-torsors on ${\mathfrak R}$. As
$V$ varies, this gives a prestack $[U/R]'$ of groupoids over $X$. We
define $[U/R]$ to be  the stackification (see Remark~\ref{rem-stack} below)
of $[U/R]'$. By construction $[U/R]'(V)$ is a torsor over the tensor
category of pairs $(T_{0},j_{0})$ where $T_{0}$ is an $H$-torsor over
$U_{|V}$ and
\[
j_{0} : {\mathfrak t}^{*}T_{0} \stackrel{\cong}{\to} {\mathfrak t}^{*}T_{0}
\]
is an isomorphism of ${\mathfrak H}$-torsors on ${\mathfrak
R}_{|V}$. The stackification of the latter is identified, via descent,
with $B\mycal{H}$ and so $[U/R]$ is indeed an $\mycal{H}$-gerbe on
$X$.

\

\begin{rem} \label{rem-stack} The necessity of taking stackification
in the above construction is dictated by the subtlety of the
conditions required to have a `sheaf of categories'.  Let $\mycal{X}
\to X$ be a category fibered in groupoids. Recall that there are two
types of sheaf-like conditions on can impose on $\mycal{X}$:
\begin{description}
\item[(1)] For any open $V \subset X$  and any two objects $\xi, \eta
\in \mycal{X}(V)$, the presheaf of sets \linebreak
$\left( \text{open } W
\subset V \right)
\mapsto \op{Hom}_{\mycal{X}(W)}(\xi_{|W},\eta_{|W})$
is required to be a sheaf.
\item[(2)] If $V = \cup_{i} W_{i}$ is an open covering of $V$ and  we
have
\begin{itemize}
\item $\xi_{i} \in \op{ob}(\mycal{X}(W_{i}))$;
\item $\varphi_{ij} : \xi_{j|W_{ij}} \widetilde{\to} \xi_{i|W_{ij}}$
isomorphisms satisfying the cocycle condition;
\end{itemize}
then we require the existence of an object $\xi \in
\op{ob}(\mycal{X}(V))$ together with isomorphisms $\psi_{i} :
\xi_{|W_{i}} \widetilde{\to} \xi_{i}$, so that $\varphi_{ij} =
\psi_{i}\circ\psi_{j}^{-1}$.
\end{description}
Now if $\mycal{X}$ satisfies $(1)$ we say that $\mycal{X}$ is a
{\em prestack} over $X$ and if it satisfies $(1)+(2)$, then we say that
$\mycal{X}(X)$ is a {\em stack}.

Given any prestack $\mycal{X}$ one shows (see
\cite[\S~3]{laumon-stacks} for details) that there is a unique (up
to equivalence) stack $\mycal{X}^{a} \to X$ together with a map
$\mycal{X} \to \mycal{X}^{a}$ which is fully faithful and locally on
$X$ is essentially bijective. The stack $\mycal{X}^{a}$ is called the
stackification of $\mycal{X}$ and is completely analogous to the sheaf
one associates with a presheaf of sets.
\end{rem}

\

\begin{rem} \label{rem-presentations}
{\bf (i)} We say that a groupoid $(\xymatrix@1{R\,
\ar@<.5ex>[r]^-{{\bs}} \ar@<-.5ex>[r]_-{{\bt}} & U }, \bm,\bi,\be)$ of
algebraic (or analytic) spaces is {\em smooth} (respectively {\em
etale}) if the structure maps $\bs$ and $\bt$ are smooth
(respectively etale). The stacks $\mycal{X} \to X$ over $X$ which
admit a smooth or etale groupoid presentation (lifting a
presentation for $X$) are the stacks which are closest to schemes and
on which one can do geometry in essentially the same way as on
spaces. In fact if $\mycal{X}$ is a stack 
which admits a smooth (respectively etale) presentation, then
$\mycal{X}$ is called an {\em Artin algebraic stack} (respectively
{\em Deligne-Mumford stack}) and is the main object of study in the
algebraic geometry of stacks \cite{laumon-stacks}.

We consider only stacks for which the diagonal map
$\mycal{X} \to \mycal{X}\times_{X}\mycal{X}$ is affine.
Note that for an $\mycal{H}$-gerbe $\leftidx{_{\alpha}}{X}{} \to X$, the condition
of having an affine diagonal is equivalent to $H \to X$ being an
affine group scheme. In particular $\leftidx{_{\alpha}}{X}{}$ is an algebraic stack
(in the sense of Artin) if and only if $\leftidx{_{\alpha}}{X}{}$ has a groupoid
presentation.

\medskip

\noindent
{\bf (ii)} The case of main interest for us is when $H = {\mathbb
G}_{m}$, so $\mycal{G} = {\mathcal O}_{X}^{\times}$. In this case $R$
is the total space of a (punctured) line bundle on ${\mathfrak R}$ and
$j$, $j_{0}$ are isomorphisms of line bundles.

\medskip

\noindent
{\bf (iii)} The above discussion has an obvious analogue where schemes
are replaced by (algebraic or analytic) spaces or manifolds.

\medskip

\noindent
{\bf (iv)} Not every presentation of $X$ will lift to a presentation
of a given gerbe $\leftidx{_{\alpha}}{X}{}$. For example, only the trivial gerbe
$_{0}X = B\mycal{H} \to X$ can be presented by a lift of the trivial
presentation $\xymatrix@1{X\, \ar@<.5ex>[r]
\ar@<-.5ex>[r] & X \to X}$ generated by $\op{id}_{X} : X \to
X$. However if $\underline{\alpha} = \{\alpha_{ijk} \}$ is an
$\mycal{H}$-valued \v{C}ech cocycle w.r.t. an open covering $\{ U_{i}
\}$ of $X$, then the presentation of $X$ generated by $\gamma: U = \coprod
U_{i} \to X$ can be lifted to a presentation
\[
\xymatrix@1{R\, \ar@<.5ex>[r]^-{{\bs}}
\ar@<-.5ex>[r]_-{{\bt}} & U \ar[r]^-{\gamma} & X}
\]
for an $\mycal{H}$-gerbe $\leftidx{_{\alpha}}{X}{}$, whose classifying element is
the class $\alpha = [\underline{\alpha}] \in H^{2}(X,\mycal{H})$.  To
define this presentation we take $R := \coprod_{i,j} U_{ij}\times
\mycal{H}$ with its natural projections $\bs$ and $\bt$ onto $U =
\coprod_{i} U_{i}$. The multiplication map
$\bm : R\times_{U} R \to R$ sends a point
\[
(x;a,b) \in
U_{ijk}\times_{X} H \times_{X}  H =
(U_{ij}\times_{X} H)\times_{U_{i}} (U_{jk}\times_{X}
H) \subset R\times_{U} R
\]
to the point $(x; \alpha_{ijk}\cdot a\cdot b) \in U_{ik}\times_{X} H$.
The inversion $\bi : R \to R$ sends  $(x,a) \in U_{ij}\times_{X} H$ to
the point $(x,a^{-1})$ and the identity $\be : U \to R$ sends
$x \in U_{i}$ to the point $(x,x;1) \in U_{ii}\times_{X} H$.

Slightly more generally: the same reasoning shows that
for any map of complex spaces $\gamma : U \to
X$ and any $\alpha \in \check{H}^{2}(X,\mycal{H})$ it follows that the
presentation of $X$ generated by $\gamma$ can be lifted to a
presentation for the $\mycal{H}$-gerbe $\leftidx{_{\alpha}}{X}{}$ if and only if
$\gamma^{*}\alpha = 0 \in \check{H}^{2}(U,\gamma^{*}\mycal{H})$.
\end{rem}

\

\bigskip

\noindent
{\bf Basic construction, continued:} Assume we are given a short exact sequence
of sheaves of groups
\[
1 \to \mycal{H} \to \mycal{G} \to \mycal{K} \to 1
\]
on an algebraic (analytic) space $X$. Suppose that $\mycal{H}$ is
commutative and that the sheaves $\mycal{H}$, $\mycal{G}$ and
$\mycal{K}$ are represented by group schemes $H$, $G$ and $K$
respectively. In section~\ref{ssubsec-gerbes} we associated to every
$\mycal{K}$ torsor $\mycal{T} \to X$ an $\mycal{H}$ gerbe
$\delta(\mycal{T})$ with class $\delta([\mycal{T}]) \in
H^{2}(X,\mycal{H})$. In this situation the gerbe $\delta(\mycal{T})$
comes with a natural presentation:
\[
\xymatrix@1{R\, \ar@<.5ex>[r]^-{{\bs}} \ar@<-.5ex>[r]_-{{\bt}} & U
\ar[r] & \delta(\mycal{T})}.
\]
Here $\gamma : U \to X$ is the scheme representing $\mycal{T}$ and $R$
is defined as follows. The presentation of $X$ generated by $\gamma$
has ${\mathfrak R} = U\times_{X} U = U\times_{X} K$ since $U$ is an
$H$-torsor. Furthermore, under the identification ${\mathfrak R} =
U\times_{X} K$ the structure maps ${\mathfrak s}$ and ${\mathfrak t}$
become the projection on $U$ and the action of $K$ respectively. In
other words $\xymatrix@1{{\mathfrak R}\, \ar@<.5ex>[r]^-{{\mathfrak
s}} \ar@<-.5ex>[r]_-{{\mathfrak t}} & U}$ is the transformation
groupoid for the action of $K$ on $U$ and $X = U/{\mathfrak R} =
U/K$. To get the presentation $\xymatrix@1{R\, \ar@<.5ex>[r]^-{{\bs}}
\ar@<-.5ex>[r]_-{{\bt}} & U}$ of $\delta(\mycal{T})$ we can just take
the transformation groupoid for the action of $G$ on $U$ (where $H$
acts trivially), i.e. take $R := U\times_{X} G$ with $\bs$ and $\bt$
being again the projection and the action maps. Equivalently we may
take $R \to {\mathfrak R}$ to be the trivial ${\mathfrak
G}$-torsor and check that it satisfies the biextension and
trivialization conditions. In particular we get that
$\delta(\mycal{T})$ is a quotient gerbe - it is identified as the
quotient
\[
\delta(\mycal{T}) = [U/R] = [U/G]
\]
of the space $U$ by the group scheme $G$, where $G$ acts with a
stabilizer $H$ at each point.

\

\medskip

\begin{ex} \label{ex-azumaya}
As a special case of the above we obtain the Azumaya presentation of
an ${\mathcal O}_{X}^{\times}$-gerbe. Let $P$ be a ${\mathbb
P}^{n-1}$ bundle on a scheme $X$ and let ${\mathfrak P}$ denote the
corresponding sheaf of sections. The bundle $P$ is associated to a unique
${\mathbb P}GL_{n}$-bundle $U \to X$ (the frame bundle of
$P$) whose sheaf of sections $\mycal{U}$ is naturally a torsor
over ${\mathbb P}GL_{n}({\mathcal O}_{X})$. The image of $\mycal{U}$
under the coboundary map for the sequence
\[
1 \to {\mathcal O}_{X}^{\times} \to GL_{n}({\mathcal O}_{X}) \to
{\mathbb P}GL_{n}({\mathcal O}_{X})  \to 1
\]
is an ${\mathcal O}_{X}^{\times}$-gerbe $_{{\mathfrak P}}X$ on $X$ which
comes together with a {\em right Azumaya presentation} of $_{{\mathfrak P}}X$:
\[
\xymatrix@1{R^{r}\, \ar@<.5ex>[r]^-{{\bs}} \ar@<-.5ex>[r]_-{{\bt}} & U
\ar[r] & _{{\mathfrak P}}X}
\]
where $R^{r} := U\times GL_{n}({\mathbb C})$, $\bs : R^{r} \to U$ is the
projection and $\bt : R^{r} \to U$ is the (right) action of
$GL_{n}({\mathbb C})$ on $U$. Alternatively one may consider the sheaf
${\mycal A}_{{\mathfrak P}}$ of Azumaya algebras corresponding to
${\mathfrak P}$. The subsheaf ${\mycal A}_{{\mathfrak P}}^{\times} \subset
{\mycal A}_{{\mathfrak P}}$ of invertible elements in ${\mycal
A}_{{\mathfrak P}}$ is representable by an affine group scheme
$A_{{\mathfrak P}}^{\times} \to X$ which acts simply transitively on the
left on the frame bundle $U \to X$. Using this group scheme we get a
{\em left Azumaya presentation} of $_{{\mathfrak P}}X$:
\[
\xymatrix@1{R^{l}\, \ar@<.5ex>[r]^-{{\bs}} \ar@<-.5ex>[r]_-{{\bt}} & U
\ar[r] & _{{\mathfrak P}}X}
\]
where $R^{l} := A_{{\mathfrak P}}^{\times}\times_{X} U$, $\bs : R^{l} \to
U$ is the projection and $\bt : R^{l} \to U$ is the (left) action of
$A_{{\mathfrak P}}^{\times}$ on $U$.

The same gerbe $_{{\mathfrak P}}X$ has yet another presentation, called
the {\em Brauer-Severi presentation}. Here the atlas is $P$ itself and
the relations are the total space of the punctured line bundle
${\mathcal O}(1,-1)^{\times}$ on $P\times_{X} P$.
\end{ex}

\

\medskip

\noindent
It is often useful to describe the sheaves on a gerbe as cartesian
sheaves on the simplicial space generated by a presentation or
equivalently as descent datum for a presentation.
Concretely, given a flat presentation
\begin{equation} \label{eq-present}
\xymatrix@1{R\, \ar@<.5ex>[r]^-{{\bs}}
\ar@<-.5ex>[r]_-{{\bt}} & U \ar[r] & \leftidx{_{\alpha}}{X}{}}
\end{equation}
of an $\mycal{H}$-gerbe $\leftidx{_{\alpha}}{X}{}$ on $X$ we have a
simple interpretation (see e.g. \cite[Section~3.3]{deligne-tannaka} or
\cite[Propositions 12.8.2 and 13.2.4]{laumon-stacks} for details and
proofs) for the category of sheaves on $\leftidx{_{\alpha}}{X}{}$: a
sheaf of ${\mathcal O}_{\leftidx{_{\alpha}}{X}{}}$-modules
(respectively: a line bundle, a vector bundle) on
$\leftidx{_{\alpha}}{X}{}$ can be identified with a pair $(F,j)$,
where $F$ is a sheaf of ${\mathcal O}_{U}$-modules (respectively: a
line bundle, a vector bundle, a complex of sheaves, ect.) on $U$, and
$j : \bs^{*}F \widetilde{\to} \bt^{*}F$ is an isomorphism of sheaves
on $R$ satisfying the cocycle condition\footnote{To write this
condition one uses the natural identifications
\[
p_{1}^{*}\bs^{*}F \xrightarrow{p_{1}^{*}j}  p_{1}^{*}\bt^{*}F =
p_{2}^{*}\bs^{*}F \xrightarrow{p_{2}^{*}j} p_{2}^{*}\bt^{*} F,
\]
and
\[
p_{1}^{*}\bs^{*}F = \bm^{*}\bs^{*} F \xrightarrow{\bm^{*}j}
\bm^{*}\bt^{*} F = p_{2}^{*}\bt^{*} F,
\]
provided by the groupoid axioms.
}
\begin{equation} \label{eq-cocycle0}
(p_{2}^{*}j)\circ (p_{1}^{*}j) = \bm^{*}j
\end{equation} 
on
$R\times_{\bt, U, \bs} R$ and the normalization $\be^{*}(j) =
\op{id}_{F}$ on $U$. Assume now
that $\mycal{H} = H({\mathcal O}_{X})$ for some complex
reductive abelian group $H$. In this situation we can recast the above
description
of sheaves on $\leftidx{_{\alpha}}{X}{}$ in terms of the presentation
$\xymatrix@1{{\mathfrak
R}\, \ar@<.5ex>[r]^-{{\mathfrak s}}
\ar@<-.5ex>[r]_-{{\mathfrak t}} & U}$  generated by $\gamma : U \to X$
and the $H$-torsor $\pi : R \to {\mathfrak R}$. Given a sheaf of
modules $(F,j)$ on $\leftidx{_{\alpha}}{X}{}$ we can use the map $\pi : R \to
{\mathfrak R}$ to  push the isomorphism $j$ down to ${\mathfrak
R}$. Decomposing according to the characters $\widehat{H}$ of $H$ we see that
$\pi_{*}(j)$ corresponds to a family $\{ {\mathfrak j}_{\chi} \}_{\chi
\in \widehat{H}}$ of isomorphisms, where
\[
{\mathfrak j}_{\chi} : {\mathfrak s}^{*}F\otimes (\pi_{*}{\mathcal
O}_{R})_{\chi} \widetilde{\to} {\mathfrak t}^{*}F.
\]
The category of sheaves of modules on $\leftidx{_{\alpha}}{X}{}$ is
therefore ``graded'' by the character group $\widehat{H}$. A sheaf of
weight $0 \in \widehat{H}$ is just a sheaf of modules on $X$. In case
$H = {\mathbb G}_{m}$ we have $\widehat{H} = {\mathbb Z}$ and the
sheaves of weight $n$ are precisely the sheaves of weight $n$ in the
sense of section~\ref{ssubsec-gerbes}. In particular, the sheaves of
weight $1$ are the $\alpha$-twisted sheaves on $X$. This observation
leads to a very concrete description of the weight one sheaves on a
${\mathbb G}_{m}$-gerbe. Starting with a presentation
\eqref{eq-present} of a ${\mathbb G}_{m}$-gerbe
$\leftidx{_{\alpha}}{X}{}$ on $X$ write $\mycal{L} \to {\mathfrak R} =
U\times_{X} U$ for the line bundle associated to the ${\mathbb
G}_{m}$-torsor $R \to {\mathfrak R}$ via the tautological character
$\op{id} : {\mathbb G}_{m} \to {\mathbb G}_{m}$. The groupoid
condition on the presentation \eqref{eq-present} gives us a
biextension isomorphism $p_{12}^{*}\mycal{L}\otimes
p_{23}^{*}\mycal{L} = p_{13}^{*}\mycal{L}$ on $U\times_{X} U\times_{X}
U$ and so a sheaf on $\leftidx{_{\alpha}}{X}{}$ is the same thing as a
sheaf $F$ on $U$ equipped with an $\mycal{L}$-twisted descent datum on
$U\times_{X} U$, i.e. with an isomorphism
\[
p_{1}^{*}F \xrightarrow{j} p_{2}^{*}F \otimes {\mycal L},
\]
of sheaves on $U\times_{X} U$, satisfying the cocycle condition
\begin{equation} \label{eq-cocycle}
p_{13}^{*}j = (p_{23}^{*}j\otimes \op{id}_{p_{23}^{*}\mycal{L}})\circ
p_{12}^{*}j
\end{equation}
on $U\times_{X} U\times_{X} U$. Note that in writing
\eqref{eq-cocycle} we had to use the biextension isomorphism for $\mycal{L}$.

\begin{ex} \label{ex-azumaya-sheaves} Specializing the previous
discussion to the case of the Azumaya gerbe $_{\mathfrak P}X$ of
example~\ref{ex-azumaya} we get a natural identification of the
category $D^{b}_{n}(_{\mathfrak P}X)$ with the derived category of
complexes of quasi-coherent sheaves on $X$ equipped with an action of
the Azumaya algebra $_{\mathfrak P}A$ and such that the center
${\mathcal O}_{X}^{\times}$ of $_{\mathfrak P}A^{\times}$ acts
on the cohomology sheaves with character $n$.
\end{ex}

\subsubsection{Brauer groups} \label{ssubsec-Brauer}

Since the ${\mathcal O}^\times$ gerbes are naturally classified by
elements in cohomological Brauer groups, it will be helpful to have an
overview of the different variants of the Brauer group of a complex
space before discussing properties of individual gerbes.

Below we are going to discuss three versions of the Brauer group of a
ringed space $Z$: Azumaya ($Br(Z)$), geometric ($Br_{\op{geom}}(Z)$),
and cohomological ($Br'(Z)$). Each of these makes sense in either the
etale or the analytic topology on $Z$. In particular, for a complex
algebraic space $Z$ we have a diagram:
\[
\xymatrix@M+4pt{ Br(Z) \ar[r] \ar[d] & Br(Z)_{\op{geom}}
\ar[r]  \ar[d] & Br'(Z)
\ar[d] \\
Br_{an}(Z) \ar[r] & Br_{an}(Z)_{\op{geom}} \ar[r] & Br'_{an}(Z).
}
\]
When $Z$ is {\em smooth}, the following facts are known:
\begin{itemize}
\item[(i)] all the maps in this diagram are injective;
\item[(ii)] $Br'(Z)$ is torsion by the purity theorem from \cite{groth-bg3};
\item[(iii)] $\op{im}[Br'(Z) \to Br'_{an}(Z)]$ coincides (see
\cite{milne-book}) with the torsion subgroup of $Br'_{an}(Z)$;
\end{itemize}

\

Grothendieck has conjectured that the inclusion $Br(Z)
\hookrightarrow Br'(Z)$ is an isomorphism for all smooth
quasi-projective schemes. This may hold also for separated normal
$Z$. The validity of the conjecture was established in the algebraic
setting in \cite{gabber,hoobler,schroer} for arbitrary curves, for
normal separated algebraic surfaces, for abelian varieties, for smooth
toric varieties and for separated unions of two affine varieties. The
analogous conjecture in the analytic case is virtually unexplored. The
only general result to date \cite{hyubrechts-schroeer} concerns
analytic K3 surfaces and asserts that every torsion class in
$Br_{an}'(X)$ of an analytic K3 surface $X$ comes from an Azumaya
algerba on $X$.

\begin{rem} \label{rem-quotient-gerbes} As a corollary of fact (i) and
Grothendieck's conjecture, we get that $Br(Z) = Br(Z)_{\op{geom}}$ for
a smooth $Z$. This corollary is known to hold \cite{ehkv} in many
cases in which the Grothendieck conjecture is still unknown. In
fact, for a normal Noetherian scheme, the result of \cite[Theorem~3.6]{ehkv}
characterizes the image of $Br(Z)$ in $Br(Z)_{\op{geom}}$ as the
algebraic-geometric gerbes for which one can find a flat presentation
\[
\left(\xymatrix@1{R\ar@<.5ex>[r]^-{\bs}
\ar@<-.5ex>[r]_-{\bt} & U}, \bm,\bi,\be,\bgamma \right)
\]
with a projective structure map $\bgamma :
U \to Z$, or equivalently as the classes of ${\mathbb G}_{m}$
gerbes of quotient type (i.e. a quotient of an algebraic space by an
affine algebraic group). In general this characterization seems to be
optimal since there are examples of quotient gerbes on non-separated
surfaces whose isomorphism class is not represented by an Azumaya
algebra, and examples of infinite order elements in $Br'(Z)$ for a
normal separated $Z$ which are represented by algebraic-geometric
gerbes but are not quotient gerbes \cite[Examples 2.21 and
3.12]{ehkv}.
\end{rem}

\

If $Z$ is a complex scheme, then the {\em Azumaya Brauer group}
$Br(Z)$, is defined \cite{groth-bg1} as the group of Morita
equivalence classes of sheaves of Azumaya algebras on $Z$. Recall
\cite{groth-bg1} that an Azumaya algebra on $Z$ is a coherent sheaf of
algebras which locally in the etale topology on $Z$ is isomorphic
to the endomorphisms algebra of an algebraic vector bundle on $Z$.
Two Azumaya algebras ${\mathcal A}$ and ${\mathcal B}$ are called
Morita equivalent if etale locally on $Z$ we can find vector
bundles $E$ and $F$ so that the sheaves of algebras ${\mathcal
A}\otimes \mycal{E}nd(E)$ and ${\mathcal B}\otimes \mycal{E}nd(F)$ are
isomorphic. Morita equivalence classes of Azumaya algebras form a
commutative group under the operation of tensoring over
${\mathcal O}_{Z}$; the inverse is given by the opposite algebra.

The Skolem-Noether theorem \cite[Proposition~2.3]{milne-book}
implies that the Azumaya algebras of rank $n^{2}$ are classified by
elements in $H^{1}_{\et}(Z,{\mathbb P}GL(n))$. The short exact
sequence of groups schemes over $Z$:
\[
1 \to {\mathbb G}_{m} \to GL(n) \to {\mathbb P}GL(n) \to 1,
\]
gives rise to a coboundary map
\begin{equation} \label{eq-coboundary-brauer-map}
H^{1}_{\et}(Z,{\mathbb P}GL(n)) \to H^{2}_{\et}(Z,{\mathbb G}_{m}).
\end{equation}
The image of $a \in H^{1}_{\et}(Z,{\mathbb P}GL(n))$ under this
coboundary map is an $n$-torsion class in $H^{2}_{\et}(Z,{\mathbb
G}_{m})$ which is the obstruction to representing $a$ by the
endomorphism algebra of a rank $n$ vector bundle. In particular the
map \eqref{eq-coboundary-brauer-map} induces a homomorphism
\begin{equation} \label{eq-brauer-map}
Br(Z) \to H^{2}_{\et}(Z,{\mathbb G}_{m})_{\op{tor}} \subset
H^{2}_{\et}(Z,{\mathbb G}_{m}).
\end{equation}
When $Z$ is smooth, the homomorphism \eqref{eq-brauer-map} is known to
be injective \linebreak \cite[Theorem~IV.2.5]{milne-book}. This
suggests that the Brauer classes are intimately related to elements in
$H^{2}_{\et}(Z,{\mathbb G}_{m})$ and so one defines the {\em algebraic
cohomological Brauer group}:
\[
Br'(Z) :=
H^{2}_{\et}(Z,{\mathbb G}_{m}).
\]
Recall that by fact (ii) the  group $H^{2}_{\et}(Z,{\mathbb G}_{m})$ is
purely torsion. 
As explained in section~\ref{ssubsec-presentations}, Azumaya algebras
give rise to groupoid presentations of ${\mathbb G}_{m}$-gerbes on
$Z$. In other words, for a smooth $Z$ the inclusion
\eqref{eq-brauer-map} can be refined to a sequence of inclusions:
\[
Br(Z) \hookrightarrow Br(Z)_{\op{geom}}
{\hookrightarrow} Br'(Z) = H^{2}_{\et}(Z,{\mathbb G}_{m}) =
H^{2}_{\et}(Z,{\mathbb G}_{m})_{\op{tor}},
\]
where $Br(Z)_{\op{geom}}$ denotes the group of equivalence
classes of algebraic-geometric ${\mathbb G}_{m}$-gerbes on $Z$. Recall
that a ${\mathbb G}_{m}$-gerbe is {\em algebraic geometric} if it is
an algebraic stack in the sense of Artin, i.e. if it admits a flat
(equivalently, a smooth) groupoid presentation \cite{artin-stacks}.

By analogy we define the {\em analytic Azumaya Brauer group}
$Br_{an}(Z)$, and the {\em analytic geometric Brauer group}
$Br_{an}(Z)_{\op{geom}}$ of an analytic space $Z$, as the groups of
Morita equivalent classes of analytic Azumaya algebras on $Z$ and of
isomorphism classes of analytic geometric ${\mathcal
O}_{Z_{an}}^{\times}$-gerbes respectively. The isomorphism type of an
${\mathcal O}_{Z_{an}}^{\times}$-gerbe is determined by a class in the
{\em analytic cohomological Brauer group}:
\[
Br_{an}'(Z) := H^{2}_{an}(Z,{\mathcal O}_{Z}^{\times}).
\]
In other words, the classifying map
\[
Br_{an}(Z)_{\op{geom}} \hookrightarrow Br_{an}'(Z)
\]
is injective.

\

The analytic cohomological Brauer group can be studied via the
exponential sequence:
\[
0 \to {\mathbb Z} \to {\mathcal O} \stackrel{\exp}{\to} {\mathcal
O}^{\times} \to 1.
\]
The corresponding cohomology sequence gives
\[
\xymatrix@1{
0 \ar[r] & H^{2}_{an}(Z,{\mathcal O}_{Z})/H^{2}(Z,{\mathbb
Z}) \ar[r] & Br'_{an}(Z) \ar[r] & \ker\left[H^{3}(Z,{\mathbb Z}) \to
H^{3}_{an}(Z,{\mathcal O}_{Z})\right] \ar[r] & 0.
}
\]
This of course is analogous to the usual description of the Picard
group:
\[
\xymatrix@M+3pt{ 0 \ar[r] & H^{1}_{an}(Z,{\mathcal O}_{Z})/H^{1}(Z,{\mathbb
Z}) \ar[r] & \op{Pic}(Z) \ar[r] & \ker\left[H^{2}(Z,{\mathbb Z}) \to
H^{2}_{an}(Z,{\mathcal O}_{Z})\right] \ar[r]
\ar@{=}[d]^-*+{\; \frame{\txt<9pc>{
{\small when $Z$ is compact and K\"{a}hler}
}}} & 0. \\ & & &
H^{1,1}_{{\mathbb Z}}(Z) & }
\]
In addition, if $Z$ is compact and K\"{a}hler, the Hodge theorem
implies that
\[
\op{im}[H^{1}(Z,{\mathbb Z}) \to H^{1}_{an}(Z,{\mathcal
O}_{Z})] \subset H^{1}_{an}(Z,{\mathcal O}_{Z}),
\]
is a discrete subgroup of maximal rank. Hence, we can identify the
connected component of $\op{Pic}(Z)$ with the quotient of its tangent
space $H^{1}_{an}(Z,{\mathcal O}_{Z})$ by $H^{1}(Z,{\mathbb Z})$.
In the case of $Br'_{an}(Z)$, there is still a `tangent space':
$H^{2}_{an}(Z,{\mathcal O}_{Z})$, but it is divided by the typically
non-discrete subgroup
\[
\op{im}[H^{2}(Z,{\mathbb Z}) \to H^{2}_{an}(Z,{\mathcal O}_{Z})]
\subset H^{2}_{an}(Z,{\mathcal O}_{Z}),
\]
and so there is no good (=separated) topology on $Br'_{an}(Z)$.

In the special case when $Z$ is a $K3$ surface,
we get  that $Br_{an}(Z)_{\op{geometric}} = Br_{an}'(Z)$ is the
quotient of the one dimensional vector space $H^{2}_{an}(Z,{\mathcal
O}_{Z})$ by the lattice dual to the transcendental lattice of $Z$,
i.e. by $H^{2}(Z,{\mathbb Z})/H^{1,1}_{{\mathbb Z}}(Z)$. Notice that
for a very general analytic $K3$ this lattice has rank
$22$ and for a very general algebraic $K3$ it has rank $21$.

More precisely, one defines the transcendental lattice $\Tr_{Z}$ of a $K3$
surface $Z$ by the short exact sequence:
\[
0 \to \Tr_{Z} \to H^{2}(Z,{\mathbb Z}) \to H^{1,1}_{{\mathbb
Z}}(Z)^{\vee} \to 0.
\]
In other words, $\Tr_{Z}$ is the sublattice of $H^{2}(Z,{\mathbb Z})$
consisting of classes perpendicular to all classes of curves in
$Z$. The dual sequence reads:
\[
0 \to H^{1,1}_{{\mathbb Z}}(Z) \to H^{2}(Z,{\mathbb Z}) \to
\Tr_{Z}^{\vee} \to 0,
\]
and we have a natural map
\[
\op{Hom}_{{\mathbb Z}}(\Tr_{Z},{\mathbb R}) \cong H^{2}(Z,{\mathbb
R})/(H^{1,1}_{{\mathbb Z}}(Z)\otimes {\mathbb R}) \twoheadrightarrow
H^{2}(Z,{\mathbb R})/H^{1,1}_{{\mathbb R}}(Z) \cong
H^{2}_{an}(Z,{\mathcal O}_{Z}).
\]
This leads to the following commutative diagram with exact rows and columns:
\[
\xymatrix{
& & 0 \ar[d] & 0 \ar[d] & \\
& & \op{Hom}_{{\mathbb Z}}(\Tr_{Z},{\mathbb Z}) \ar[r]^-{\cong} \ar[d]
& \Tr_{Z}^{\vee} \ar[d] & \\
0 \ar[r] & H^{1,1}_{{\mathbb R}}(Z)/(H^{1,1}_{{\mathbb Z}}(Z)\otimes
{\mathbb R}) \ar[r] \ar@{=}[d] & \op{Hom}_{{\mathbb
Z}}(\Tr_{Z},{\mathbb R}) \ar[r] \ar[d] & H^{2}_{an}(Z,{\mathcal O}_{Z})
\ar[r] \ar[d] &  0 \\
0 \ar[r] & H^{1,1}_{{\mathbb R}}(Z)/(H^{1,1}_{{\mathbb Z}}(Z)\otimes
{\mathbb R}) \ar[r]  & \op{Hom}_{{\mathbb
Z}}(\Tr_{Z},{\mathbb R}/{\mathbb Z}) \ar[r] \ar[d] & Br_{an}'(Z)
\ar[r] \ar[d] &  0 \\
& & 0 & 0 &
}
\]
The bottom row explicates $Br_{an}(Z) = Br_{an}'(Z)$ as the quotient
of the real torus \linebreak $\op{Hom}_{{\mathbb Z}}(\Tr_{Z},{\mathbb
R}/{\mathbb Z})$ by the vector space $H^{1,1}_{{\mathbb
R}}(Z)/(H^{1,1}_{{\mathbb Z}}(Z)\otimes {\mathbb R})$, embedded in it
as a (usually dense) subgroup. Note that this vector space does not
contain any torsion points of the torus. Equivalently the restricted
map
\[
\op{Hom}_{{\mathbb Z}}(\Tr_{Z},{\mathbb Q}/{\mathbb Z})
\widetilde{\to} Br_{an}(Z)_{\op{tor}}
\]
is an isomorphism. When $Z$ happens to be an algebraic $K3$ surface we
have a natural identification $Br(Z) = Br_{an}(Z)_{\op{torsion}}$
and so we recover the standard interpretation of elements of the
algebraic Brauer group of $Z$ as a homomorphism from the
transcendental lattice of $Z$ to ${\mathbb Q}/{\mathbb Z}$ (see
e.g. \cite[Lemma 5.4.1]{caldararu-thesis} or \cite{caldararu-k3}).

\subsection{Tate-Shafarevich groups and genus one fibrations}
\label{subsec-TSh-generalities}

In this section we review some basic facts about twisted forms of a
given elliptic fibration over an analytic space $B$. For more details
the reader is referred to the excellent references
\cite{dolgachev-gross} and \cite{nakayama-global}.
First we recall some terminology and set up the notation.

For us a {\em genus one fibration} will always mean a holomorphic map
$\pi : X \to B$ between normal analytic varieties whose generic fiber
is a smooth curve of genus one. We define an {\em elliptic fibration}
to be a genus one fibration equipped with a
holomorphic section $\sigma : B \to X$ of $\pi$. Note that this is
slightly more restrictive than the conventional notion of an elliptic
fibration used in say \cite{dolgachev-gross}, \cite{nakayama-global},
where only the existence of a meromorphic section of $\pi$ is
required. A genus one fibration will be called (relatively) minimal if
$X$ has at most terminal singularities and if the canonical class
$K_{X}$ is $\pi$-nef.

Let now $X$ and $B$ be normal analytic varieties and let
\[
\xymatrix@1{X
\ar[r]^-{\pi} & B \ar@/^0.5pc/[l]^-{\sigma}}
\]
be an elliptic fibration on $X$. Let $D \subset B$ denote the
discriminant divisor of $\pi$ and let $B^{o} := B-D$, $B^{oo} := B -
\op{Sing}(D)$. The corresponding inclusions are denoted by $\imath : D
\hookrightarrow B$, $\jmath^{o} : B^{o} \hookrightarrow B$ and
$\jmath^{oo} : B^{oo} \hookrightarrow B$. We also put $X^{o} :=
X\times_{B} B^{o}$, $X^{oo} := X\times_{B} B^{oo}$ and $\pi^{o} :=
\pi_{|X^{o}}$, $\pi^{oo} := \pi_{|X^{oo}}$.

Sometimes we may need to require the additional genericity
assumption that $X$ is smooth and
that  $\pi : X \to B$ is Weierstrass.

\begin{rem} \label{rem-minimal model} {\bf (i)} When $X$ is a surface,
  the genericity assumption implies in 
particular that all the singular fibers of $\pi$ are of Kodaira types
$I_{1}$ or $II$, i.e. they are nodes and cusps.

\

\noindent
{\bf (ii)} In this paper we will always deal with a situation in
which $X$ is smooth and either $\pi$ is smooth or $X$ is a surface and
$\pi$ has at worst
$I_{1}$ fibers. We have included in the present
discussion the more general case of an arbitrary Weierstrass
$\pi$ with a smooth total space, because of the potential applications
of our duality construction to genus one fibered Calabi-Yau manifolds
of arbitrary dimension. This however goes beyond the scope of the
present work and will be the subject of a future paper.
\end{rem}

\

\medskip

Let $X^{\sharp}\subset X$ denote the regular locus of $\pi$, viewed as
an abelian group scheme over $B$. Denote by $\grX_{an}$ the
corresponding sheaf of abelian groups in the analytic topology on
$B$. When $B$ and $X$ happen to underly complex algebraic varieties we
will write $\grX$ for the etale sheaf of sections of $X^{\sharp}
\to B$.

\medskip

\

\noindent
The {\em analytic Weil-Ch\^{a}telet group} $WC_{an}(X)$ of $X$ is the
group of bimeromorphism classes of analytic genus one fibrations $Y
\to B$ such that:
\begin{itemize}
\item $Y\times_{B} B^{o} \to B^{o}$ is bimeromorphic to a smooth genus
one fibration;
\item The relative Jacobian fibration $\op{Pic}^{0}(Y/B)$ is
bimeromorphic to $X^{\sharp}$ (and hence to $X$). Note that this
definition makes sense since for a suitably chosen dense open subset
$U \subset B$ the (sheafification of the) presheaf
$\mycal{P}ic^{0}(Y/U)$ of relative Picard groups along the fibers of
$Y\times_{B} U \to U$ is representable by an analytic space.
\end{itemize}

\medskip

\

\noindent
The {\em analytic Tate-Shafarevich group } $\TSh_{an}(X)$ of $X$ is the
subgroup of $WC_{an}(X)$ consisting of elements $\alpha \in
WC_{an}(X)$ such that for any representative $Y \to B$
of $\alpha$ and any point $b \in B$ one can find an analytic
neighborhood $b \in U \subset B$ so that
$Y\times_{B}U \to U$ has a meromorphic
section. This implies that $Y \to B$ has no multiple fibers in
codimension one.

The group $\TSh_{an}(X)$ can be described
cohomologically \cite{nakayama-global} as follows. Assume that
$X^{oo}$ is a smooth space. Then by
\cite[Proposition~5.5.1]{nakayama-global}, the
natural classifying map
\begin{equation} \label{eq-Sh-classifying}
\TSh_{an}(X) \to  H^{1}_{an}(B,j^{oo}_{*}j^{oo*}\grX_{an}).
\end{equation}
is injective. Furthermore  if $B^{oo} = B$, or if $\pi$ is Weierstrass
with a smooth total space, then the map 
\eqref{eq-Sh-classifying} is an isomorphism
\cite[Proposition~5.5.1]{nakayama-global}. In addition one knows (see
e.g. \cite[Theorem~5.4.9]{nakayama-global}) that
under the same assumptions,  the sheaf $j^{oo}_{*}j^{oo*}\grX_{an}$
fits in a short exact sequence
\[
0 \to \grX_{an}  \to j^{oo}_{*}j^{oo*}\grX_{an} \to
(R^{2}\pi_{*}{\mathbb Z}_{X}/\imath_{*}\imath^{!} R^{1}\pi_{*}{\mathcal
O}_{X}^{\times})_{\text{torsion}} \to 0.
\]
Since by definition the sheaf $(R^{2}\pi_{*}{\mathbb
Z}_{X}/\imath_{*}\imath^{!} R^{1}\pi_{*}{\mathcal
O}_{X}^{\times})_{\text{torsion}}$ is supported on the multiple fiber
sublocus of $D$, it follows that in the absence of multiple fibers,
i.e. under our definition of an elliptic fibration we have an isomorphism:
\begin{equation} \label{eq-coho-interpretation}
\TSh_{an}(X) \cong H^{1}_{an}(B,\grX_{an}).
\end{equation}
In the remainder of this paper we will always assume tacitly that the
isomorphism \eqref{eq-coho-interpretation} holds, in fact we will
assume that either $\pi$ is smooth or that $X$ is a surface.

Because of
this cohomological interpretation  we can view the elements in
$\TSh_{an}(X)$ simply as $\grX_{an}$-torsors. This definition of
$\TSh_{an}(X)$ is consistent with the usual definition of the algebraic
Tate-Shafarevich group \cite{groth-bg1,groth-bg2,groth-bg3} and
\cite{dolgachev-gross}. When $X$ is a surface we can also interpret
(by  compactifying a genus one fibration  and then
choosing a Weierstrass smooth model over $B$)
the elements in $\TSh_{an}(X)$  as
smooth analytic surfaces equipped with a genus one fibration over
$B$.

\medskip

\

\noindent
The {\em algebraic Weil-Ch\^{a}telet and Tate-Shafarevich groups} $WC(X)$
and $\TSh(X)$ are defined in a similar manner
\cite{groth-bg1,groth-bg2,groth-bg3} and \cite{dolgachev-gross} with
the etale topology replacing the analytic one. Furthermore,
the analysis carried out in \cite[Section~1]{dolgachev-gross} implies,
that under the assumption that $X$ and $B$ are both smooth and
that $\pi$ has a regular section, the algebraic Tate-Shafarevich group
can be interpreted cohomologically as
\[
\TSh(X) = H^{1}_{\et}(B,\grX),
\]
i.e. the elements in $\TSh(X)$ can be viewed as algebraic spaces $Y
\to B$ which are $\grX$-torsors.

\

\medskip

\noindent
Given an element $\alpha \in \TSh_{an}(X)$ (or $\alpha \in \TSh(X)$)
we denote by $X_{\alpha}^{\sharp}$ the analytic (or algebraic) space
representing the torsor $\alpha$ and by $\pi_{\alpha}^{\sharp} :
X_{\alpha}^{\sharp} \to B$ the corresponding projection. Following
\cite{dolgachev-gross} we say that a morphism of analytic (algebraic)
spaces $Y \to B$ is a {\em good model} for $\alpha$ if $Y \to B$ is
bimeromorphic to $X_{\alpha}^{\sharp} \to B$, $Y$ is smooth and the
map $Y \to B$ is proper and flat.

\begin{rem} \label{rem-good}
Note that when
$\pi$ is smooth $X_{\alpha}^{\sharp}$ is itself a good model for
$\alpha$ and when $X$ is an arbitrary smooth surface we always have a
preferred good model for $\alpha$, namely the relatively minimal model
of a compactification of $X_{\alpha}^{\sharp}$. When $X$ is of
dimension three the good models of elements in $\TSh(X)$ have been analyzed
in detail, see e.g. \cite{miranda,antonella,dolgachev-gross}. In this
case the good model exists (possibly after blowing up $B$ at finitely
many points) but is not unique. However all good models of a given
$\alpha$ are related by flops and in particular have equivalent
derived categories of coherent sheaves (see
e.g. \cite{bondal-orlov-flops,bridgeland-flops,kawamata-flops}).
\end{rem}

In the cases when $\pi : X \to B$ is smooth or $X$ is a surface we put
\[
\pi_{\alpha} : X_{\alpha} \to B
\]
for the canonical good model of $\alpha$. In particular, if $\pi : X
\to B \cong {\mathbb P}^{1}$ is an elliptic $K3$ surface we have that
$X_{\alpha}$ is a well defined analytic (respectively algebraic) $K3$
surface for any element $\alpha \in \TSh_{an}(X)$ (respectively
$\alpha \in \TSh(X)$).

The meromorphic action of the
analytic group space $X^{\sharp} \to B$ on $X_{\alpha}$
induces a natural meromorphic action map
\[
a_{\alpha} :
X\times_{B} X_{\alpha} \dashrightarrow X.
\]
Furthermore,
given a positive integer $n$ we can consider the sheaf of groups
$\grX_{an}[n] \to B$ consisting of the $n$-torsion
points in $\grX_{an}$.
The sheaf $\grX_{an}[n]$ is represented by a group space
$X^{\sharp}[n]$ which is quasi-finite over $B$.  We
will write $X[n]$ for the closure of $X^{\sharp}[n]$ in $X$ and by an
abuse of notation we will denote the meromorphic map
\[
X[n]\times_{B} X_{\alpha} \dashrightarrow X_{\alpha}
\]
again by $a_{\alpha}$.

Since $X^{\sharp}[n]$ is finite over a dense open set in
$B$ we can form the quotient $X_{\alpha}/X^{\sharp}[n]$
which as an analytic space is well defined up to a
bimeromorphism which respects the genus one fibration. Moreover
$X_{\alpha}/X^{\sharp}[n]$ is naturally a $\grX_{an}$-torsor at the
general point and so represents an element in $\TSh_{an}(X)$. It is
not hard to calculate this element in terms of $\alpha$ and $n$
only. In fact it is clear that $X_{\alpha}/X^{\sharp}[n]$ is
tautologically the same as the quotient
\[
X_{\alpha}^{\times_{B} n}/K,
\]
which by definition represents the element $n\alpha \in
\TSh_{an}(X)$.

Here
\[K = \ker[(X^{\sharp})^{\times_{B} n} \xrightarrow{\op{mult}} X^{\sharp}]
\]
is the kernel of the natural product map corresponding the group law
on $X^{\sharp}$, and the action of $K$ is induced from the
component-wise action of $(X^{\sharp})^{\underset{B}{\times} n}$ on
$(X_{\alpha})^{\underset{B}{\times} n}$.

In particular we have a bimeromorphism $X_{\alpha}/X^{\sharp}[n]
\widetilde{\dashrightarrow} X_{n\alpha}$ which is unique up to an
auto-bimeromorphism of $X_{n\alpha}$, compatible with the genus one
fibration. However, as one can see from the proof of
\cite[Lemma~5.3.3]{nakayama-global}, if we assume that $\pi$ is
relatively minimal with a smooth total space, then all such
auto-bimeromorphisms are holomorphic and are translations by sections
in $\grX_{an}$. So, under this the genericity assumption, we will
have a bimeromorphic identification $X_{n\alpha} =
X_{\alpha}/X^{\sharp}[n]$ and hence a well defined meromorphic map
\[
q_{\alpha}^{n} : X_{\alpha} \dashrightarrow X_{n\alpha}.
\]
If in addition we assume that the fibration $\pi : X \to B$ has a
trivial Mordel-Weil group, then the meromorphic map $q_{\alpha}^{n}$
is canonical and does not depend on any choices.

\

The {\em
index} of an element $\alpha \in \TSh_{an}(X)$ is defined to be the
minimal degree of a global multisection of $\pi_{\alpha}$. We will
denote the index by $\op{ind}(\alpha)$.

Assume now that $B$ and $X$ are quasi-projective.
Since
the element $0 \in \TSh_{an}(X)$ is represented by the algebraic
elliptic fibration  $\pi : X \to B$, it follows that for each $\alpha$
of finite index the space $X_{\alpha}$ admits a dominant meromorphic map
\[
q_{\alpha}^{\op{ind}(\alpha)} : X_{\alpha} \dashrightarrow X
\]
to the algebraic variety $X$. In fact Nakayama
\cite[Proposition~5.5.4]{nakayama-global} proves that such a
$X_{\alpha}$ is bimeromorphic to an algebraic variety and so must be
an algebraic space. Furthermore in the case of surfaces Kodaira shows
\cite{kodaira} that $X_{\alpha}$ is quasi-projective if and only if
$\alpha$ is torsion in $\TSh_{an}(X)$.

\subsection{Complementary fibrations}
\label{subsec-T}

Let $X$ be smooth and let
\[
\xymatrix@1{X
\ar[r]^-{\pi} & B \ar@/^0.5pc/[l]^-{\sigma}}
\]
be a relatively minimal elliptic fibration.
Consider an element $\alpha \in \TSh_{an}(X)$ and a good
representative $\pi_{\alpha} : X_{\alpha} \to B$ for $\alpha$. Our
goal in this section is to describe the cohomological Brauer group
$Br'_{an}(X_{\alpha})$ in terms of the Tate-Shafarevich group
$\TSh_{an}(X)$. For this we need to analyze the relationship between
the sheaf $\mycal{X}_{an}$ and the relative Picard sheaf of
$\pi_{\alpha}$.

If all the fibers of $\pi$ are integral, then
$\mycal{P}ic(X_{\alpha}/B)$ is representable and we have a short exact
sequence of abelian sheaves in the analytic topology:
\begin{equation}
\label{eq-XtoPic}
0 \to \mycal{X} \to \mycal{P}ic(X_{\alpha}/B)
\xrightarrow{\deg_{\alpha}} {\mathbb Z} \to 0,
\end{equation}
where $\deg_{\alpha}$ is the map assigning to each $L \in
\op{Pic}(\pi_{\alpha}^{-1}(U))/\pi_{\alpha}^{*}\op{Pic}(U)$ its degree
along a smooth fiber.

\begin{rem} \label{rem-nonintegral} If we want to allow non-integral
fibers for $\pi$, then $\mycal{P}ic(X_{\alpha}/B)$ becomes
non-representable, but it has a maximal representable quotient
$\mycal{Q}_{\alpha}$ as shown in
e.g. \cite{raynaud-picard} and \cite{dolgachev-gross} in the
algebraic case and \cite{nakayama-global} in the analytic case. The
sheaf of groups $\mycal{Q}_{\alpha}$  is defined as:
\[
\mycal{Q}_{\alpha} := \mycal{P}ic(X_{\alpha}/B)/\mycal{E}_{\alpha},
\]
where $\mycal{E}_{\alpha} \subset \mycal{P}ic(X_{\alpha}/B)$ is a
subsheaf generated by local components of the preimage
$\pi_{\alpha}^{-1}(D)$ of the discriminant $D \subset B$ (see
\cite[Proposition~1.13]{dolgachev-gross} for the precise
statement). Note that when all fibers of $\pi$ are integral we have
$\mycal{E}_{\alpha} = 0$.

In this generality, the short exact sequence \eqref{eq-XtoPic} is
replaced by a commutative diagram with exact rows and columns:
\[
\xymatrix{
& & 0 & 0 &  \\
& & {\mathbb Z} \ar[u] \ar@{=}[r] & {\mathbb Z} \ar[u] & \\
0 \ar[r] & \mycal{X}_{an} \ar[r] & \mycal{Q}_{\alpha}  \ar[r]
\ar[u]^-{\deg_{\alpha}} & R^{2}\pi_{\alpha*} {\mathbb
Z}/\mycal{E}_{\alpha} \ar[u] \ar[r] & 0 \\
0 \ar[r] & \mycal{X} \ar@{=}[u] \ar[r] & \ker(\deg_{\alpha}) \ar[u]
\ar[r] & (R^{2}\pi_{\alpha*} {\mathbb
Z}/\mycal{E}_{\alpha})_{\op{torsion}} \ar[u] \ar[r] & 0 \\
& & 0 \ar[u] & 0 \ar[u] &
}
\]
Note also that $(R^{2}\pi_{\alpha*} {\mathbb
Z}/\mycal{E}_{\alpha})_{\op{torsion}}$ is supported on $D$ and that
its fibers at smooth points of a component of $D$ parameterizing
Kodaira fibers of type $I_{n}$ are isomorphic to ${\mathbb Z}/n$.
\end{rem}

\

\bigskip

\noindent
Fix now  an element $\alpha \in \TSh_{an}(X)$. Under some mild
assumptions on $\alpha$ we will construct a natural map
$T_{\alpha} : \TSh_{an}(X) \to Br'_{an}(X_{\alpha})$, which will allow
us to compare the Tate-Shafarevich and Brauer groups. The
existence of $T_{\alpha}$ is established in the following lemma.

\begin{lem} \label{lem-Talpha} Assume that $X$ is smooth, $\pi$  has
integral fibers and $Br'_{an}(B) = 0$. Assume further that
\begin{equation} \label{eq-vanishing-kernel}
\ker\left(H^{3}_{an}(B,{\mathcal O}_{B}^{\times})
\xrightarrow{\pi_{\alpha}^{*}} H^{3}_{an}(X_{\alpha},{\mathcal
O}_{X_{\alpha}}^{\times})\right) = 0
\end{equation}
Then there is a canonical homomorphism $T_{\alpha}$ which fits in an
exact sequence of abelian groups:
\[
0 \to {\mathbb Z}/\op{ind}(\alpha) \to \TSh_{an}(X)
\xrightarrow{T_{\alpha}} Br'_{an}(X_{\alpha}) \to H^{1}(B,{\mathbb
Z}).
\]
\end{lem}
{\bf Proof.}  The long exact sequence of \eqref{eq-XtoPic} gives
\[
\xymatrix@R-4pt{
{H^{0}(B,{\mathbb Z})}/{H^{0}_{an}(B,\mycal{P}ic(X_{\alpha}/B))}
\ar@{^{(}->}[r] \ar@{=}[d] & H^{1}_{an}(B,\mycal{X}) \ar@{=}[d] \ar[r] &
H^{1}_{an}(B,\mycal{P}ic(X_{\alpha}/B)) \ar[r] & H^{1}(B,{\mathbb Z})
\\
{\mathbb Z}/\op{ind}(\alpha) & \TSh_{an}(X) & &
}
\]
So it suffices to find an identification
\begin{equation} \label{eq-Br}
H^{1}_{an}(B,\mycal{P}ic(X_{\alpha}/B)) \cong
Br'_{an}(X_{\alpha}).
\end{equation}
Consider now
the Leray spectral sequence for $\pi_{\alpha} :
X_{\alpha} \to B$ and the sheaf ${\mathcal O}^{\times}_{X_{\alpha}}$,
which
has only two non-zero rows, so it also becomes a long exact sequence:
\[
\xymatrix@1@C-2pt{ Br'_{an}(B) \ar[r] & Br'_{an}(X_{\alpha}) \ar[r] &
H^{1}_{an}(B,\mycal{P}ic(X_{\alpha}/B)) \ar[r] &
\ker\left(H^{3}_{an}(B,{\mathcal O}_{B}^{\times})
\xrightarrow{\pi_{\alpha}^{*}} H^{3}_{an}(X_{\alpha},{\mathcal
O}_{X_{\alpha}}^{\times})\right).}
\]
The assumption $\ker\left(H^{3}_{an}(B,{\mathcal O}_{B}^{\times})
\xrightarrow{\pi_{\alpha}^{*}} H^{3}_{an}(X_{\alpha},{\mathcal
O}_{X_{\alpha}}^{\times})\right) = 0$ thus immediately implies the
identification \eqref{eq-Br} and so the lemma is proven. \hfill $\Box$

\

\bigskip

The lemma has the following immediate corollary:

\begin{cor} \label{cor-surface-3fold}
Assume that $X$ is a smooth projective surface
and $\pi$ has integral fibers. Then the map
\[
T_{\alpha} : \TSh_{an}(X) \to Br'_{an}(X_{\alpha})
\]
exists for all $\alpha \in  \TSh_{an}(X)$.

\end{cor}
{\bf Proof.} The existence of $T_{\alpha}$ is an immediate
consequence of Lemma~\ref{lem-Talpha} since in this case $B$ is a
smooth curve and so $H^{2}_{an}(B,{\mathcal O}^{\times}_{B}) =
H^{3}_{an}(B,{\mathcal O}^{\times}_{B}) = 0$ for 
dimension reasons. \hfill $\Box$

\

%
%
%

\

\bigskip

\noindent
If the vanishing assumption \eqref{eq-vanishing-kernel} does not hold,
we can still construct a variant of the map $T_{\alpha}$ which is
defined only on a part of the group $\TSh_{an}(X)$:

\begin{lem} \label{lem-m-alpha} Assume that $X$ is smooth, $\pi$  has
integral fibers and $Br'_{an}(B) = 0$. Then:
\begin{itemize}
\item[(i)] if $\alpha$ is $m$-torsion
in $\TSh_{an}(X)$, then there is a group homomorphism (compatible
with $T_{\alpha}$ when the latter exists)
\[
m\TSh_{an}(X) \to Br'_{an}(X_{\alpha}),
\]
from the subgroup $m\TSh_{an}(X) \subset \TSh_{an}(X)$ of
$m$-divisible elements in $\TSh_{an}(X)$ to the cohomological Brauer
group of $X_{\alpha}$;
\item[(ii)] if $\alpha$ is $m$-divisible in $\TSh_{an}(X)$, then
there is a group homomorphism (compatible with $T_{\alpha}$ when the
latter exists)
\[
\TSh_{an}(X)[m] \to Br'_{an}(X_{\alpha})
\]
from the subgroup $\TSh_{an}(X)[m] \subset \TSh_{an}(X)$ of
$m$-torsion elements in $\TSh_{an}(X)$ to the cohomological Brauer
group of $X_{\alpha}$;
\end{itemize}
\end{lem}
{\bf Proof.} For any given $\alpha \in \TSh_{an}(X)$ we have a
composition map
\[
\xymatrix{
\TSh_{an}(X) \ar[rr]^-{d_{\alpha}} \ar[rd] & & H^{3}_{an}(B,{\mathcal
O}_{B}^{\times})  \\
& H^{1}_{an}(B,\mycal{P}ic(X_{\alpha}/B)) \ar[ur] &
}
\]
The assignment $\alpha \mapsto d_{\alpha}$ gives rise to a group
homomorphism
\begin{equation} \label{eq-d}
d : \TSh_{an}(X) \to \op{Hom}_{\mathbb
Z}(\TSh_{an}(X),H^{3}_{an}(B,{\mathcal O}_{B}^{\times})).
\end{equation}
In
particular, if $\alpha$ is $m$-torsion, then $d_{\alpha}(m\cdot \xi) =
d_{m\cdot \alpha}(\xi) = d_{0}(\xi) = 0 \in H^{3}_{an}(B,{\mathcal
O}_{B}^{\times})$ and so the image of $m\TSh_{an}(X)$ in
$H^{1}_{an}(B,\mycal{P}ic(X_{\alpha}/B))$ must be contained in
$Br'_{an}(X_{\alpha})$. This proves {\em (i)}.

Similarly, if $\alpha = m\cdot \varphi$ is $m$-divisible, then for any
$\xi$ we have $d_{\alpha}(\xi) = d_{m\cdot \varphi}(\xi) = d_{\varphi}
(m\cdot \xi)$ and so $d_{\alpha}$ vanishes identically on
$\TSh_{an}(X)[m]$. Thus the image of $\TSh_{an}(X)[m]$ in
$H^{1}_{an}(B,\mycal{P}ic(X_{\alpha}/B))$ must be contained in
$Br'_{an}(X_{\alpha})$ which completes the proof of {\em (ii)} and
the lemma.
\ \hfill
$\Box$

\

\bigskip

\noindent
We will denote the maps in items {\em (i)} and {\em (ii)} of
Lemma~\ref{lem-m-alpha} again by $T_{\alpha}$. Since by construction
these maps are compatible with the map $T_{\alpha}$ from
Lemma~\ref{lem-Talpha}, whenever the latter exists,
this abuse of notation can not lead to any confusion.

\

\bigskip

\noindent
Let us examine in more detail the map $d : H^{1}_{an}(B,\mycal{X}) \to
\op{Hom}_{{\mathbb Z}}(H^{1}_{an}(B,\mycal{X}),H^{3}_{an}(B,{\mathcal
O}_{B}^{\times}))$ given in \eqref{eq-d}. This map can be rewritten as
a bilinear pairing
\[
\langle \bullet, \bullet \rangle :
H^{1}_{an}(B,\mycal{X})\otimes_{{\mathbb Z}} H^{1}_{an}(B,\mycal{X})
\to H^{3}_{an}(B,{\mathcal O}_{B}^{\times}).
\]

\noindent The proof of Lemma~\ref{lem-m-alpha} shows that for
every $\alpha \in \TSh_{an}(X)$ we have a well defined
homomorphism
\[
T_{\alpha} : \alpha^{\perp} \to Br'_{an}(X_{\alpha}),
\]
where
$\alpha^{\perp} \subset \TSh_{an}(X)$ is the orthogonal complement of
$\alpha$ with respect to $\langle \bullet, \bullet \rangle$.

\begin{defi} \label{defi-orthogonal} Two genus one
fibrations $\alpha, \beta \in \TSh_{an}(X)$ will be called {\em
complementary} if $\langle \alpha, \beta \rangle = 0$. We will call
$\alpha$ and $\beta$ {\em $m$-compatible} if one of them is $m$-divisible
and the other one is $m$-torsion.
\end{defi}

\

Note that using the pairing
$\langle \bullet,\bullet \rangle,$
Lemma~\ref{lem-m-alpha} follows from the obvious observation that
every $m$-compatible pair $\alpha$, $\beta$ is
complementary.

For future reference we spell out the special case when $\alpha = 0$:

\begin{cor} \label{cor-T0} Assume that $X$ is smooth, the fibers of
$\pi$ are integral, and $Br'_{an}(B) = 0$. Then we have an isomorphism
$H^{1}_{an}(B,\mycal{P}ic(X/B)) \cong Br'_{an}(X)$ and we have an
exact sequence of abelian groups
\[
0 \to \TSh_{an}(X) \xrightarrow{T_{0}} Br'_{an}(X) \to H^{1}(B,{\mathbb Z}).
\]
\end{cor}
{\bf Proof.} Since $\sigma : B \to X$ is a section of $\pi$ it follows
that the composition
\[
H^{i}_{an}(B,{\mathcal O}_{B}^{\times}) \xrightarrow{\pi^{*}}
H^{i}_{an}(X,{\mathcal O}_{X}^{\times}) \xrightarrow{\sigma^{*}}
H^{i}_{an}(B,{\mathcal O}_{B}^{\times})
\]
is the identity. Thus $\ker\left( H^{i}_{an}(B,{\mathcal
O}_{B}^{\times}) \xrightarrow{\pi^{*}} H^{i}_{an}(X,{\mathcal
O}_{X}^{\times})\right) = 0$ and so $H^{1}_{an}(B,\mycal{P}ic(X/B))
\cong Br'_{an}(X)$. Combined with the fact that $\op{ind}(0) = 1$ this
gives the short exact sequence of groups above. The corollary is proven.
\ \hfill $\Box$

\medskip

Our pairing $\langle \bullet, \bullet \rangle$ can be explicitly 
described as follows. Every element $\alpha
\in \TSh_{an}(X) = H^{1}_{an}(B,\mycal{X})$ has two different
incarnations:
\begin{itemize}
\item $\alpha$ can be interpreted as a group extension of ${\mathbb
Z}_{B}$ by $\mycal{X}$. Concretely this is just the sheaf of groups
$\mycal{P}ic(X_{\alpha}/B)$ as it fits in the extension
\eqref{eq-XtoPic} viewed as an element $e(\alpha)$ in
$\op{Ext}_{{\mathbb Z}_{B}}^{1}({\mathbb Z}_{B},\mycal{X})$.
\item $\alpha$ can be interpreted as an extension of $\mycal{X}$ by
${\mathcal O}_{B}^{\times}[1]$. Concretely this is the amplitude one object
$_{\alpha}\mycal{X}$ in the derived category of abelian sheaves on $B$
which is the pullback of the extension class of
\[
1 \to {\mathcal O}_{B}^{\times}[1] \to R\pi_{\alpha*}{\mathcal
O}_{X_{\alpha}}^{\times}[1] \to R^{1}\pi_{\alpha*}{\mathcal
O}_{X_{\alpha}}^{\times} \to 1
\]
via the natural inclusion $\mycal{X} \to \mycal{P}ic(X_{\alpha}/B) =
R^{1}\pi_{\alpha*}{\mathcal O}_{X_{\alpha}}^{\times}$. Alternatively,
$_{\alpha}\mycal{X}$ can be thought of as a sheaf of commutative group
stacks on $B$ which is just the sheaf of all maps from $B$ to the
${\mathcal O}^{\times}$-gerbe on $X$ whose characteristic class is
$T_{0}(\alpha)$. Note that this gerbe is well defined in view of
Corollary~\ref{cor-T0}. We will write $g(\alpha) \in
\op{Ext}^{1}_{{\mathbb Z}_{B}}(\mycal{X},{\mathcal O}_{B}^{\times}[1])
= \op{Ext}^{2}_{{\mathbb Z}_{B}}(\mycal{X},{\mathcal O}_{B}^{\times})$
for the extension class of $_{\alpha}\mycal{X}$. For more on the
relevance of commutative group stacks see Section~\ref{dualstacks}.
\end{itemize}
With this notation it is now clear that $\langle \alpha, \beta \rangle$
is just the Yoneda product $g(\beta)\circ e(\alpha)$.

\begin{lem} The bilinear pairing
\[
\langle \bullet, \bullet \rangle :
H^{1}_{an}(B,\mycal{X})\bigotimes_{{\mathbb Z}} H^{1}_{an}(B,\mycal{X})
\to H^{3}_{an}(B,{\mathcal O}_{B}^{\times})
\]
is skew-symmetric. \label{lem-skewsymmetric}
\end{lem}
{\bf Proof.} The Poincare sheaf $\mycal{P} \to X\times_{B} X$
satisfies the biextension property and so can be interpreted functorially (see
\cite[Expos\'{e}~VII,Corollary~{3.6.5}]{sga7}) as an object
${\mathfrak L}({\mathcal P}) \in \op{ob} D^{b}({\mathbb
Z}_{B}\op{-mod})$ in the derived category
of abelian sheaves on $B$, which is an extension of
$\mycal{X}\overset{L}{\otimes}\mycal{X}$ by ${\mathcal
O}^{\times}_{B}$. In other words, ${\mathfrak L}({\mathcal P})$ fits
in a distinguished triangle
\[
{\mathcal O}_{B}^{\times} \to {\mathfrak L}({\mathcal P}) \to
\mycal{X}\overset{L}{\otimes} \mycal{X} \to {\mathcal O}_{B}^{\times}[1]
\]
of complexes of abelian sheaves. Let
$
{\mathfrak p} \in
\op{Ext}^{1}_{{\mathbb
Z}_{B}}(\mycal{X}\overset{L}{\otimes}\mycal{X},{\mathcal
O}_{B}^{\times}) =
\op{Hom}(\mycal{X}\overset{L}{\otimes}\mycal{X},{\mathcal
O}_{B}^{\times}[1])
$
be the corresponding extension class. From the definition of the
homomorphisms
\[
\begin{split}
e(\alpha) & \in \op{Hom}_{D^{b}({\mathbb
Z}_{B}\op{-mod})}({\mathbb Z}_{B},\mycal{X}[1]), \text{ and } \\
g(\alpha) & \in \op{Hom}_{D^{b}({\mathbb
Z}_{B}\op{-mod})}(\mycal{X},{\mathcal O}_{B}^{\times}[2])
\end{split}
\]
 one can
easily check that $g(\alpha)$ can be identified with the composition
\[
\mycal{X} = {\mathbb Z}_{B}\otimes \mycal{X}
\xrightarrow{e(\alpha)\otimes \op{id}_{\mycal{X}}}
\mycal{X}\overset{L}{\otimes}\mycal{X}[1] \xrightarrow{{\mathfrak p}}
{\mathcal O}_{B}^{\times}[1].
\]
Indeed, observe that both $g(\alpha)$ and ${\mathfrak p}\circ
(e(\alpha)\otimes \op{id}_{\mycal{X}})$ can naturally be interpreted as
amplitude one objects in the derived category of abelian sheaves on
$B$. Since any amplitude one object in $D^{b}({\mathbb
Z}_{B}\op{-mod})$ can be viewed as a stack over $B$ it suffices to
show the equivalence of the categories fibered in groupoids
corresponding to $g(\alpha)$ and ${\mathfrak p}\circ (e(\alpha)\otimes
\op{id}_{\mycal{X}})$ respectively.  Let as before $_{\alpha}\mycal{X}
\to B$ denote the fibered category corresponding to $g(\alpha)$. Since
by construction $_{\alpha}\mycal{X}$ comes from the push-forward
$R\pi_{\alpha*}{\mathcal O}_{X_{\alpha}}^{\times}$ we can identify
explicitly the groupoid of sections of $_{\alpha}\mycal{X}$ over an
open set $U$ in $B$ as the groupoid of all line bundles $L$ on
$(X_{\alpha}\times_{B} X)_{|U}$ having the property that for any point
$b \in U$ and any $x \in X_{b}$ we have that
$L_{|(X_{\alpha})_{b}\times \{x\}} \cong {\mathcal O}_{X_{b}}(x -
\sigma(b))$. Finally, using the description of the complex ${\mathfrak
L}(\mycal{P})$ in terms of fibered categories given in
\cite[Expos\`{e}~VII]{sga7} we see immediately that this groupoid is
precisely the groupoid of sections over $U$ of the fibered category
corresponding to ${\mathfrak p}\circ (e(\alpha)\circ \op{id}_{\mycal{X}})$.

Now taking into account that $g(\alpha) = {\mathfrak p}\circ
(e(\alpha)\otimes \op{id}_{\mycal{X}})$, we  see
 that for any two elements
$\alpha, \beta \in H^{1}_{an}(B,\mycal{X})$ the product $\langle
\alpha, \beta \rangle \in H^{3}_{an}(B,{\mathcal O}_{B}^{\times})$ can
be rewritten as the Yoneda product
\[
\langle \alpha, \beta \rangle = {\mathfrak p}\circ (\alpha\mcup \beta),
\]
where $\alpha\mcup\beta \in
H^{2}_{an}(B,\mycal{X}\overset{L}{\otimes} \mycal{X})$ is the external cup
product of $\alpha$ and $\beta$.

To understand the symmetry properties of $\langle \bullet,
\bullet \rangle$ it only remains to notice that
\[
\langle \beta, \alpha \rangle = {\mathfrak p}\circ (\beta\mcup
\alpha) = {\mathfrak p}\circ \sw(\alpha\mcup \beta),
\]
where $\sw : \mycal{X}\overset{L}{\otimes}\mycal{X} \to
\mycal{X}\overset{L}{\otimes}\mycal{X}$ is the involution switching the two
factors. However recall that  $\mycal{P}$ is a normalized Poincare
bundle and so can be explicitly described as
the rank one divisorial sheaf
\[
\mycal{P} = {\mathcal O}_{X\times_{B} X}(\Delta - \sigma\times_{B}X -
X\times_{B}\sigma - \varpi^{*}N_{\sigma/X})
\]
where $\varpi : X\times_{B} X \to B$ is the natural projection and
$N_{\sigma/X}$ is the normal bundle to the section $\sigma \subset X$.
In particular $\sw^{*}(\mycal{P}) = \mycal{P}$ and so $\mycal{P} \to
X\times_{B}X$ is a symmetric biextension. This shows that ${\mathfrak
p}\circ \sw = {\mathfrak p}$. Combined with the fact that
$\sw(\beta\mcup\alpha) = (-1)^{|\alpha|\cdot |\beta|}\alpha\mcup\beta
= -\alpha\mcup\beta$ we conclude that $\langle \beta, \alpha \rangle =
- \langle \alpha, \beta \rangle$. The lemma is proven. \ \hfill $\Box$

\

\bigskip

An immediate corollary of the skew-symmetry of
$\langle\bullet,\bullet\rangle$ is that $T_{\alpha}(\beta) \in
Br'_{an}(X_{\alpha})$  is well
defined iff  $T_{\beta}(\alpha) \in
Br'_{an}(X_{\beta})$  is well
defined.

In the case of surfaces we get the following:

\begin{cor} \label{cor-surface} Suppose that $X$ is a smooth surface
and that $\pi$ is non-isotrivial with all fibers integral. Then
$\TSh(X)$ is infinitely divisible
and so any $\alpha \in \TSh_{an}(X)$ is $m$-compatible with all
elements in $\TSh_{an}(X)[m]$.
\end{cor}
{\bf Proof.}
To show that $\TSh(X)$ is infinitely divisible note that since $B$ is
a curve we can apply Corollary~\ref{cor-T0} to conclude that the map
$T_{0}$ fits in a short exact sequence
\[
0 \to \TSh_{an}(X) \to Br'_{an}(X) \to H^{1}(B,{\mathbb Z})
\]
where the last map is the composition of the identification
$Br'_{an}(X) \cong H^{1}_{an}(B,\mycal{P}ic(X/B))$ coming from the
Leray spectral sequence and the map $H^{1}_{an}(B,\mycal{P}ic(X/B))
\to H^{1}(B,{\mathbb Z})$ corresponding to the degree morphism $\deg :
\mycal{P}ic(X/B) \to {\mathbb Z}_{B}$.

Since by assumption $\pi$ has only integral fibers, we have natural
identifications
\[
\mycal{P}ic(X/B) = R^{1}\pi_{*}{\mathcal
O}_{X}^{\times}, \quad \text{ and } \quad {\mathbb Z}_{B} =
R^{2}\pi_{*}{\mathbb Z}_{X}
\]
under which the degree map $\deg : \mycal{P}ic(X/B) \to {\mathbb
Z}_{B}$ becomes the coboundary homomorphism $\delta : R^{1}\pi_{*}{\mathcal
O}_{X}^{\times} \to R^{2}\pi_{*}{\mathbb Z}_{X}$ in the long exact
sequence of higher direct images associated to the exponential
sequence
\[
0 \to {\mathbb Z}_{X} \to {\mathcal O}_{X} \stackrel{\exp}{\to}
{\mathcal O}_{X}^{\times} \to 1
\]
and the map $\pi : X \to B$.

In particular, this implies that the map $Br'_{an}(X) \to
H^{1}(B,{\mathbb Z})$
fits in the commutative diagram:
\[
\xymatrix{
Br'_{an}(X) \ar[r] \ar[d]^-{\cong} &
H^{1}(B,{\mathbb Z}) \ar[d]_-{\cong} \\
H^{1}_{an}(B,R^{1}\pi_{*}{\mathcal O}^{\times}_{X}) \ar[r]^-{\delta}  &
H^{1}(B,R^{2}\pi_{*}{\mathbb Z}_{X}) \\
H^{2}_{an}(X,{\mathcal O}_{X}^{\times}) \ar[u]_-{\theta}
\ar[r]^-{\delta} & H^{3}(X,{\mathbb Z}). \ar[u]^-{\eta}
}
\]
Here the maps $\theta$ and $\eta$ between the third and second rows
come from the Leray 
spectral sequences for the map $\pi :  X \to B$ and the sheaves
${\mathcal O}^{\times}_{X}$ and ${\mathbb Z}_{X}$, which give:
\[
\xymatrix@1{
H^{2}_{an}(B,{\mathcal O}_{B}^{\times}) \ar[r] &
H^{2}_{an}(X,{\mathcal O}_{X}^{\times})  \ar[r]^-{\theta} &
H^{1}_{an}(B,R^{1}\pi_{*}{\mathcal O}_{X}^{\times}) \ar[r] & 0
}
\]
and
\[
\xymatrix@1{
H^{2}(B,R^{1}\pi_{*}{\mathbb Z}_{X}) \ar[r] & H^{3}(X,{\mathbb Z})
\ar[r]^-{\eta} & H^{1}(B,R^{2}\pi_{*}{\mathbb Z}_{X}) \ar[r] & 0.
}
\]
Now $H^{2}_{an}(B,{\mathcal O}_{B}^{\times}) = 0$ since $B$ is one
dimensional, and $H^{2}(B,R^{1}\pi_{*}{\mathbb Z}_{X})  = 0$ by the
irreducibility of monodromy.  This implies that
$\TSh_{an}(X) = \ker[Br'_{an}(X) \to H^{1}(B,{\mathbb Z})] =
\op{im}[H^{2}_{an}(X,{\mathcal O}_{X}) \to H^{2}_{an}(X,{\mathcal
O}_{X}^{\times})]$ and so $\TSh(X)$ is divisible. The corollary is
proven. \ \hfill $\Box$

\

\begin{defi} \label{defi-alphaXbeta} For any complementary pair $\alpha,
\beta \in \TSh_{an}(X)$ we denote the ${\mathcal O}^{\times}$-gerbe on
$X_{\beta}$
classified by $T_{\beta}(\alpha)$ by $\leftidx{_{\alpha}}{X}{}_{\beta}$.
\end{defi}

\begin{con} \label{con-duality} For any complementary pair $\alpha,
\beta \in \TSh_{an}(X)$, there exists an equivalence
\[
D^{b}_{1}(\leftidx{_{\alpha}}{X}{}_{\beta}) \cong D^{b}_{-1}(_{\beta}X_{\alpha})
\]
of the bounded derived categories of sheaves of pure weights $\pm 1$ on
$\Xc{\alpha}{\beta}$ and $\Xc{\beta}{\alpha}$ respectively.
\end{con}

In section~\ref{sec-duality} we will prove this conjecture in any
dimension under the additional assumptions that $\pi$ is smooth and
that $\alpha$ and $\beta$ are $m$-compatible. In section~\ref{sec-K3}
we will prove it unconditionally when $X$ is a surface.

\section{Smooth genus one fibrations} \label{sec-smooth}

In this section we will consider smooth genus one fibrations over
smooth bases of arbitrary dimension and ${\mathcal O}^{\times}$-gerbes over
them.

\subsection{${\mathcal O}^{\times}$-gerbes} \label{subsec-Oxgerbes}

In this section we work with a fixed smooth elliptic fibration
\[
\xymatrix@1{X
\ar[r]^-{\pi} & B \ar@/^0.5pc/[l]^-{\sigma}},
\]
and two genus one fibrations $X_{\alpha}, X_{\beta}$ corresponding to
two $m$-compatible elements $\alpha, \beta \in \TSh_{an}(X)$. Recall
that $m$-compatibility means that one of the elements, say $\beta$, is
actually algebraic, i.e. $\beta$ is a torsion element of some order $m$ in
$\TSh_{an}(X)$, while $\alpha$ is an $m$-divisible
element. Choose an element $\varphi \in
\TSh_{an}(X)$ such that $m\varphi = \alpha$.  We will use
this data to construct presentations for gerbes $\Xe{\beta}{\alpha}$ over
$X_{\alpha}$ and $\Xl{\alpha}{\beta}$ over
$X_{\beta}$. Different choices of the root $\varphi$ give rise to
different but Morita equivalent presentations of the same gerbes.

\subsubsection{The lifting presentation} \label{ssubsec-lifting}

Recall from Section~\ref{sec-TSh} that a gerbe presentation over a
variety $X$ is a
diagram
\[
\gp{R}{U}{X}{p_{1}}{p_{2}}{u}
\]
where $R \to U\times_{X} U$ is a ${\mathbb C}^{\times}$-bundle
satisfying the biextension condition.

\medskip

\noindent
We define a gerbe $\Xl{\alpha}{\beta}$ on $X_{\beta}$ via the {\em
lifting presentation}:
\begin{equation} \label{eq-lifting1}
\gp{\rl{\alpha}{\beta}}{\ul{\alpha}{\beta}}
{X_{\beta}}{p_{1}}{p_{2}}{p_{2}}
\end{equation}
where
\[
\begin{split}
\ul{\alpha}{\beta} & := X_{\varphi}\times_{B}
X_{\beta}, \\
p_{2} & : \ul{\alpha}{\beta} = X_{\varphi}\times_{B}
X_{\beta}
\to X_{\beta} \quad  \text{is the second projection, and } \\
\rl{\alpha}{\beta} & := \op{tot}(\mycal{P}^{\times}_{\ot}) \to
(\ul{\alpha}{\beta})\times_{X_{\beta}} (\ul{\alpha}{\beta}).
\end{split}
\]
Here $\mycal{P}_{\ot} \to (\ul{\alpha}{\beta})\times_{X_{\beta}}
(\ul{\alpha}{\beta})$ is the pullback  via the map
\[
\xymatrix@R-14pt@C-4pt{
p_{\ot} : & X_{\varphi}\times_{B} X_{\varphi} \times_{B} X_{\beta}
\ar[r] & X\times_{B} X \\
& (a,b,x) \ar@{|->}[r] & (a-b,m\cdot x),
}
\]
of the Poincare bundle
\[
\mycal{P} := {\mathcal O}_{X\times_{B} X}(\Delta - \sigma\times_{B} X
- X\times_{B} \sigma - \varpi^{*}c_{1}(B))
\]
on $X\times_{B} X$. As usual we denote the natural
projection $X\times_{B} X \to B$ by $\varpi$.
The required biextension property for $\mycal{P}_{\ot}$ follows immediately
from the see-saw principle .

For future reference we note that under the obvious identification
\[
(\ul{\alpha}{\beta})\times_{X_{\beta}} (\ul{\alpha}{\beta}) =
X_{\varphi}\times_{B} X_{\varphi} \times_{B} X_{\beta}
\]
the lifting presentation \eqref{eq-lifting1} can be rewritten as
\begin{equation} \label{eq-lifting2}
\xymatrix@R+10pt@C+15pt{
\op{tot}(\mycal{P}_{\ot}^{\times}) \ar[d] \ar@<.5ex>[drr]^-{s}
\ar@<-.5ex>[drr]_-{t} & & & \\
X_{\varphi}\times_{B} X_{\varphi} \times_{B} X_{\beta}
\ar@<.5ex>[rr]^-{p_{13}}  \ar@<-.5ex>[rr]_-{p_{23}} & &
X_{\varphi}\times_{B} X_{\beta}
\ar[r]^-{p_{2}} & X_{\beta}.
}
\end{equation}

\subsubsection{The extension presentation} \label{ssubsec-extension}

Similarly, we define a gerbe $\Xe{\beta}{\alpha}$ on $X_{\alpha}$ via
the {\em extension presentation}:
\begin{equation} \label{eq-extension1}
\gp{\re{\beta}{\alpha}}{\ue{\beta}{\alpha}}{X_{\alpha}}
{p_{1}}{p_{2}}{q_{\varphi}}
\end{equation}
where
\[
\begin{split}
\ue{\beta}{\alpha} & := X_{\varphi}, \\
q_{\varphi} & : X_{\varphi} \to X_{\alpha} \quad \text{is the
multiplication by $m$ map, and} \\
\re{\beta}{\alpha} & := \op{tot}(\Phi_{\beta}^{\times}) \to
X_{\varphi}\times_{X_{\alpha}} X_{\varphi}.
\end{split}
\]
Here $\Phi_{\beta}$ could be taken as the line bundle
$d^{*}(M)$ on
$X_{\varphi}\times_{X_{\alpha}} X_{\varphi}$, where
\[
d : X_{\varphi}\times_{X_{\alpha}} X_{\varphi} \to X[m]
\]
is the difference map and $M$ is any line bundle on $X[m]$ whose
punctured total space gives a group extension:
\[
\xymatrix@1{
1 \ar[r] & {\mathcal O}^{\times}_{B}  \ar[r] & \op{tot}(M^{\times})
\ar[r] & X[m]
\ar[r] & 0.
}
\]
We will explain below how to construct a relative line bundle
$\Sigma_{\beta} \in \Gamma(B,\mycal{P}ic^{m}(X_{\beta}/B))$ and a
global line bundle $M_{\beta} \to X[m]$, determined by the condition
that its punctured total space is the theta group $G_{\beta}$:
\[
\xymatrix@1{
1 \ar[r] & {\mathcal O}^{\times}_{B}  \ar[r] & G_{\beta} \ar[r] & X[m]
\ar[r] & 0 }
\]
of $\Sigma_{\beta}$.
The simplest choice would be to take $M := M_{\beta}$.
However, we will see later that in order to achieve duality with the
lifting gerbe, the correct choice is to take
$M:=M_{\beta} \otimes M_0^{-1}$, where
$M_{0}$ is determined by the condition that its
punctured total space is the
theta group $G_{0}$:
\[
\xymatrix@1{
1 \ar[r] & {\mathcal O}^{\times}_{B}  \ar[r] & G_{0} \ar[r] & X[m]
\ar[r] & 0}
\]
corresponding to the similarly defined relative line bundle
$\Sigma_{0} \in \Gamma(B,\mycal{P}ic^{m}(X/B))$.

To define $\Sigma_{\beta}$, consider first two genus one curves $E'$
and $E''$ with the same Jacobian $E$. Let $q : E' \to E''$ be a map
which induces the multiplication by $m$ map $E \to E$. For any two
points $a, b \in E'$ such that $q(a) = q(b)$, we have that ${\mathcal
O}_{E'}(m\cdot a) \cong {\mathcal O}_{E'}(m\cdot b)$. This determines
a map $E'' \to \op{Pic}^{m}(E')$. Applying this to our map $q_{\beta}
: X_{\beta} \to X$ we get a well defined morphism $X \to
\op{Pic}^{m}(X_{\beta}/B)$. The image of $\sigma \subset X$ is our
relative line bundle $\Sigma_{\beta}$. Similarly we get $\Sigma_{0}$
from the multiplication by $m$ map $q_{0}^{m} : X \to X$. Note that by
construction the relative line bundle $\Sigma_{0} \in
\op{Pic}^{m}(X/B)$ is induced by a global line bundle ${\mathcal
O}_{X}(m\sigma)$. Similarly, if $Br_{an}'(B) = 0$, then the relative
line bundle $\Sigma_{\beta} \in \op{Pic}^{m}(X/B)$ is induced by some
global line bundle ${\mathcal O}_{X_{\beta}}(\sigma_{\beta})$, where
$\sigma_{\beta} \subset X_{\beta}$ is an $m$-section of $\pi_{\beta}$.

Given any relative line bundle $\Sigma \in
\Gamma(B,\op{Pic}^{m}(Y/B))$, where $\pi_{Y} : Y \to B$ is in
$\TSh(X)$, we get a global line bundle $M$ on $X[m]$ equipped with an
isomorphism
\[
\mu^{*}M \cong \op{pr}_{1}^{*}M\otimes \op{pr}_{2}^{*}M
\]
satisfying the usual biextension property. Here $\mu : X[m]\times_{B}
X[m] \to X[m]$ is multiplication, and $\op{pr}_{1}$, $\op{pr}_{2}$ are the
projections.

Locally on $B$, choose an actual line bundle $\Sigma' \in \op{Pic}(Y)$
lifting $\Sigma$. We define $M' := p_{1*}(a^{*}\Sigma'\otimes
p_{2}^{*}\Sigma^{'-1})$. Here $a : X[m]\times_{B} Y \to Y$ is the
action and $p_{1}$, $p_{2}$ are the projections. Note that
$a^{*}\Sigma' \otimes p_{2}^{*}\Sigma^{'-1}$   is trivial on each
fiber of $p_{1}$, so $M'$ is a line bundle. Let $M''$ be the line
bundle determined by another lift $\Sigma''$ of $\Sigma$. Any such
$\Sigma''$  is necessarily of the form $\Sigma'\otimes \pi_{Y}^{*}L$
for some line bundle $L$ on $B$. The equality $\pi_{Y}\circ a =
\pi_{Y}\circ p_{2}$ therefore induces a canonical isomorphism
\[
a^{*}\Sigma' \otimes p_{2}^{*}\Sigma^{'-1} \to  a^{*}\Sigma'' \otimes
p_{2}^{*}\Sigma^{''-1},
\]
so the locally defined line bundles $M'$, $M''$ glue to a global line
bundle $M$ on $X[m]$.

The biextension property of $M$ is equivalent to saying that its
punctured total space $G := \op{tot}(M^{\times})$ has the structure of a
group scheme over $B$ called the {\em theta group} of $\Sigma$. It is
a central extension of $X[m]$ by ${\mathbb G}_{m}$. Explicitly a local
section of $G$ is a pair $(x,\lambda)$, where $x$ is a local section
of $X[m] \to B$ and $\lambda : t_{x}^{*}\Sigma \to \Sigma$ is an
isomorphism.

Applying these constructions to our relative line bundles
$\Sigma_{\beta}$ and $\Sigma_{0}$ produces the desired line bundles
$M_{\beta}$ and $M_{0}$ and theta groups $G_{\beta}$ and $G_{0}$.

For future reference we note that under the obvious isomorphism
\[
d\times p_{2}: X_{\varphi}\times_{X_{\alpha}} X_{\varphi} \widetilde{\to}
X[m]\times_{B} X_{\varphi}
\]
the extension presentation \eqref{eq-extension1} can be rewritten in
the equivalent form:
\begin{equation} \label{eq-extension2}
\xymatrix@R+10pt@C+15pt{
\op{tot}(\Phi_{\beta}^{\times}) \ar[d] \ar@<.5ex>[drr]^-{s}
\ar@<-.5ex>[drr]_-{t} & & & \\
X[m]\times_{B} X_{\varphi} \ar@<.5ex>[rr]^-{a_{\varphi}}
\ar@<-.5ex>[rr]_-{p_{2}} & & X_{\varphi}
\ar[r]^-{q_{\varphi}} & X_{\alpha},
}
\end{equation}
where $a_{\varphi} : X[m]\times_{B} X_{\varphi} \to X_{\varphi}$
denotes the action and $p_{2}$ denotes the second projection.

\subsubsection{Coboundary realizations} \label{ssubsec-coboundary}

Although we constructed $\Xl{\alpha}{\beta}$ and $\Xe{\beta}{\alpha}$
directly, it is worth noting that the lifting and extension
presentations are both special cases of the coboundary construction
described in Sections \ref{ssubsec-gerbes} and
\ref{ssubsec-presentations}. Recall that the input for the coboundary
construction for a gerbe presentation on a variety $Y$ consists of a
short exact sequence of group schemes over $Y$
\[
1 \to {\mathbb G}_m \to G \to K \to 1
\]
together with a $K$ torsor $U$.

The lifting presentation is obtained from the short exact sequence
\[
1 \to {\mathbb G}_m \to \op{tot}(p_{1,m \cdot 2}^*
\mycal{P}^{\times}) \to \pi_{\beta}^*(X) \to 1
\]
and the
$\pi_{\beta}^*(X) = X \times_B X_{\beta}$-torsor
$X_{\varphi} \times_B X_{\beta}$. Note that
the group structure on
$\op{tot}(p_{1,m \cdot 2}^* \mycal{P}^{\times})$ in the above sequence
comes from the biextension property
of the Poincare bundle (for the group law on $\pi_{\beta}^*(X)$).

The extension presentation is obtained from the short exact sequence
\[
1 \to {\mathbb G}_m \to \op{tot}(\Phi_{\beta}^{\times}) \to
\pi_{\alpha}^*(X[m]) \to 1
\]
and the $\pi_{\alpha}^*(X[m])= X[m] \times_B X_{\alpha}$-torsor
$X_{\varphi}$.

In other words, we can write $\Xl{\alpha}{\beta}$ and  $\Xe{\beta}{\alpha}$
as quotient gerbes:

\[
\begin{split}
\Xl{\alpha}{\beta} & = [X_{\varphi}\times_B X_{\beta}/\op{tot}(p_{1,m
\cdot 2}^* \mycal{P}^{\times})] \\
\Xe{\beta}{\alpha} & = [X_{\varphi}/\op{tot}(\Phi_{\beta}^{\times})].
\end{split}
\]

\subsection{The class of the lifting gerbe}

In this section we continue to assume that $\pi : X \to B$ is smooth, $\beta \in
\TSh_{an}(X)$ is of finite order $m$ and $\alpha \in
\TSh_{an}(X)$ is $m$-divisible.

\begin{theo} \label{theo-lifting-vs-coho} The class
$[\Xl{\alpha}{\beta}]$ of the lifting gerbe equals
$T_{\beta}(\alpha)$. In other words $\Xl{\alpha}{\beta}$ is a model
for $\Xc{\alpha}{\beta}$.
\end{theo}
{\bf Proof.} The proof is in two steps. In step (1) we show that a
cocycle representing the class $T_{\beta}(\alpha)$ which defines
$\Xc{\alpha}{\beta}$ becomes a
coboundary $\delta(c)$ when pulled back to $L :=
\ul{\alpha}{\beta}$. In step (2) we
check that the line bundle defined by $c$ on $L\times_{X_{\beta}}
L$ coincides with the Poincare bundle $\rl{\alpha}{\beta}$.

We need to show the isomorphism of two gerbes on the smooth space
$X_{\beta}$.

We will be working with the Cartesian product
\[
\xymatrix{L \ar[r]^-{\lambda_{\varphi}} \ar[d]_-{\lambda_{\beta}}
\ar[dr]^{\lambda_{B}} & X_{\varphi} \ar[d]^-{\pi_{\varphi}} \\
X_{\beta} \ar[r]_-{\pi_{\beta}} & B
}
\]
Recall from section~\ref{subsec-Brauer} that the class of
$\Xc{\alpha}{\beta}$ is $T_{\beta}(\alpha) \in
H^{2}(X_{\beta},{\mathcal O}^{\times})$. By the Leray spectral
sequence for $\pi_{\beta} : X_{\beta} \to B$, this group equals
$H^{1}(B,\mycal{P}ic(X_{\beta}/B))$. Explicitly, $T_{\beta}(\alpha)$ is
the class of the $\mycal{P}ic(X_{\beta}/B)$-torsor induced from the
$\mycal{P}ic^{0}(X_{\beta}/B) = X$-torsor $X_{\beta}$.

Similarly, the Leray spectral sequence for $\lambda_{B} : L \to B$
gives an injection
\begin{equation}
\label{eq-injection}
H^{1}(B,\mycal{P}ic(L/B)) \hookrightarrow H^{2}(L,{\mathcal O}^{\times}).
\end{equation}
The $\lambda_{\beta}$-pullback of $T_{\beta}(\alpha)$ is in the image
of \eqref{eq-injection} and is the $\mycal{P}ic(L/B)$-torsor induced
from $X_{\beta}$ via $\lambda_{\beta}^{*} :
\mycal{P}ic^{0}(X_{\beta}/B) \to \mycal{P}ic(L/B)$.

For step (1), consider the short exact sequence of sheaves of groups
on $B$:
\[
0 \to \mycal{P}ic(X_{\beta}/B) \stackrel{\lambda_{\beta}^{*}}{\to}
\mycal{P}ic(L/B) \stackrel{\op{ev}}{\to}
\mycal{H}om_{B}(X_{\beta},\op{Pic}(X_{\varphi}/B)) \to 0.
\]
Here $\op{ev}$ sends a line bundle on $L$ to the family of its
restrictions on $\op{pt}\times_{B} X_{\varphi}$, and it is surjective
because $\pi_{\beta} : X_{\beta} \to B$ is smooth.

We are claiming that $\lambda_{\beta}^{*}(T_{\beta}(\alpha)) = 0$, so
we need to show that $T_{\beta}(\alpha)$ is in the image of the
coboundary
\[
\partial : H^{0}(B,\mycal{H}om_{B}(X_{\beta},\op{Pic}(X_{\varphi}/B)))
\to H^{1}(B,\mycal{P}ic(X_{\beta}/B)).
\]
In $H^{0}(B,\mycal{H}om_{B}(X_{\beta},\op{Pic}(X_{\varphi}/B)))$ we
have a natural element
\[
q : X_{\beta} \to X_{m\beta} = X = \op{Pic}^{0}(X_{\varphi}/B)\subset
\op{Pic}(X_{\varphi}/B),
\]
depending on the choice of a trivialization $\Sigma_{\beta}$ of
$X_{m\beta}$.

We will see that $\partial(q) = T_{\beta}(\alpha)$.

Let ${\mathfrak U} = \{ U_{i} \}_{i \in I}$ be an analytic open cover
of $B$ for which we have trivializations $s_{i} : U_{i} \to
X_{\varphi}$ of the $X$-torsor $X_{\varphi}$. In order to calculate
$\partial(q)$ we first lift $q$ to an element $c \in
C^{0}({\mathfrak U},\mycal{P}ic(L/B))$. This lift is given in terms of
the map
\[
t_{-s_{i}}\times q : L_{|U_{i}} = X_{\varphi}\times_{U_{i}} X_{\beta} \to
X\times_{U_{i}} X
\]
by $c_{i} := (t_{-s_{i}}\times q)^{*}\mycal{P}$, where
$\mycal{P}$ is the standard Poincare bundle on $X\times_{B} X$.

The \v{C}ech differential $\delta(c) \in Z^{1}({\mathfrak
U},\mycal{P}ic(L/B))$ is given by $\{ c_{i}\otimes
c_{j}^{-1} \}_{i,j \in I}$. It comes from $Z^{1}({\mathfrak
U},\mycal{P}ic(X_{\beta}/B))$, and is represented there by $\{
{\mathcal O}_{\pi_{\beta}^{-1}(U_{ij})}(m(s_{j} -s_{i})) \}_{i,j \in
I}$. On the other hand, $ms_{i}$ can be interpreted as a section  of
$X_{\alpha} = X_{m\varphi}$ over
$U_{i}$, so this cocycle represents our
$T_{\beta}(\alpha)$.

For step (2), consider the cochain
\[
\{ p_{1}^{*}c_{i}\otimes p_{2}^{*}c_{i}^{-1} \}_{i \in I}
\in C^{0}({\mathfrak U},\mycal{P}ic(L\times_{X_{\beta}}L/B)).
\]
>From the discussion in section~\ref{subsec-Oxgerbes}  we know that
this cochain is in fact a global section ${\mathbb L}$ of
$\mycal{P}ic(L\times_{X_{\beta}}L/B)$. We need to show that
${\mathbb L} = \mycal{P}_{\ot}$.
As usual we identify $L\times_{X_{\beta}} L$ with
$X_{\varphi}\times_{B} X_{\varphi}\times_{B} X_{\beta}$, so $p_{1}$ and
$p_{2}$ become $p_{13}$ and $p_{23}$.

It suffices to show the equality ${\mathbb
L}_{|\lambda_{B}^{-1}(U_{i})} = p_{1}^{*}c_{i}\otimes
p_{2}^{*}c_{i}^{-1}$ for each open set $U_{i}$. This follows by
the theorem of the cube from the identifications:
\[
\begin{split}
p_{1}^{*}c_{i} & = p_{13}^{*}\circ (t_{-s_{i}}\times q)^{*}
\mycal{P} \\
p_{2}^{*}c_{i} & = p_{23}^{*}\circ (t_{-s_{i}}\times q)^{*}
\mycal{P} \\
{\mathbb L} & = p_{\ot}^{*}\mycal{P}.
\end{split}
\]
This finishes the proof of the theorem. \hfill $\Box$

\subsection{The class of the extension gerbe}

In this section we again assume that $\pi : X \to B$ is a smooth
elliptic fibration, that $Br_{an}'(B) = 0$ and that $\alpha, \beta \in
\TSh_{an}(X)$ are $m$-compatible with $\beta$ of order $m$.

\begin{theo} \label{theo-extension-vs-coho} The class
$[\Xe{\beta}{\alpha}]$ of the extension gerbe equals
$T_{\alpha}(\beta)$. In other words $\Xe{\beta}{\alpha}$ is a model
for $\Xc{\beta}{\alpha}$.
\end{theo}
{\bf Proof.} Recall that the assumption $Br_{an}'(B) = 0$ together
with the Leray spectral sequence for $\pi_{\alpha} : X_{\alpha} \to B$
give us an injection
\[
H^{2}_{an}(X_{\alpha},{\mathcal O}^{\times}) \hookrightarrow
H^{1}_{an}(B, R^{1}\pi_{\alpha*}{\mathcal O}^{\times}) =
H^{1}(B,\mycal{P}ic(X_{\alpha}/B)).
\]
In terms of this inclusion, $T_{\alpha}(\beta)$ can be identified with
the isomorphism class of the \linebreak $\mycal{P}ic(X_{\alpha}/B)$-torsor
associated to the $\mycal{P}ic^{0}(X_{\alpha}/B) = \mycal{X}$ torsor
$X_{\beta}$. In order to show that $[\Xe{\beta}{\alpha}] =
T_{\alpha}(\beta)$  we must first check that $T_{\alpha}(\beta)$ pulls
back to the trivial element in $H^{2}_{an}(\ue{\beta}{\alpha},
{\mathcal O}^{\times})$.

Recall from Section~\ref{ssubsec-extension} that the atlas
$\ue{\beta}{\alpha}$ for the extension presentation is defined by
fixing an element $\varphi \in \TSh_{an}(X)$ such that $m\cdot
\varphi =
\alpha$ and then taking $\ue{\beta}{\alpha} := X_{\varphi}$. With this
definition the
structure morphism $\ue{\beta}{\alpha} \to X_{\alpha}$ is identified
with the map $q:=q_{\varphi}^{m} : X_{\varphi} \to X_{\alpha}$ of
multiplication by $m$ along the fibers.

The pullback via $q$ of relative line bundles defined on the fibers of
$\pi_{\alpha} : X_{\alpha} \to B$  gives rise to a morphism of sheaves
of groups
\[
Q: \mycal{P}ic(X_{\alpha}/B) \longrightarrow \mycal{P}ic(X_{\varphi}/B).
\]
Since $q$ corresponds to multiplication by $m$, it follows that $Q$
fits in a commutative diagram
\[
\xymatrix@M+3pt{
\mycal{P}ic(X_{\alpha}/B) \ar[r]^-{Q} & \mycal{P}ic(X_{\varphi}/B) \\
\mycal{P}ic^{0}(X_{\alpha}/B) \ar@{^{(}->}[u] \ar@{=}[d] &
\mycal{P}ic^{0}(X_{\varphi}/B) \ar@{^{(}->}[u] \ar@{=}[d] \\
\mycal{X} \ar[r]^-{\op{mult}_{m}} & \mycal{X}
}
\]
Also, since $\pi_{\alpha}\circ q = \pi_{\varphi}$, it follows that the
pullback map
\[
q^{*} : H^{2}_{an}(X_{\alpha},{\mathcal O}^{\times}) \longrightarrow
H^{2}_{an}(X_{\varphi},{\mathcal O}^{\times})
\]
is compatible with the Leray spectral sequences for $\pi_{\alpha}$ and
$\pi_{\varphi}$, and so fits in a commutative diagram
\[
\xymatrix@M+3pt{H^{2}_{an}(X_{\alpha},{\mathcal O}^{\times})
\ar[r]^-{q^{*}} \ar@{^{(}->}[d] & H^{2}_{an}(X_{\varphi},{\mathcal
O}^{\times})\ar@{^{(}->}[d] \\ H^{1}_{an}(B,\mycal{P}ic(X_{\alpha}/B))
\ar[r]^{h^{1}(Q)} & H^{1}_{an}(B,\mycal{P}ic(X_{\varphi}/B)) \\
H^{1}(B,\mycal{X}) \ar[u] \ar[r]^-{h^{1}(\op{mult}_{m})} &
H^{1}(B,\mycal{X}) \ar[u] }
\]
Thus we can identify $q^{*}(T_{\alpha}(\beta))$ with the class of the
$\mycal{P}ic(X_{\varphi}/B)$-torsor which is induced from the
$\mycal{P}ic^{0}(X_{\varphi}/B) = \mycal{X}$-torsor
$h^{1}(\op{mult}_{m})(X_{\beta}) = X_{m\cdot \beta}$. However by
assumption $m\cdot \beta = 0$ and so $X_{m\cdot \beta}$ is trivial as
a $\mycal{X}$-torsor. Therefore $q^{*}(T_{\alpha}(\beta)) = 0$ as
promised.

To complete the proof of the theorem we need to realize the cocycle
$q^{*}(T_{\alpha}(\beta))$ as a coboundary:
\[
q^{*}(T_{\alpha}(\beta)) = \partial(\psi)
\]
for some $\psi \in C^{1}_{an}(X_{\varphi},{\mathcal O}^{\times})$,
and then check that the line bundle defined by $\psi$ on
$X_{\varphi}\times _{X_{\alpha}} X_{\varphi}$ is isomorphic to
$\Phi_{\beta}$.

In terms of the inclusion $H^{2}_{an}(X_{\varphi},{\mathcal
O}^{\times}) \subset H^{1}_{an}(B,\mycal{P}ic(X_{\varphi}/B))$ this
amounts to writing the class $q^{*}(T_{\alpha}(\beta)) \in
H^{1}_{an}(B,\mycal{P}ic(X_{\varphi}/B))$ as the coboundary of some
\v{C}ech cochain $\psi \in C^{0}_{an}(b,\mycal{P}ic(X_{\varphi}/B))$
and then showing that the global section of
$\mycal{P}ic(X_{\varphi}\times_{X_{\alpha}} X_{\varphi}/B)$ determined
by $\psi$ coincides with the global section given by
$\Phi_{\beta}$.
To carry this out we will need to first choose a  cocycle
representating of $T_{\beta}(\alpha) \in
H^{1}_{an}(B,\mycal{P}ic(X_{\alpha}/B))$ or equivalently a cocycle
representating for the $\mycal{X}$-torsor $X_{\beta}$.

Let ${\mathfrak U} = \{ U_{i} \}$ be an analytic open covering of $B$
which trivializes $X_{\beta}$ as an $X$ torsor. Choose trivializing
sections $s_{i} \in \Gamma(U_{i},X_{\beta})$ over each $U_{i}$. Then
$T_{\alpha}(\beta) \in H^{1}_{an}(B,\mycal{P}ic(X_{\alpha}/B))$ is
represented by the \v{C}ech cocycle
\[
\{{\mathcal O}_{X_{\beta}}(s_{j} - s_{i}) \} \in Z^{1}({\mathfrak
U},\mycal{X}) = Z^{1}({\mathfrak U}, \mycal{P}ic^{0}(X_{\alpha}/B))
\to
Z^{1}({\mathfrak U}, \mycal{P}ic(X_{\alpha}/B)).
\]
Here ${\mathcal O}_{X_{\beta}}(s_{j} - s_{i})$ is viewed as a line
bundle of degree zero along the fibers of \linebreak $\pi_{\alpha} :
X_{\alpha|U_{ij}} \to U_{ij}$ via the canonical identification
$\mycal{P}ic^{0}(X_{\beta}/B)= \mycal{X} =
\mycal{P}ic^{0}(X_{\alpha}/B)$. In particular
$q^{*}(T_{\alpha}(\beta)) \in
H^{1}_{an}(B,\mycal{P}ic(X_{\varphi}/B))$ is represented by the
cocycle
\[
\{{\mathcal O}_{X_{\beta}}(s_{j} - s_{i})^{\otimes m} \} \in
Z^{1}({\mathfrak U},\mycal{X}) = Z^{1}({\mathfrak U},
\mycal{P}ic^{0}(X_{\varphi}/B)) \to Z^{1}({\mathfrak U},
\mycal{P}ic(X_{\varphi}/B)).
\]
In order to write this cocycle as a coboundary we will have to
trivialize the $\mycal{X}$-torsor $X_{m\cdot \beta}$.

Recall that in the construction of the line bundle $\Phi_{\beta}$ we
used a particular trivialization of $X_{m\cdot \beta}$, namely the
relative line bundle
\[
\Sigma_{\beta} \in \Gamma_{an}(B,\mycal{P}ic^{m}(X_{\beta}/B)).
\]
Using $\Sigma_{\beta}$ we can construct a cochain
\[
\psi = \{ \psi_{i} \} \in C^{0}({\mathfrak
U},\mycal{P}ic^{0}(X_{\varphi}/B))
\]
with $\psi_{i} := {\mathcal O}_{X_{\beta}}(-m\cdot s_{i})\otimes
\Sigma_{\beta}$. By construction, $\psi$ determines a global section
\[
p_{1}^{*}\psi\otimes p_{2}^{*}\psi^{-1} : B \to
\op{Pic}^{0}(X_{\varphi}\times_{X_{\alpha}} X_{\varphi}/B),
\]
namely, the section locally given by $p_{1}^{*}\psi_{i}\otimes
p_{2}^{*}\psi^{-1}_{i} \in \Gamma_{an}(U_{i},
\mycal{P}ic^{0}(X_{\varphi}\times_{X_{\alpha}} X_{\varphi}/B))$.
On the other hand the section corresponding to $\Phi_{\beta} \to
X_{\varphi}\times_{X_{\alpha}} X_{\varphi}$ can be described as
follows.

Recall from section \ref{ssubsec-extension} that
\[
\Phi_{\beta}  = d^{*}(M_{\beta}\otimes M_{0}^{-1})
\]
where $d : X_{\varphi}\times_{X_{\alpha}} X_{\varphi} \to X[m]$ is the
natural difference map and $M_{\beta}$, $M_{0}$ are line bundles on
$X[m]$ satisfying
\[
\begin{split}
d_{\beta}^{*}M_{\beta} & \cong p_{1}^{*}\Sigma_{\beta}\otimes
p_{2}^{*}\Sigma_{\beta}^{-1} \\
d_{0}^{*}M_{0} & \cong p_{1}^{*}\Sigma_{0}\otimes
p_{2}^{*}\Sigma_{0}^{-1}.
\end{split}
\]
Here again
\[
\xymatrix@R-2pc{ X_{\beta}\underset{q_{\beta}^{m}, X, q_{\beta}^{m}}{\times}
X_{\beta} \ar[r]^-{d_{\beta}} & X[m] \\ X\underset{q^{m}_{0}, X,
q^{m}_{0}}{\times} X \ar[r]^-{d_{0}} & X[m] }
\]
stand for the difference maps.

By construction $p_{1}^{*}\psi_{i}\otimes p_{2}^{*}\psi_{i}^{-1}$
lives naturally in $\Gamma(U_{i}, \mycal{P}ic^{0}(X_{\beta}\times_{B}
X_{\beta}/B))$ (which we have identified with $\Gamma(U_{i},
\mycal{P}ic^{0}(X_{\varphi}\times_{B} X_{\varphi}/B))$).  In view of
this, it will be convenient if we rewrite all objects as line bundles
on $X_{\beta}\times_{X} X_{\beta}$. To that end, choose a local
section $s : U_{i} \to X_{\beta}$ and let $t_{-s}$ be the induced
isomorphism
\[
\xymatrix{
X_{\beta|U_{i}} \ar[rr]^-{t_{-s}}\ar[rd] & & X_{|U_{i}} \ar[ld] \\
& U_{i} &
}
\]
of translation by $s$ along the fibers. With this notation we have a
commutative diagram
\[
\xymatrix{ X[m]_{|U_{i}} \ar@{=}[d] & (X_{\beta}\times_{X}
X_{\beta})_{|U_{i}} \ar[l]_-{d_{\beta}} \ar[d]^-{t_{-s}\times t_{-s}}
\ar@<.5ex>[r]^-{p_{1}} \ar@<-.5ex>[r]_-{p_{2}} & X_{\beta_|U_{i}}
\ar[d]^-{t_{-s}} \\
X[m]_{|U_{i}} & (X\times_{X} X)_{|U_{i}} \ar[l]_-{d_{0}}
\ar@<.5ex>[r]^-{p_{1}} \ar@<-.5ex>[r]_-{p_{2}} & X_{|U_{i}} }
\]
and thus
\[
d_{\beta}^{*} M_{0} = (t_{s}\times t_{-s})^{*}\circ d_{0}^{*}
M_{0} = p_{1}^{*}{\mathcal O}_{X_{\beta}}(m s)\otimes
p_{2}^{*}{\mathcal O}_{X_{\beta}}(-m s).
\]
Therefore, in order to compare $p_{1}^{*}\psi_{i}\otimes
p_{2}^{*}\psi_{i}^{-1}$ and $d_{\beta}^{*}(M_{\beta}\otimes
M_{0}^{-1})$, we need to show that on $(X_{\beta}\times_{X}
X_{\beta})_{|U_{i}}$ we have
\[
p_{1}^{*}(\Sigma_{\beta}(-m s_{i}))\otimes p_{2}^{*}(\Sigma_{\beta}(m
s_{i})) \cong p_{1}^{*}(\Sigma_{\beta}(-m s))\otimes
p_{2}^{*}(\Sigma_{\beta}(m s))
\]
for every section $s : U_{i} \to X_{\beta}$.

Equivalently, it suffices to show that
\begin{equation} \label{eq-pullback-identity}
p_{1}^{*}{\mathcal O}_{X_{\beta}}(m(s - s_{i}))\otimes
p_{2}^{*}{\mathcal O}_{X_{\beta}}(m(s_{i} - s)) \cong {\mathcal O}.
\end{equation}
On the other hand we have a commutative diagram
\[
\xymatrix{
& X_{\beta}\times_{X} X_{\beta} \ar[dl]_-{p_{1}} \ar[dr]^-{p_{2}}
\ar[dd]^-{\mu} & \\
X_{\beta} \ar[dr]_-{q_{\beta}^{m}} & & X_{\beta}
\ar[dl]^-{q_{\beta}^{m}} \\
& X &
}
\]
and since ${\mathcal O}_{X_{\beta}}(m(s - s_{i}))$ is of degree zero
along the fibers of $\pi_{\beta}$ we have
\begin{equation} \label{eq-pullback}
{\mathcal O}_{X_{\beta}}(m(s - s_{i})) = (q_{\beta}^{m})^{*} {\mathcal
O}_{X}(s - s_{i}).
\end{equation}
Here ${\mathcal O}_{X}(s - s_{i}) \in \Gamma(U_{i},
\mycal{P}ic^{0}(X/B))$ denotes the relative line bundle on $X$
corresponding to $ {\mathcal O}_{X_{\beta}}(s - s_{i}) \in
\Gamma(U_{i}, \mycal{P}ic^{0}(X_{\beta}/B))$ under the canonical
identification $\mycal{P}ic^{0}(X/B)) = \mycal{X} =
\mycal{P}ic^{0}(X_{\beta}/B)$.

The formula \eqref{eq-pullback} implies that
\[
\begin{split}
p_{1}^{*}{\mathcal O}_{X_{\beta}}(m(s - s_{i})) & \cong
\mu^{*}{\mathcal O}_{X}(s - s_{i}) \\
p_{2}^{*}{\mathcal O}_{X_{\beta}}(m(s_{i} - s)) & \cong
\mu^{*}{\mathcal O}_{X}(s_{i} - s)
\end{split}
\]
and so \eqref{eq-pullback-identity} holds. The theorem is
proven. \hfill $\Box$

\subsection{Duality between the lifting and extension presentations}
\label{sec-duality}

We are now ready to prove Theorem\ref{Main-smooth} for a smooth
elliptic fibration
\[
\xymatrix@1{X
\ar[r]^-{\pi} & B \ar@/^0.5pc/[l]^-{\sigma}},
\]
over a smooth space $B$ satisfying $Br'_{an}(B) = 0$.

Let $\alpha, \beta, \varphi \in \TSh_{an}(X)$ satisfy $m
\beta = 0$, $m \varphi = \alpha$ (in particular $\alpha$ and $\beta$
are $m$-compatible). We want to compare the derived
categories of coherent sheaves on $\Xc{\alpha}{\beta}$ and
$\Xc{\beta}{\alpha}$.

\subsubsection{The gerby Fourier-Mukai transform} \label{ssubsec-FM}

Let $D^{b}_{1}(\Xc{\alpha}{\beta})$ and
$D^{b}_{-1}(\Xc{\beta}{\alpha})$ denote the derived categories of
coherent sheaves of weight one and minus one on the gerbes
$\Xc{\alpha}{\beta}$ and $\Xc{\beta}{\alpha}$
respectively. Alternatively, as explained at the end of
Section~\ref{ssubsec-presentations}, we can view
$D^{b}_{1}(\Xc{\alpha}{\beta})$ and $D^{b}_{-1}(\Xc{\beta}{\alpha})$
as derived categories of $T_{\beta}(\alpha)$-twisted sheaves on
$X_{\beta}$ and $T_{\alpha}(-\beta)$ twisted sheaves on $X_{\alpha}$
respectively.

We want to construct a Fourier-Mukai functor
\[
\FM : D^{b}_{1}(\Xc{\alpha}{\beta}) \to D^{b}_{-1}(\Xc{\beta}{\alpha})
\]
which is an equivalence. To achieve this we will work with the models
$\Xl{\alpha}{\beta}$ and $\Xe{\beta}{\alpha}$ for $\Xc{\alpha}{\beta}$
and $\Xc{\beta}{\alpha}$ respectively. The idea is to use the explicit
presentations for these models of the gerbes and construct the functor
$\FM$ in terms of data on the atlases. 

Since $\Xl{\alpha}{\beta} = [\ul{\alpha}{\beta}/\rl{\alpha}{\beta}]$
and $\Xe{\beta}{\alpha} = [\ue{\beta}{\alpha}/\re{\beta}{\alpha}]$ we
have natural structure morphisms
\[
\xymatrix@R-1pc{
\ggl : & \ul{\alpha}{\beta} \ar[r] & \Xl{\alpha}{\beta} \\
\gge : & \ue{\beta}{\alpha} \ar[r] & \Xe{\beta}{\alpha}
}
\]
for the lifting and extension
presentations. The (derived) pullback by $\ggl$ gives a natural functor
\[
\ggl^{*} : D^{b}_{1}(\Xl{\alpha}{\beta}) \to D^{b}(\ul{\alpha}{\beta}),
\] 
which sends complexes of sheaves on $\Xl{\alpha}{\beta}$ to objects in 
$D^{b}(\ul{\alpha}{\beta})$ preserved by the relations.

Explicitly, for a ${\mathcal L} \in D^{b}_{1}(\Xl{\alpha}{\beta})$,
the pullback
$\ggl^{*}{\mathcal L}$ is given by a pair $(L,\boldsymbol{f})$ where:
\begin{itemize}
\item $L$ is a bounded complex of sheaves on the atlas
$\ul{\alpha}{\beta} =X_{\varphi}\times_{B} X_{\beta}$.
\item ${\boldsymbol f} : p_{13}^{*}L \stackrel{q.i.}{\to}
p_{23}^{*}L\otimes \mycal{P}_{\ot}$ is a
quasi-isomorphism of complexes on $X_{\varphi}\times_B
X_{\varphi}\times_B X_{\beta}$,
satisfying the cocycle condition \eqref{eq-cocycle}.
\end{itemize}
Here $p_{ij}$ is the projection of $X_{\varphi}\times_{B}
X_{\varphi}\times_{B} X_{\beta}$ onto the
product of the $i$-th and $j$-th components.

Under the gerby Fourier-Mukai transform, ${\mathcal L}$ should go to
an object ${\mathcal Q} \in D^{b}_{-1}(\Xe{\beta}{\alpha})$. To
produce this object we will perform an integral transform from the derived
category of the atlas $\ul{\alpha}{\beta}$ to the derived category of
the atlas $\ue{\beta}{\alpha}$. Again, we would like to use the fact that
the pullback by $\gge$ gives a functor
\[
\gge^{*} : D^{b}_{-1}(\Xe{\beta}{\alpha}) \to D^{b}(\ue{\beta}{\alpha}),
\]
which sends complexes on $\Xe{\beta}{\alpha}$ to objects in
$D^{b}(\ue{\beta}{\alpha})$ preserved by the relations. In principle
this is all we can say about the images of $\gge^{*}$ since even for
schemes the derived categories of coherent sheaves do not necessarily
glue, see e.g. \cite{hartshorne-rd}. However in the case of
$\Xe{\beta}{\alpha}$ we can be much more precise. It is known
\cite[Theorem~A]{polishchuk-biextensions} that given a scheme $S$ and
a finite flat morphism $p : U \to S$, the derived category of coherent
sheaves on $S$ is equivalent to the category of pairs
$(F,\boldsymbol{\phi})$, where $F \in D^{b}(U)$ and $\boldsymbol{\phi}
: p_{1}^{*}F \widetilde{\to} p_{2}^{*}F$ is an isomorphism in
$D^{b}(U\times_{S} U)$ satisfying the cocycle condition in
$D^{b}(U\times_{S} U\times_{S} U)$. If in addition we are given a
${\mathbb C}^{\times}$ bundle $R \to U\times_{S} U$ equipped with a
biextension isomorphism, so that $[U/R]$ is a ${\mathcal
O}^{\times}$-gerbe on $S$, we can repeat the reasoning of
\cite[Theorem~A]{polishchuk-biextensions} verbatim to conclude that
the category $D^{b}_{1}([U/R])$ is equivalent to the category of pairs
$(G,\boldsymbol{\psi})$, where $G \in D^{b}(U)$ and $\boldsymbol{\psi}
: p_{1}^{*}G \widetilde{\to} p_{2}^{*}G\otimes R$ is an isomorphism in
$D^{b}(U\times_{S} U)$ satisfying the cocycle condition
\eqref{eq-cocycle} in $D^{b}(U\times_{S} U\times_{S} U)$. In
particular, since by construction the morphism $\ue{\beta}{\alpha} =
X_{\varphi} \to X_{\alpha}$ is finite and flat, we conclude that
$\gge^{*}$ identifies $D^{b}_{-1}(\Xe{\beta}{\alpha})$ with the
category of pairs $(Q,\boldsymbol{g})$ where:
\begin{itemize}
\item $Q$ is a bounded complex of sheaves on the atlas
$\ue{\beta}{\alpha} = X_{\varphi}$.
\item $\boldsymbol{g} : a_{\varphi}^{*}Q \stackrel{q.i.}{\to}
p_{2}^{*}Q\otimes p_{1}^{*}(M_{\beta}^{-1}\otimes M_{0})$ is a
quasi-isomorphism of complexes on $X[m]\times_B
X_{\varphi}$, satisfying the cocycle condition \eqref{eq-cocycle}.
\end{itemize}
Here $p_{i}$ is the projection of $X[m]\times_{B}
X_{\varphi}$ onto the $i$-th component.

\begin{rem}
The reader may wish to focus on the case when ${\mathcal L}$ is a line
bundle on $\Xl{\alpha}{\beta}$ of fiber degree zero, i.e. $L \to
X_{\varphi}\times_{B} X_{\beta}$ is a line bundle with
$\deg(L_{|\{x\}\times_{B} X_{\beta}}) = 0$ and
the existence of $\boldsymbol{f}$ is equivalent to having isomorphisms
$L_{|X_{\varphi}\times_{B} \{ y\}} \cong {\mathcal O}(my
)\otimes\Sigma_{\beta}^{-1}$ for all $y \in X_{\beta}$.

In this case the object ${\mathcal Q}$ should be a spectral datum on
the gerbe $\Xe{\beta}{\alpha}$ whose support is of degree one over
$B$, i.e. $Q$ is a torsion sheaf on $X_{\varphi}$
supported on $(q_{\varphi}^{n})^{-1}(s)$ for some section  $s \subset
X_{\alpha}$.
\end{rem}

\

\medskip

\noindent
We will construct the functor $\FM$ by first constructing a functor
between the derived categories on the atlases and then checking
that this functor preserves the relations.

On the level of atlases, consider the functor
\[
p_{1*} : D^{b}(X_{\varphi}\times_{B} X_{\beta}) \to D^{b}(X_{\varphi}).
\]
We now have the following:

\begin{prop} \label{prop-gerbyFM} The functor $p_{1*} :
D^{b}(X_{\varphi}\times_{B} X_{\beta}) \to D^{b}(X_{\varphi})$
preserves the relations defining $\Xl{\alpha}{\beta}$ and
$\Xe{\beta}{\alpha}$ and so descends to a well defined exact functor
\[
\FM := (\gge^{*})^{-1}\circ p_{1*} \circ \ggl^{*} :
D^{b}_{1}(\Xc{\alpha}{\beta}) \to D^{b}_{-1}(\Xc{\beta}{\alpha}). 
\]
\end{prop}
{\bf Proof.} Let $\ggl^{*}{\mathcal L} = (L,\boldsymbol{f})$ be as above. Let
$Q:= p_{1*}L$. We need to construct a quasi-isomorphism
of complexes on $X[m]\times_B X_{\varphi}$
\begin{equation} \label{eq-g}
\boldsymbol{g} : a_{\varphi}^{*}Q \to p_{2}^{*}Q\otimes
p_{1}^{*}(M_{\beta}^{-1}\otimes M_{0})
\end{equation}
which depends functorially on $\boldsymbol{f}$.

We start by noting that there is a natural commutative diagram
\[
\xymatrix@M+5pt{ X_{\varphi}\times_{B} X_{\varphi} \times_{B} X_{\beta}
\ar@<.5ex>[r]^-{p_{13}} \ar@<-.5ex>[r]_-{p_{23}} & X_{\varphi}
\times_{B} X_{\beta} \ar@{=}[d] \\
X_{\varphi}\times_{X_{\alpha}}
X_{\varphi} \times_{B} X_{\beta} \ar@{^{(}->}[u]
\ar@<.5ex>[r]^-{\op{pr}_{13}} \ar@<-.5ex>[r]_-{\op{pr}_{23}}
\ar[d]^-{\op{pr}_{12}} &
X_{\varphi} \times_{B} X_{\beta} \ar[d]^-{p_{1}} \\
X_{\varphi}\times_{X_{\alpha}}
X_{\varphi} \ar@<.5ex>[r]^-{\op{pr}_{13}}
\ar@<-.5ex>[r]_-{\op{pr}_{23}} &
X_{\varphi}
}
\]
Since the two bottom squares are fiber products we have the base
change formulas:
\begin{equation} \label{eq-base-change}
\begin{split}
\op{pr}_{1}^{*}p_{1*}L & = p_{12*}\op{pr}_{13}^{*}L \\
\op{pr}_{2}^{*}p_{1*}L & = p_{12*}\op{pr}_{23}^{*}L.
\end{split}
\end{equation}
Also, using the isomorphism $a_{\varphi}\times \op{id} : X[m]\times_{B}
X_{\varphi} \to X_{\varphi}\times_{X_{\alpha}} X_{\varphi}$ we reduce
the problem of finding the map \eqref{eq-g} to the equivalent problem
of constructing a map
\begin{equation} \label{eq-gg}
\op{pr}_{1}^{*}p_{1*}L \to \op{pr}_{2}^{*}p_{1*}L\otimes \Phi_{\beta}.
\end{equation}
of complexes on
$X_{\varphi}\times_{X_{\alpha}} X_{\varphi}$.
Now using \eqref{eq-base-change} and adjunction, this becomes
\begin{equation} \label{eq-ytf}
\op{pr}_{13}^*L \to \op{pr}_{23}^*L\otimes
\op{pr}_{12}^{*}\Phi_{\beta},
\end{equation}
on
$X_{\varphi}\times_{X_{\alpha}} X_{\varphi}\times_B X_{\beta}$.
Since
$\op{pr}_{13}^*L$ and $\op{pr}_{23}^*L$ are the restrictions of
$p_{13}^*L$ and $p_{23}^*L$, respectively, from
$X_{\varphi}\times_{B} X_{\varphi}\times_B X_{\beta}$
to
$X_{\varphi}\times_{X_{\alpha}} X_{\varphi}\times_B X_{\beta}$, the
restriction of our map ${\boldsymbol f}$ gives
\[
{\boldsymbol f}: \op{pr}_{13}^*L \to \op{pr}_{23}^*L\otimes
(\mycal{P}_{\ot|X_{\varphi}\times_{X_{\alpha}} X_{\varphi}\times_B
X_{\beta}}).
\]
Therefore, in order to reconstruct from ${\boldsymbol f}$
a map \eqref{eq-ytf} (equivalently ${\boldsymbol g}$), it suffices to
exhibit a canonical isomorphism
\begin{equation} \label{eq-iso}
\mycal{P}_{\ot|X_{\varphi}\times_{X_{\alpha}} X_{\varphi}\times_B
X_{\beta}} \cong \op{pr}_{12}^{*}\Phi_{\beta},
\end{equation}
on $X_{\varphi}\times_{X_{\alpha}} X_{\varphi}\times_B X_{\beta}$.  As
a first step in establishing \eqref{eq-iso} we note that both sides
are pullbacks of sheaves on $X[m]$. On the right hand side,
$\Phi_{\beta}$ was defined as $d^{*}(M_{\beta}^{-1}\otimes M_{0})$ for
the difference map $d : X_{\varphi}\times_{X_{\alpha}} X_{\varphi} \to
X[m]$. On the left hand side, it suffices (in view of the see-saw
principle) to argue that, for a point $\xi \in X[m]$, the restriction
\[
\mycal{P}_{\ot|\{\xi\} \times_{B} X_{\varphi}\times_B
X_{\beta}}
\]
is trivial. But by the definition of $\mycal{P}_{\ot}$ (see
Section~\ref{ssubsec-lifting}), this restriction can be identified
with $\xi^{\otimes m}$. Since $\xi$ has order $m$, we are done.

To conclude the construction of the map \eqref{eq-iso} and the proof
of the proposition, we need to show that the direct image
$R^{0}p_{1-2*}(\mycal{P}_{\ot|X_{\varphi}\times_{X_{\alpha}}
X_{\varphi}\times_B X_{\beta}})$ is isomorphic to the line bundle
$M_{\beta}^{-1}\otimes M_{0}$ on $X[m]$. For this it is useful to
identify $X_{\varphi}\times_{X_{\alpha}} X_{\varphi}\times_B
X_{\beta}$ with $X[m]\times_{B} X_{\varphi} \times_{B}
X_{\beta}$. Under that identification
$\mycal{P}_{\ot|X_{\varphi}\times_{X_{\alpha}} X_{\varphi}\times_B
X_{\beta}}$ becomes the pullback of $\mycal{P}_{|X[m]\times_{B} X}$
under the natural map
\[
\xymatrix@R-25pt{X[m]\times_{B} X_{\varphi} \times_{B} X_{\beta}
\ar[r] & X[m] \times_{B} X \\
(\xi,x,y) \ar@{|->}[r] & (\xi,q_{\beta}(y)).
}
\]
Thus, the existence of the isomorphism \eqref{eq-iso} is equivalent to
the existence of an isomorphism
\begin{equation} \label{eq-iso2}
\bp_{1,m\cdot 2}^{*}(\mycal{P}_{|X[m]\times_{B} X_{\beta}}) \cong
\bp_{1}^{*}(M_{\beta}^{-1}\otimes M_{0}),
\end{equation}
where $\bp_{1} : X[m]\times_{B} X_{\beta} \to X[m]$ is the projection
and $\bp_{1,m\cdot 2} : X[m]\times_{B} X_{\beta} \to X[m]\times_{B} X$ is
the map given by $(\xi,x) \mapsto (\xi,q_{\beta}(x))$.

Recall from Section~\ref{ssubsec-extension} that by definition we have
\[
\bp_{1}^{*}M_{\beta} = a_{\beta}^{*}\Sigma_{\beta}\otimes \bp_{2}^{*}
\Sigma_{\beta}^{-1}.
\]
Here $a_{\beta} : X[m]\times_{B} X_{\beta}$ is the action,
$\bp_{2} : X[m]\times_{B}  X_{\beta} \to X_{\beta}$ is the projection
on the second factor and $\Sigma_{\beta}$ is a line bundle on
$X_{\beta}$ of fiber degree $m$, which corresponds to the
`multiplication by $m$' map $q_{\beta} : X_{\beta} \to X$.

Look at the embedding $X[m]\times_{B} X_{\beta} \subset X\times_{B}
X_{\beta}$. The projections $\bp_{1}, \bp_{2}$ and the maps
$\bp_{1,m\cdot 2}$ and $a_{\beta}$ extend to the natural projections
$X\times_{B} X_{\beta} \to X$ and $X\times_{B} X_{\beta} \to
X_{\beta}$ and maps $X\times_{B} X_{\beta} \to X\times_{B} X$ and
$X\times_{B} X_{\beta} \to X_{\beta}$, which we will denote by the
same letters.

With this notation we have

\begin{lem}
\label{lem-difference} Let $\mycal{P} \to X\times_{B} X$ denote the
standard Poincare bundle. Then we have a natural isomorphism
\[
\bp_{1,m\cdot 2}^{*}\mycal{P}\otimes a_{\beta}^{*} \Sigma_{\beta}\otimes
\bp_{2}^{*} \Sigma_{\beta}^{-1} \cong \bp_{1}^{*}{\mathcal
O}_{X}(m\sigma).
\]
\end{lem}
{\bf Proof of the lemma.} We will use the see-saw principle. Let $\xi
\in X$ and let $b = \pi(\xi) \in B$. Then by viewing $\xi$ as a line
bundle of degree zero on $(X_{\beta})_{b}$ and using the fact that
$\Sigma_{\beta|X_{b}}$ is of degree $m$  we compute
\[
\begin{split}
\bp_{1,m\cdot 2}^{*}\mycal{P}_{|\{\xi \}\times (X_{\beta})_{b}} & =
q_{\beta}^{*}(\xi) = \xi^{\otimes m}, \\
(a_{\beta}^{*} \Sigma_{\beta}\otimes
\bp_{2}^{*} \Sigma_{\beta}^{-1})_{|\{\xi \}\times (X_{\beta})_{b}} & =
t_{\xi}^{*}L_{\beta}\otimes L_{\beta}^{-1} = \xi^{\otimes -m}.
\end{split}
\]
Thus for every $\xi$ we have
\[
(\bp_{1,m\cdot 2}^{*}\mycal{P}\otimes a_{\beta}^{*} \Sigma_{\beta}\otimes
\bp_{2}^{*} \Sigma_{\beta}^{-1})_{|\{\xi \}\times (X_{\beta})_{b}}
\cong {\mathcal O}_{(X_{\beta})_{b}} = {\mathcal O}_{(X_{\beta})_{b}}
\]
and so by the see-saw principle
$D_{\beta} := R^{0}\bp_{1*}(\bp_{1,m\cdot 2}^{*}\mycal{P}\otimes a_{\beta}^{*}
\Sigma_{\beta}\otimes \bp_{2}^{*} \Sigma_{\beta}^{-1})$ is a line
bundle on $X$ satisfying
\[
\bp_{1,m\cdot 2}^{*}\mycal{P}\otimes a_{\beta}^{*} \Sigma_{\beta}\otimes
\bp_{2}^{*} \Sigma_{\beta}^{-1} \cong
\bp_{1}^{*}D_{\beta}.
\]
To compute the bundle $D_{\beta}$ we consider a point $x \in X_{\beta}$. Let
$b = \pi_{\beta}(x)$. Restricting to $X_{b}\times \{x\}$ we get
\[
\begin{split}
\bp_{1,m\cdot 2}^{*}\mycal{P}_{|X_{b}\times \{x\}} & = q_{\beta}(x) \quad
\text{(considered as a line bundle of degree zero on $X_{b}$)} \\
a_{\beta}^{*}\Sigma_{\beta|X_{b}\times \{x\}} & =
t_{x}^{*}\Sigma_{\beta} \\
\bp_{2}^{*}\Sigma_{\beta|X_{b}\times \{x\}} & = {\mathcal O}_{X_{b}}.
\end{split}
\]
Next, by the defining relationship between $q_{\beta}$ and
$\Sigma_{\beta}$ we have that $q_{\beta}(x)$ is the line bundle of
degree zero on $X$ corresponding to ${\mathcal
O}_{(X_{\beta})_{b}}(mx)\otimes \Sigma_{\beta}^{-1}$ under the
identification $\op{Pic}^{0}((X_{\beta})_{b}) = \op{Pic}^{0}(X_{b})$.
Also, by the definition of a translation we have that
$t_{x}^{*}\Sigma_{\beta}$ is the tensor product of ${\mathcal
O}_{X_{b}}(m\sigma(b))$ with the line bundle of degree zero on $X_{b}$
corresponding to ${\mathcal O}_{(X_{\beta})_{b}}(-mx)\otimes
\Sigma_{\beta}$ under the identification
$\op{Pic}^{0}((X_{\beta})_{b}) = \op{Pic}^{0}(X_{b})$.

In other words, we have
\[
(\bp_{1,m\cdot 2}^{*}\mycal{P}\otimes a_{\beta}^{*} \Sigma_{\beta}\otimes
\bp_{2}^{*} \Sigma_{\beta}^{-1})_{X_{b}\times \{x\}} \cong {\mathcal
O}_{X_{b}}(m\sigma(b)),
\]
for all $x \in X_{\beta}$. This implies that up to a twist by a
pullback of a line bundle on $B$ we have $D_{\beta} \cong {\mathcal
O}_{X}(m\sigma)$. Finally, to fix the choice of this line bundle on
$B$ we look at the restriction of $\bp_{1,m\cdot 2}^{*}\mycal{P}\otimes
a_{\beta}^{*} \Sigma_{\beta}\otimes \bp_{2}^{*} \Sigma_{\beta}^{-1}$
on $\sigma\times_{B} X_{\beta}$ which is clearly isomorphic to
${\mathcal O}_{X_{\beta}}$. Hence the line bundle on $B$ is trivial
and the lemma is proven. \ \hfill $\Box$

\

\medskip

\noindent
In view of Lemma~\ref{lem-difference}, the only thing left to check in
order to establish the isomorphism \eqref{eq-iso2} is that the line
bundle $M_{0}$ on $X[m]$ is isomorphic to the restriction ${\mathcal
O}_{X}(m\sigma)_{|X[m]}$. However applying the same reasoning we used
in the proof of Lemma~\ref{lem-difference} to the projections $p_{1},
p_{2} : X\times_{B} X \to X$ and the obvious maps $p_{1,m\cdot 2} :
X\times_{B} X \to X\times_{B} X$ and $a_{0} : X\times X \to X$, we see
that $p_{1,m\cdot 2}^{*}\mycal{P}\otimes a_{0}^{*}{\mathcal
O}_{X}(m\sigma)\otimes p_{2}^{*}{\mathcal O}_{X}(-m\sigma) \cong
p_{1}^{*}{\mathcal O}_{X}(m\sigma)$. On the other hand, from the
definition of the Poincare bundle we have that
$p_{1,m\cdot 2}^{*}\mycal{P}_{|X[m]\times_{B} X} \cong {\mathcal O}$ and so
$M_{0} \cong {\mathcal O}_{X}(m\sigma)_{|X[m]}$. This finishes the
proof of the existence of \eqref{eq-iso}. To complete the proof of the
proposition, it only remains to note that since $(Q,\boldsymbol{g})$
was constructed from $(L,\boldsymbol{f})$ by means of the pushforward
via $X_{\varphi}\times_{B} X_{\beta} \to X_{\beta}$ and the fixed
isomorphism \eqref{eq-iso}, it follows that $\boldsymbol{g}$ will
satisfy the cocycle condition whenever $\boldsymbol{f}$ does. The
proposition is proven. \ \hfill $\Box$

\subsubsection{Categorical yoga for equivalences} \label{ssubsec-yoga}

We have constructed a functor $\FM : D^{b}_{1}(\Xc{\alpha}{\beta})
\to D^{b}_{-1}(\Xc{\beta}{\alpha})$. We are going to prove that it
is an equivalence. In order to do this, it is convenient to recall
some general criteria, due to Bondal-Orlov and Bridgeland, for
equivalences of triangulated categories.

Throughout this subsection we let $\bF : \mycal{A} \to \mycal{B}$ be
an exact functor between triangulated categories. A class $\Omega$  of
objects in $\mycal{A}$ is called a {\em spanning class} if for every
$a \in \op{ob} \mycal{A}$, the left orthogonality condition
\[
\op{Hom}^{i}_{\mycal{A}}(a,\omega) = 0,
\quad \text{for all} \quad i \in {\mathbb Z}, \omega \in \Omega
\]
implies that $a = 0$, and similarly on the right. Recall the following

\

\medskip

\noindent
{\bf Theorem \cite{bondal-orlov-flops}} {\em Assume that $\Omega$ is a
spanning class for $\mycal{A}$ and that the functor $\bF : \mycal{A}
\to \mycal{B}$ has left and right adjoints. Then $\bF$ is fully
faithful if and only if it is orthogonal:
\[
\bF : \op{Hom}_{\mycal{A}}^{i}(\omega_{1},\omega_{2}) \widetilde{\to}
\op{Hom}_{\mycal{B}}^{i}(\bF\omega_{1},\bF\omega_{2}), \quad \text{for
all } i \in {\mathbb Z}, \omega_{1}, \omega_{2} \in \Omega.
\]
}
\

\medskip

\noindent
Assume now that our triangulated category $\mycal{A}$ is linear. A
functor $\bS_{\!\mycal{A}} : \mycal{A} \to \mycal{A}$ is called
\cite{bondal-kapranov} a {\em
Serre functor} for $\mycal{A}$ if it is an exact equivalence and
induces bifunctorial isomorphisms
\[
\op{Hom}_{\mycal{A}}(a,b) \to
\op{Hom}_{\mycal{A}}(b,\bS_{\!\mycal{A}}a)^{\vee}, \quad
\text{for all } a,b \in \op{ob}\mycal{A},
\]
satisfying compatibility with compositions.
The basic example of a Serre functor is $\bS : D^{b}(X) \to D^{b}(X)$,
$\bS a := a\otimes \omega_{X}[n]$, where $X$ is a smooth
$n$-dimensional projective variety and $\omega_{X}$ is the canonical
bundle of $X$. If a Serre functor exists, it is unique up to an
isomorphism of functors. We are now ready to state the main
equivalence criterion we will be using:

\

\medskip

\noindent
{\bf Theorem \cite{bridgeland,bkr}} {\em Assume that $\mycal{A}$ and
$\mycal{B}$ have Serre functors $\bS_{\!\mycal{A}}$,
$\bS_{\!\mycal{B}}$, that $\mycal{A} \neq 0$, $\mycal{B}$ is
indecomposable, and that $\bF : \mycal{A} \to \mycal{B}$ has a left
adjoint. Then $\bF$ is an equivalence if it is fully faithful and it
intertwines the Serre functors: $\bF\circ \bS_{\!\mycal{A}}(\omega) =
\bS_{\!\mycal{B}}\circ \bF(\omega)$ on all elements $\omega \in
\Omega$ in the spanning class. }

\

\medskip

\noindent
We want to show that our gerby Fourier-Mukai functor
$\FM : D^{b}_{1}(\Xc{\alpha}{\beta})
\to D^{b}_{1}(\Xc{\beta}{\alpha})$
is an equivalence. The results above suggest that we should first
exhibit Serre functors for $D^{b}_{1}(\Xc{\alpha}{\beta})$ and
$D^{b}_{1}(\Xc{\beta}{\alpha})$,
and find a suitable spanning class for
$D^{b}_{1}(\Xc{\alpha}{\beta})$. These results, which do not involve
$\FM$, are carried out in section
\ref{ssubsec-serre} below. In section \ref{ssubsec-ort} we then
complete the argument by showing that  our $\FM$ preserves
orthogonality and intertwines the Serre functors.

\subsubsection{Serre functors and spanning classes for ${\mathcal
O}_{X}^{\times}$-gerbes}
\label{ssubsec-serre}

Let $X$ be an $n$-dimensional smooth projective variety. Let $c :
{\leftidx{_{\alpha}}{X}{}} \to X$ be an ${\mathcal O}_{X}^{\times}$-gerbe
corresponding to an element $\alpha \in H^{2}(X,{\mathcal
O}^{\times}_{X})$.

\begin{claim} \label{claim-serre} The functor

\[
\xymatrix@R-25pt{
\bS : & D^{b}_{1}(\leftidx{_{\alpha}}{X}{}) \ar[r] & D^{b}_{1}(\leftidx{_{\alpha}}{X}{})   \\
& a \ar@{|->}[r] & a {\otimes} c^{*}\omega_{X}[n]
}
\]
is a Serre functor.
\end{claim}
{\bf Proof.}
For any $a,b \in D^{b}_{1}(\leftidx{_{\alpha}}{X}{})$, we need a natural isomorphism

\[
\op{Hom}_{D^{b}_{1}(\leftidx{_{\alpha}}{X}{})}(a,b) \to
\op{Hom}_{D^{b}_{1}(\leftidx{_{\alpha}}{X}{})}(b,\bS a)^{\vee}.
\]
Since $a,b$ have weight 1, $\mycal{R}Hom_{\leftidx{_{\alpha}}{X}{}}(a,b)$
has weight $0$, so there exists a unique
$\mycal{H}(a,b) \in D^{b}(X)$ such that

\[
\mycal{R}Hom_{\leftidx{_{\alpha}}{X}{}}(a,b)=c^{*}\mycal{H}(a,b).
\]
It follows that

\[
\op{Hom}_{D^{b}_{1}(\leftidx{_{\alpha}}{X}{})}(a,b)=R{\Gamma}_{X}(\mycal{H}(a,b)).
\]
Similarly,

\[
\mycal{R}Hom_{\leftidx{_{\alpha}}{X}{}}(b,\bS a)=
c^{*}(\mycal{H}(b,a) \otimes \omega_{X}[n]),
\]
so

\[
\op{Hom}_{D^{b}_{1}(\leftidx{_{\alpha}}{X}{})}(b,\bS a)^{\vee}=
R{\Gamma}_{X}(\mycal{H}(b,a)\otimes \omega_{X}[n])^{\vee}=
R{\Gamma}_{X}(\mycal{H}(b,a)^{\vee}),
\]
where the last step uses the usual Serre duality.
So all we need is the identification

\[
\mycal{H}(a,b) \widetilde{\to} \mycal{H}(b,a)^{\vee},
\]

or a non-degenerate pairing on
$\mycal{H}(a,b) \times \mycal{H}(b,a).$
But since

\[
c^*: D^b(X) \widetilde{\to} D^{b}_{0}(\leftidx{_{\alpha}}{X}{})
\]
is an equivalence of categories, this follows immediately from the non-
degenerate pairing on
$\mycal{R}Hom_{\leftidx{_{\alpha}}{X}{}}(a,b) \times
\mycal{R}Hom_{\leftidx{_{\alpha}}{X}{}}(b,a)$
given by the trace.
\ \hfill $\Box$

\

\medskip

\noindent
Since our functor $\bF = \FM$ was constructed as a push-forward on the
atlases, it has an obvious left adjoint $\bG$ corresponding to the
pullback functor on the atlases. Therefore $\FM$ also has a right
adjoint, namely $S_{D^{b}_{1}(\Xc{\alpha}{\beta})}\circ \bG \circ
S_{D^{b}_{1}(\Xc{\beta}{\alpha})}^{-1}$. 

Fix a point $p \in \Xc{}{\beta}$. Since the restriction of
$\Xc{\alpha}{\beta}$ to $p$ is the trivial gerbe on $p$ for any
$\alpha$, the torsion sheaf ${\mathcal O}_{p}$ can be considered as a
weight one sheaf on $\Xc{\alpha}{\beta}$ for any $\alpha$. For our
spanning class $\Omega$ we take the structure sheaves of points on the
space $\Xc{}{\beta}$, viewed as sheaves of weight one on the stack
$\Xc{\alpha}{\beta}$.

\begin{claim} \label{claim-spanning} Let $c :
  {\leftidx{_{\alpha}}{X}{}} \to X$ be an  
${\mathcal O}_{X}^{\times}$-gerbe on a smooth projective $X$. Then the
class $\Omega$ consisting of structure sheaves ${\mathcal O}_{p}$
of points on $X$, viewed as sheaves of weight one on 
$\leftidx{_{\alpha}}{X}{}$, 
is a spanning class for $D^{b}_{1}(\leftidx{_{\alpha}}{X}{})$.
\end{claim}
{\bf Proof.} In order to show that the class $\Omega$  is a
spanning class, we need to show that for every
$a \in \op{ob} D^{b}_{1}(\leftidx{_{\alpha}}{X}{})$, the left orthogonality
condition
\[
\op{Hom}^{i}_{D^{b}_{1}(\leftidx{_{\alpha}}{X}{})}(a,{{\mathcal O}}_{p}) = 0,
\quad \text{for all} \quad i \in {\mathbb Z}, p \in X
\]
implies that $a = 0$. We also need the analogous result on the right,
but this follows using the Serre functor. We can also reduce to the
case that $a$ is represented by a single sheaf on
$\leftidx{_{\alpha}}{X}{}$, i.e. $a$ is an $\alpha$-twisted sheaf on
$X$. Now such an $a$ is specified in terms of its sections on an
appropriate etale atlas $U$ plus some $\alpha$-twisted gluing. In
order to conclude that $a=0$, it suffices to show that every $p \in X$
has a neighborhood $U'$ on which $a=0$.  But by restricting to a small
enough neighborhood $U'$ of $p$ in $U$, we can get the class $\alpha$
to vanish. The restriction of $a$ to $U'$ and the ${{\mathcal O}}_{p}$
for $p \in U'$ become ordinary sheaves.  The group
$\op{Hom}^{i}_{D^{b}_{1}(\leftidx{_{\alpha}}{X}{})}(a,{{\mathcal
O}}_{p})$ can be computed on either $U$ or $U'$.  Therefore, the
orthogonality condition forces $a$ to vanish on $U'$, which is what we
need.  \ \hfill $\Box$

\subsubsection{Orthogonality and intertwining} \label{ssubsec-ort}

Now that we have a Serre functor and a spanning class, we are ready
to apply the general results of subsection \ref{ssubsec-yoga}
to our gerby Fourier-Mukai functor $\FM$.

\begin{claim} \label{claim-ortho}
The gerby Fourier-Mukai functor $\FM : D^{b}_{1}(\Xc{\alpha}{\beta})
\to D^{b}_{-1}(\Xc{\beta}{\alpha})$ is orthogonal on $\Omega$:
\[
\FM : \op{Hom}_{D^{b}_{1}(\Xc{\alpha}{\beta})}^{i}
({{\mathcal O}}_{x_1},{{\mathcal O}}_{x_2}) \widetilde{\to}
\op{Hom}_{D^{b}_{1}(\Xc{\beta}{\alpha})}^{i}
(\bF{{\mathcal O}}_{x_1},\bF{{\mathcal O}}_{x_2}), \quad \text{for
all } i \in {\mathbb Z}, x_1, x_2 \in X_{\beta}.
\]

\end{claim}
{\bf Proof.}
Recall (\ref{prop-gerbyFM})  that our functor $\FM$ descends from
$p_{1*}: D^{b}(X_{\varphi}\times_{B} X_{\beta}) \to D^{b}(X_{\varphi})$.
Let $b_1, b_2 \in B$ be the images of $x_1, x_2 \in X_{\beta}$. If $b_1
\neq b_2$ then the supports are disjoint, so the $\op{Hom}^i$  on
both sides clearly vanish. Assume then that $b_1 = b_2 = b$.
In this case, the structure sheaf
${{\mathcal O}}_{{\Xl{\alpha}{\beta}} \vert {x_i}}$
is supported on the fiber
$C_{x_i}=X_{\varphi}\times_{B}(x_i)$,
and both supports map isomorphically to
$C_{b}=X_{\varphi}\times_{B}(b) \subset X_{\varphi}$.
Now $\FM({{\mathcal O}}_{{\Xl{\alpha}{\beta}} \vert {x_i}})$
is the line  bundle
${\mycal L}_{\beta}(-mx_i)$ on $C_b$, so
\[
\begin{split}
\op{Hom}_{D^{b}_{1}(\Xc{\beta}{\alpha})}^{i}
(\FM{{\mathcal O}}_{x_1},\FM{{\mathcal O}}_{x_2}) & =
\op{Hom}_{C_b}({\mycal L}_{\beta}(-mx_1), {\mycal L}_{\beta}(-mx_2)) \\
& = \op{Hom}_{D^{b}_{1}(\Xc{\alpha}{\beta})}^{i}
({{\mathcal O}}_{x_1},{{\mathcal O}}_{x_2}),
\end{split}
\]
completing the proof.
\ \hfill $\Box$

\

\medskip

\noindent
We note that both sides of the claim vanish unless $x_1, x_2$ differ by
a point of $X[m]$, in which case they define isomorphic sheaves. The
spanning class $\Omega$ may therefore be taken to be parametrized by
$X=X_{\beta}/X[m]$ rather than by $X_{\beta}$.

\begin{claim} \label{claim-intertwine}
The gerby Fourier-Mukai functor
$\FM : D^{b}_{1}(\Xc{\alpha}{\beta}) \to D^{b}_{-1}(\Xc{\beta}{\alpha})$
intertwines the Serre functors, i.e.:
$\FM\circ \bS_{\!\Xc{\alpha}{\beta}}({{\mathcal O}}_{p}) =
\bS_{\!\Xc{\beta}{\alpha}}\circ \FM({{\mathcal O}}_{p})$
for all points $p \in X_{\beta}$.

\end{claim}
{\bf Proof.} Follows immediately from the fact that $\FM {\mathcal
O}_{x}$ is supported on $C_{b}$ and that the canonical sheaf of
$X_{\beta}$ restricts to the trivial line bundle on $C_{b}$. \ \hfill $\Box$

\section{Surfaces} \label{sec-K3}

In case $X$ is a surface, we can refine the previous results to
include the singular fibers. On a surface, any pair of classes
$\alpha, \beta \in \TSh_{an}(X)$ are complementary in the sense
of
subsection \ref{subsec-T}, by Corollary \ref{cor-surface}, so the
gerbes $\Xc{\alpha}{\beta}, \Xc{\beta}{\alpha}$ are always well-
defined. When $m\beta =0$ we will construct the lifting presentation
of $\Xc{\alpha}{\beta}$ and the extension presentation of
$\Xc{\beta}{\alpha}$. Then we will exhibit a Fourier-Mukai transform
$\FM$ between these presentations. Finally, we will show that $\FM$ is
an equivalence of categories by verifying the criterion of
Bondal-Orlov and Bridgeland.

We assume throughout that $X$ is a smooth surface, $B$ is a smooth
curve, and the elliptic fibration $\pi: X \to B$ has at most singular
fibers of type $I_1$. Since every such elliptic surface is uniquely
determined by its monodromy representation it is clear that we can
always extend $X$ to a smooth compact relatively minimal elliptic
surface whose base curve is a suitable compactification of
$B$. Furthermore, by Kodaira's classification of compact complex surfaces
it follows that every smooth compact elliptic surface is K\"{a}hler
(in fact algebraic). Therefore $X$ must be K\"{a}hler as well.

\subsection{The lifting presentation} \label{subsec-L}
Our first goal is to construct the lifting presentation of
$\Xc{\alpha}{\beta}$, in a way that restricts to the previously
constructed presentation on the non-singular fibers. We start with
the second projection $p_2: X_{\varphi} \times_B X_{\beta} \to X_{\beta}$.
Unfortunately, this is {\em not} an atlas for the gerbe
$\Xc{\alpha}{\beta}$. The problem can be traced back to the fact that
the threefold $X_{\varphi} \times_B X_{\beta}$ is singular. So let
\[
Y:= \widehat{X_{\varphi} \underset{B}{\times} X_{\beta}}
\]
be a small resolution of $X_{\varphi} \times_B X_{\beta}$. Now $Y$ is
smooth and equipped with flat morphisms $\nu_{1} : Y \to X_{\varphi}$
and $\nu_{2} : Y \to X_{\beta}$ which lift $p_{1}$ and $p_{2}$
respectively. There is an induced map
\[
\nu_2^*: Br'_{an}(X_{\beta}) \to Br'_{an}(Y).
\]
We claim that $Y$ is an atlas for $\Xc{\alpha}{\beta}$, i.e. that
$\Xc{\alpha}{\beta}$ has a presentation:

\begin{equation} \label{LP}
\gp
{\rl{\alpha}{\beta}}
{\ul{\alpha}{\beta}}
{X_{\beta}}
{p_{1}}
{p_{2}}
{\nu_{2}}
\end{equation}
where
\[
\begin{split}
\ul{\alpha}{\beta} & := Y, \\
\nu_{2} & : Y \to X_{\beta} \quad  \text{is the second projection, and } \\
\rl{\alpha}{\beta} & := Y \times_{\Xc{\alpha}{\beta}} Y.
\end{split}
\]
We call (\ref{LP}) the {\em Lifting Presentation} of
$\Xc{\alpha}{\beta}$. By Remark~\ref{rem-presentations} (iv), the fact
that (\ref{LP}) is indeed a presentation follows from:

\begin{lem} $\nu_2^* (T_{\beta}(\alpha)) = 0.$ \label{lem-L-atlas}
\end{lem}
{\bf Proof.} Let $B^{o} \subset B$ be the complement of the discriminant,
and let $Y^{o} \subset Y, X^{o} \subset X$ be the inverse images of $B^{o}$.
The maps
\[
Y^{o} \stackrel{i}{\hookrightarrow} Y \xrightarrow{{\nu_2}} X_{\beta}
\]
lead via the exponential sequence to the diagram:

\[
\xymatrix@M+3pt{
& & H^2_{an}(X_{\beta}, {{\mathcal O}}^{\times})
\ar[r]^-{\partial}
\ar[d]^-{\nu_2^*}
& H^3(X_{\beta}, {\mathbb Z})=0
\ar[d]^-{\nu_2^*} \\
0 \ar[r] &  H^2_{an}(Y, {{\mathcal O}}_Y) / H^2(Y, {\mathbb Z})
\ar@{^{(}->}[r]^-{{\op{exp}}_Y}
\ar[d]^-{i_{{\mathcal O}}^*} &
H^2_{an}(Y, {{\mathcal O}}^\times)
\ar[r]^-{\partial}
\ar[d]^-{i_{{{\mathcal O}}^{\times}}^*} &
H^3(Y, {\mathbb Z}) \ar[d] \\
0 \ar[r] &  H^2_{an}(Y^{o}, {{\mathcal O}}_{Y^{o}}) / H^2(Y^{o}, {\mathbb Z})
\ar@{^{(}->}[r]^-{{\op{exp}}_{Y^{o}}} &
H^2_{an}(Y^{o}, {{\mathcal O}}^\times)
\ar[r]  &
H^3(Y^{o}, {\mathbb Z}).
}
\]

We have

\[
\partial \nu_2^* (T_{\beta}(\alpha)) = \nu_2^* \partial (T_{\beta}(\alpha)) =
\nu_2^* 0 =0,
\]
so
\[
\nu_2^* (T_{\beta}(\alpha)) = {\op{exp}}_Y (a)
\]
for some $a \in H^2_{an}(Y, {{\mathcal O}}_{Y}) / H^2(Y, {\mathbb Z})$. We
know from \ref{theo-lifting-vs-coho} that $Y^{o}$ is an atlas for the
restriction of $T_{\beta}(\alpha)$ to $X^{o}$, so
\[
0 =
i_{{{\mathcal O}}^{\times}}^* {\nu_2^*} (T_{\beta}(\alpha)) =
i_{{{\mathcal O}}^{\times}}^* {\op{exp}}_Y (a) =
{\op{exp}}_{Y^{o}} i_{{{\mathcal O}}}^* (a).
\]
But $H^2_{an}(Y, {{\mathcal O}}_{Y}) / H^2(Y, {\mathbb Z})$ is a birational
invariant, so $i_{{\mathcal O}}^*$ is an isomorphism. Since
${{\op{exp}}_{Y^{o}}}$ is injective, we see that $a=0$, so we are done.
\ \hfill $\Box$

\begin{contra}
Let $A^{o} \subset A$ be an open subset in a smooth variety. By the
birational invariance of cohomological Brauer groups
\cite{milne-book}, the restriction map:
\[
H^2_{\et}(A, {{\mathcal O}}^{\times}_{A}) \to
H^2_{\et}(A^{o}, {{\mathcal O}}^{\times}_{A^{o}})
\]
is injective. Nevertheless, the analogous map:
\[
H^2_{an}(A, {{\mathcal O}}^{\times}_{A}) \to
H^2_{an}(A^{o}, {{\mathcal O}}^{\times}_{A^{o}})
\]
may fail to be injective. In our situation, all we were able to prove,
and fortunately all that was needed, was that
$H^2_{an}(Y, {{\mathcal O}}^{\times}_{Y}) \to
H^2_{an}(Y^{o}, {{\mathcal O}}^{\times}_{Y^{o}})$
is injective on the image of
$H^2_{an}(X_{\beta}, {{\mathcal O}}^{\times}_{X_{\beta}}).$

As an example, let $C \subset {\mathbb P}^3$ be a smooth curve of genus
$\geq 2$, let $A$ be the blowup of ${\mathbb P}^3$ along $C$, and let
$A^{o} :={\mathbb P}^3 \setminus C$. Then by the exponential sequence,
\[
\begin{split}
H^2_{an} (A^{o},{{\mathcal O}}^{\times}_{A^{o}}) & =H^3(A^{o}; {\mathbb Z}), \\
H^2_{an}(A  , {{\mathcal O}}^{\times}_{A}  ) & =H^3(A  ; {\mathbb Z}),
\end{split}
\]
but by excision, $H^3(A^{o}; {\mathbb Z}) \cong
H^3({\mathbb P}^3, C; {\mathbb Z}) \cong {\mathbb Z},$
while $H^3(A; {\mathbb Z}) \cong H^1(C, {\mathbb Z})$.
\end{contra}

\subsection{The extension presentation} \label{subsec-E}

Next, we want to construct the extension presentation of
$\Xc{\beta}{\alpha}$, in a way that restricts to the previously
constructed extension presentation on the complement of the singular
fibers. Fix
$\varphi \in \TSh_{an}(X)$ satisfying $m \cdot \varphi = \alpha$.
Let $X_{\alpha}^{o}, X_{\varphi}^{o}$ be the inverse images in
$X_{\alpha}, X_{\varphi}$ of $B^{o}$, the complement of the discriminant
in $B$. In the non-singular case, our atlas was given by the
multiplication-by-$m$ map
$q_{\varphi}: X_{\varphi}^{o} \to X_{\alpha}^{o}$.
Unfortunately, this does not extend to a morphism
$q_{\varphi}: X_{\varphi} \to X_{\alpha}$.
Instead, we will construct another (singular!) surface
$\widehat{X}_{\varphi}$
with a birational morphism
$\widehat{X}_{\varphi} \to X_{\varphi}$
and a flat morphism
$\hat{q}_{\varphi} : \widehat{X}_{\varphi} \to X_{\alpha}$
which restricts to the previous $q_{\varphi}^{o}$.
This data gives a commutative diagram:
\[
\xymatrix@M+3pt{
H^{2}_{an}(X_{\alpha},{\mathcal O}^{\times})
\ar[r]^-{\hat{q}_{\varphi}^{*}}
\ar[d]^-{\cong} &
H^{2}_{an}(\widehat{X}_{\varphi},{\mathcal O}^{\times})
\ar[d]^-{\cong} \\
H^{2}_{an}(X_{\alpha}^{o},{\mathcal O}^{\times})
\ar[r]^{(q_{\varphi}^{o})^*} &
H^{2}_{an}({X_{\varphi}^{o}},{\mathcal O}^{\times}) }
\]
The exponential sequence shows that the two vertical maps are
isomorphisms, as in the proof of Lemma \ref{lem-L-atlas}: this uses
the fact that $H^2_{an}(\widehat{X}_{\varphi}, {{\mathcal
O}}_{\widehat{X}_{\varphi}})/H^2(\widehat{X}_{\varphi}, {\mathbb Z})$
is a birational invariant, and that $\ker[H^3(\widehat{X}_{\varphi},
{\mathbb Z}) \to H^3(\widehat{X}_{\varphi}, {\mathbb R})] = 0$,
$\ker[H^3(X_{\alpha}, {\mathbb Z})\to H^3(X_{\alpha}, {\mathbb
R})]=0$, which in turn follows from the observation that the third
cohomology of a smooth 4-manifold has no torsion and that
$\widehat{X}_{\varphi}$ and $X_{\alpha}$ are K\"{a}hler
surfaces. Since $(q_{\varphi}^{o})^*$ kills all classes of order $m$,
it follows that so does $\hat{q}_{\varphi}^{*}$, so
$\widehat{X}_{\varphi}$ is indeed an atlas.

In order to construct $\widehat{X}_{\varphi}$ we have to resolve the
rational map $q_{\varphi} : X_{\varphi} \dashrightarrow
X_{\alpha}$. For that we can work locally in the complex topology on
$B$ near a point $p \in B \setminus B^{o}$, i.e.  we can replace $B$
by a small disc centered at $p$. Over this disc the group scheme
$X^{\sharp}[m]$ has a subgroup scheme $I \subset X^{\sharp}[m]$ of
cycles invariant under the local monodromy around $p$.  Since by
assumption the singular fibers of $X$ are of type $I_{1}$, the group
scheme $I$ is isomorphic to $B\times ({\mathbb Z}/m)$ and consists of
all the sections in $X[m]$ over the disc that pass through smooth
points of the fiber $X_{p}$. Translations by such sections give rise
to a well defined action of $I$ on $X_{\varphi}$ which fixes the
singular point $x_{p}$ of the fiber $(X_{\varphi})_{p}$. Therefore
over our disk the rational
map  $q_{\varphi}$  decomposes as
\[
\xymatrix{
X_{\varphi} \ar@{-->}[rr]^-{q_{\varphi}} \ar[rd]_-{s} & & X_{\alpha}
 \\ & X_{\varphi}/I \ar@{-->}[ru] &
}
\]
where $s$ denotes the quotient map. The surface $X_{\varphi}/I$ has a
unique singularity at the image of the point $x_{p}$ and the map
$X_{\varphi}/I \to B$ is a flat genus one fibration. A straightforward
local computation at the singular point of the $I_{1}$ fiber of
$X_{\varphi}$ shows that in suitably chosen local coordinates $(z,w)$
near $x_{p}$ the generator of $I$ acts as $(z,w) \mapsto (\zeta
z,\zeta w)$, where $\zeta$ is a primitive $m$-th root of unity. This
implies that the singularity of $X_{\varphi}/I$ is of type
$A_{m-1}$. The minimal resolution $\widehat{X_{\varphi}/I} \to
X_{\varphi}/I$ of $X_{\varphi}/I$ is a flat genus one fibration over
$B$ with a single $I_{m}$ fiber over $p$.

On the other hand, over our small disk, $X_{\alpha}$ retracts to the
singular fiber $(X_{\alpha})_{p}$ whose fundamental group is ${\mathbb
Z}$. Therefore there is a unique $m$-sheeted etale cover
$\widetilde{X}_{\alpha} \to X_{\alpha}$. By construction the covering
map commutes with the projections $\tilde{\pi}_{\alpha} :
\widetilde{X}_{\alpha} \to B$ and $\pi_{\alpha} : X_{\alpha} \to B$
and so the fiber $(\widetilde{X}_{\alpha})_{p}$ of
$\tilde{\pi}_{\alpha}$ over the point $p \in B$ is a Kodaira fiber of
type $I_{m}$. The surfaces $\widetilde{X}_{\alpha}$ and
$\widehat{X_{\varphi}/I}$ are clearly isomorphic outside of the
preimages of $p \in B$ and so are birational. Since as genus one
fibrations $\widetilde{X}_{\alpha}$ and $\widehat{X_{\varphi}/I}$ are
both relatively minimal, this implies that $\widetilde{X}_{\alpha}$
and $\widehat{X_{\varphi}/I}$ are actually isomorphic due to the
uniqueness \cite{kodaira} of the relatively minimal models.  Recall
next that by the work of Ito and Nakamura \cite{ito-nakamura,bkr} the
minimal resolution $\widehat{X_{\varphi}/I}$ of $X_{\varphi}/I$ can be
identified with the Hilbert scheme of $I$-clusters in
$X_{\varphi}$. In the spirit of \cite{bkr} consider the universal
closed subscheme $Z \subset X_{\varphi}\times \widehat{X_{\varphi}/I}$
with its natural projection to $X_{\varphi}$ and
$\widehat{X_{\varphi}/I} \cong \widetilde{X}_{\alpha}$. There is a
commutative diagram of spaces
\[
\xymatrix{& Z \ar[dl] \ar[dr] \\ \widetilde{X}_{\alpha}  \ar[dr] & &
X_{\varphi} \ar[dl] \\ & X_{\varphi}/I}
\]
where the morphism $Z \to X_{\varphi}$ is birational and surjective
and the morphism $Z \to \widetilde{X}_{\alpha}$ is finite and flat. In
particular the composite map
\begin{equation} \label{eq-localqhat}
Z  \to \widetilde{X}_{\alpha} \to X_{\alpha}
\end{equation}
is a finite and flat morphism which extends the multiplication-by-$m$ map
$q_{\varphi}^{o} : X_{\varphi}^{o} \to X_{\alpha}^{o}$. It is also
helpful to observe that the two intermediate maps appearing in the
construction of \eqref{eq-localqhat} are both Galois covers with
Galois group isomorphic to ${\mathbb Z}/m$. The first one is just the
quotient of $X_{\varphi}\times_{X_{\varphi}/I}
\widetilde{X}_{\alpha}$ by $I$ and the second is the etale Galois
cover $\widetilde{X}_{\alpha} \to X_{\alpha}$.

Our extension atlas $\widehat{X}_{\varphi}$ is obtained by gluing the
local surfaces $Z$ defined over small discs centered at
discriminant points $p \in B\setminus B^{o}$ to $X_{\varphi}^{o}$. We
will write $\varepsilon : \widehat{X}_{\varphi} \to X_{\varphi}$ for
the contraction map and $\hat{q}_{\varphi} : \widehat{X}_{\varphi} \to
X_{\alpha}$ for the finite flat map gluing each \eqref{eq-localqhat}
to $q_{\varphi}^{o}$. Note that by construction the surface
$\widehat{X}_{\varphi}$ is singular. It has isolated toroidal
singularities sitting over the singular points of the singular fibers
of $\widetilde{X}_{\alpha}$. In particular $\widehat{X}_{\varphi}$ is
a normal analytic surface.

\subsection{Duality for gerby genus one fibered surfaces} \label{subsec-FM}

With the description of the global lifting and extension presentations for
gerby genus one fibered surfaces in place, we are now ready to construct the
Fourier-Mukai functor between the derived categories of pure weight
one and to show that it is an equivalence. The key property of our
gerby surfaces, which makes the construction possible is the fact
that the gerbes appearing in the picture become trivial when we restrict
our attention to a piece of the surface sitting over a sufficiently
small open disk in $B$.

We recall our convention that all our direct and inverse images, as
well as the tensor product, are taken in the derived category. Thus for
any space $Z$, we will simply write $\otimes$ for the derived tensor
product $\otimes^{L}$ on $D^{b}(Z)$ and for any map of spaces $p: Z \to
T$ we will write $p_{*}, p^{*}, p_{!}, p^{!}$, for the corresponding
derived functors (whenever these functors make sense on the bounded
derived categories).

Following the pattern of the proof in section~\ref{ssubsec-FM} we
will construct a Fourier-Mukai functor
\[
\FM : D^{b}_{1}(\Xc{\alpha}{\beta}) \to D^{b}_{-1}(\Xc{\beta}{\alpha})
\]
by exhibiting an integral transform between the derived categories on
the atlases of the presentations $\Xl{\alpha}{\beta}$ and
$\Xe{\beta}{\alpha}$ and then checking that this functor preserves the
relations.

To avoid cumbersome notation we will write
\[
Y :=
\widehat{X_{\varphi}\underset{B}{\times} X_{\beta}}
\quad \text{and} \quad
S := \widehat{X}_{\varphi}
\]
for the atlases of the
presentations $\Xl{\alpha}{\beta}$ and $\Xe{\beta}{\alpha}$ respectively.

With this notation
the relations of $\Xl{\alpha}{\beta}$ are given by the total space
\[
\xymatrix@1{\op{tot}(\mycal{P}_{\ot}^{\times}) \ar[r] &
Y\times_{X_{\beta}} Y \ar@<.5ex>[r] \ar@<-.5ex>[r] & Y}.
\]
Here
$\mycal{P}_{\ot}$ denotes an appropriate line bundle on
$Y\times_{X_{\beta}} Y$ which extends the line bundle
$p_{1-2,m3}^{*}\mycal{P} \to X_{\varphi}^{o}\times_{B^{o}}
X_{\beta}^{o}$ discussed in section~\ref{ssubsec-lifting}. Note that
we know that such a line bundle exists due to
Lemma~\ref{lem-L-atlas}. Explicitly, the total space
$\op{tot}(\mycal{P}_{1-2,m3}^{\times})$ is isomorphic to the stacky
fiber product $Y\times_{\Xc{\alpha}{\beta}} Y$ (this product is a space again
by Lemma~\ref{lem-L-atlas}). We will write $Y^{o} =
X_{\varphi}^{o}\times_{B^{o}} X_{\beta}^{o}$ for the part of $Y$ sitting
over $B^{o}$, and we put $\mycal{P}_{\ot}^{o} := \mycal{P}_{\ot|Y^{o}}
= p_{\ot}^{*}\mycal{P}$.

Similarly the relations of the presentation $\Xe{\beta}{\alpha}$ are
given by the total space
\[
\xymatrix@1{\op{tot}(\Phi^{\times}_{\beta})
\ar[r] & S\times_{X_{\alpha}} S \ar@<.5ex>[r] \ar@<-.5ex>[r] &
S}.
\]
Here $\Phi_{\beta}$ denotes an appropriate
line bundle which extends the line bundle $d^{*}(M_{\beta}\otimes
M_{0}^{-1}) \to X_{\varphi}^{o}\times_{X_{\alpha}^{o}}
X_{\varphi}^{o}$ discussed in section~\ref{ssubsec-extension}.
As before, the existence of the line bundle $\Phi_{\beta}$ is
guaranteed by the observation that the class $T_{\alpha}(\beta)$ of
the extension gerbe vanishes (see section~\ref{subsec-E}) when we pull
it back to $S$. We will write $S^{o} = X_{\varphi}^{o}$ for the part
of $S$ sitting over $B^{o}$ and we put $\Phi_{\beta}^{o} :=
\Phi_{\beta|S^{o}} = d^{*}(M_{\beta}\otimes
M_{0}^{-1})$.

Using the above setup we can now identify the category
$D^{b}_{1}(\Xc{\alpha}{\beta})$ with the category of objects $L \in
D^{b}(Y)$ equipped with descent datum $\boldsymbol{f} : p_{1}^{*}L
\xrightarrow{q.i.} p_{2}^{*}L\otimes \mycal{P}_{\ot}$ on
$Y\times_{X_{\beta}} Y$, satisfying the cocycle condition described at
the end of
section~\ref{ssubsec-presentations}. Similarly we identify
$D^{b}_{-1}(\Xc{\beta}{\alpha})$ with the category of objects $Q \in
D^{b}(S)$ equipped with descent datum $\boldsymbol{g} : p_{1}^{*}Q
\xrightarrow{q.i.} p_{2}^{*}Q\otimes \Phi_{\beta}^{-1}$ on
$S\times_{X_{\alpha}} S$ satisfying the same cocycle condition as in
section~\ref{ssubsec-presentations}. These identifications reduce the problem of
constructing the Fourier-Mukai functor $\FM :
D^{b}_{1}(\Xc{\alpha}{\beta}) \to D^{b}_{-1}(\Xc{\beta}{\alpha})$ to
the problem of constructing a functor $\bF : D^{b}(Y) \to D^{b}(S)$
which maps descent data to descent data.

We will define the functor $\bF$ as an integral transform with respect
to a suitable kernel object $\Pi \in D^{b}(Y\times S)$. We proceed to
construct $\Pi$ by gluing together certain locally defined coherent
sheaves on $Y\times S$. We carry out this gluing in the analytic
topology to obtain the general functor $\bF$ we need. Note that in the
algebraic case, the kernel $\Pi$ still produces the correct functor
$\Pi$ in view of the GAGA principle.

First we look at the smooth part of the genus one fibrations
$X_{\alpha}$, $X_{\beta}$, $X_{\varphi}$, etc. As usual, we write
$B^{o} \subset B$ for the complement of the discriminant of the map
$\pi : X \to B$. Similarly, for any space (or stack)
$Z \to B$ mapping to $B$, we put $Z^{o} := Z\times_{B} B^{o}$. The
atlases $Y^{o}$ and $S^{o}$ for the gerbes $\Xl{\alpha}{\beta}^{o}$ and
$\Xe{\beta}{\alpha}^{o}$ can be described simply as
\[
\begin{split}
Y^{o} & = X_{\varphi}^{o}\times_{B^{o}} X_{\beta}^{o} \\
S^{o} & = X_{\varphi}^{o}
\end{split}
\]
and so over $B^{o}$ we recover the setup analyzed in
section~\ref{sec-duality}. In this setup the integral transform we
need was defined as the pushforward $p_{1*} :
D^{b}(X_{\varphi}^{o}\times_{B^{o}} X_{\beta}^{o}) \to
D^{b}(X_{\varphi}^{o})$ with respect to the projection on the first
factor. Equivalently, we can view this functor as the integral
transform whose kernel is the sheaf $\Delta_{\varphi*}{\mathcal
O}_{X_{\varphi}^{o}\times_{B^{o}} X_{\beta}^{o}}$ on
$X_{\varphi}^{o}\times_{B^{o}} X_{\varphi}^{o}\times_{B^{o}}
X_{\beta}^{o}$, where $\Delta_{\varphi}(x,y) := (x,x,y)$ for all
$(x,y) \in X_{\varphi}^{o}\times_{B^{o}} X_{\beta}^{o}$. In view of
this we set:
\[
\Pi^{o} := {\mathcal O}_{X_{\varphi}^{o}\times_{B^{o}} X_{\beta}^{o}}
\in \op{Coh}(X_{\varphi}^{o}\times_{B^{o}} X_{\beta}^{o}) =
\op{Coh}(Y^{o}\times_{X_{\varphi}^{o}} S^{o}).
\]
\

Let now $p \in B\setminus B^{o}$. Choose small analytic discs $p \in
U^{p} \subset B$ around each such $p$ so that $U^{p}\cap U^{q} =
\varnothing$ for $p \neq q$ and the genus one fibrations
$X_{\varphi}$ and $X_{\beta}$ both admit analytic sections over each
$U^{p}$. For any space (or stack) $Z \to B$ mapping to $B$ we will
write $Z^{p}$ for the restriction $Z\times_{B} U^{p}$. Note that
\[
Y\times_{X_{\varphi}} S = (Y^{o}\times_{X_{\varphi}^{o}} S^{o})\bigcup
\left(\coprod_{p \in B\setminus B^{o}} \left( Y^{p}\times_{X_{\varphi}^{p}} S^{p}
\right)\right),
\]
and so, in order to extend $\Pi^{o}$ to a globally defined sheaf on
$Y\times_{X_{\varphi}} S$, it suffices to construct coherent analytic sheaves
$\Pi^{p}$ on each $Y^{p}\times_{X_{\varphi}^{p}} S^{p}$ and isomorphisms
\[
\Pi^{o}_{|(Y^{p}\times_{X_{\varphi}^{p}} S^{p})\cap
(Y^{o}\times_{X_{\varphi}^{o}} S^{o})} \cong
\Pi^{p}_{|(Y^{p}\times_{X_{\varphi}^{p}}
S^{p})\cap (Y^{o}\times_{X_{\varphi}^{o}} S^{o})}.
\]
\

Fix a $p \in B\setminus B^{o}$ and let $s_{\varphi} : U^{p} \to
X_{\varphi}^{p}$ and $s_{\beta} : U^{p} \to X_{\beta}^{p}$ be local
holomorphic sections. Since the rational map $q_{\varphi} :
X_{\varphi} \dashrightarrow X_{\alpha}$ is a well defined morphism
when restricted to the complement $X_{\varphi}^{\sharp}$ of the
singular points of the singular fibers of $\pi_{\varphi} : X_{\varphi}
\to B$, and since $s_{\varphi}(U^{p}) \subset X_{\varphi}^{\sharp}$,
it follows that $s_{\varphi}$ induces also a well defined section
$s_{\alpha} := q_{\varphi}\circ s_{\varphi} : U^{p} \to
X_{\alpha}^{p}$ of $X_{\alpha}^{p}$. Consider now the natural map
\[
\nu_{1}\times \hat{q}_{\varphi} : Y^{p}\times_{X_{\varphi}^{p}} S^{p}
\to X_{\beta}^{p}\times_{U^{p}} X_{\alpha}^{p}
\]
and the isomorphism
\[
t_{-s_{\beta}}\times t_{-s_{\alpha}} : X_{\beta}^{p}\times_{U^{p}}
X_{\alpha}^{p} \to X^{p}\times_{U^{p}} X^{p}.
\]
Let $\mycal{P}^{p}$
denote the pullback
\[
\mycal{P}^{p} := (\nu_{1}\times
\hat{q}_{\varphi})^{*}\circ (t_{-s_{\beta}}\times
t_{-s_{\alpha}})^{*}(\mycal{P}_{|X^{p}\times_{U^{p}} X^{p}}),
\]
Note that the composite map $(t_{-s_{\beta}}\times
t_{-s_{\alpha}})\circ (\nu_{1}\times \hat{q}_{\varphi}) : Y^{p}\times
S^{p} \to X^{p}\times_{U^{p}} X^{p}$ factors through a small
resolution $\widehat{X^{p}\underset{U^{p}}{\times} X^{p}} \to
X^{p}\times_{U^{p}} X^{p}$. This implies that the  derived
pullback of the torsion free sheaf $\mycal{P}_{|X^{p}\times_{U^{p}}
X^{p}}$ via $(t_{-s_{\beta}}\times
t_{-s_{\alpha}})\circ (\nu_{1}\times \hat{q}_{\varphi})$ is equal to the
pullback as a coherent sheaf. Indeed, the pullback of
$\mycal{P}_{|X^{p}\times_{U^{p}}
X^{p}}$ to the small resolution is locally free and hence it will
pullback to a coherent sheaf $\mycal{P}^{p}$ (in fact a line bundle)
$Y^{p}\times_{X_{\varphi}^{p}} S^{p}$.

The next step is to observe that without a loss of generality we may
assume that the the gerbes $\Xl{\alpha}{\beta}^{p} \to X_{\beta}^{p}$
and $\Xe{\beta}{\alpha}^{p} \to X_{\alpha}^{p}$ are
trivial. Explicitly this means that we can find line bundles
$\mycal{Q}^{\, p} \in \op{Pic}(Y^{p})$ and $\psi^{p} \in
\op{Pic}(S^{p})$ on $Y^{p}$ and $S^{p}$ respectively, together with
isomorphisms
\begin{align}
\mycal{P}_{\ot|Y^{p}\times_{X_{\varphi}^{p}} Y^{p}} & \cong
p_{1}^{*}\mycal{Q}^{\, p} \otimes p_{2}^{*}(\mycal{Q}^{\, p})^{-1}
\label{eq-ltriv}\\ \Phi_{\beta|S^{p}\times_{X_{\alpha}^{p}} S^{p}} & =
p_{1}^{*}(\psi^{p})\otimes p_{2}^{*}(\psi^{p})^{-1}, \label{eq-etriv}
\end{align}
satisfying the obvious cocycle conditions.

We will work with particular trivializations of
$\Xl{\alpha}{\beta}^{p} \to X_{\beta}^{p}$ and
$\Xe{\beta}{\alpha}^{p} \to X_{\alpha}^{p}$ which we build in terms of
the sections $s_{\varphi}$ and $s_{\beta}$ respectively.
To construct $\mycal{Q}^{\, p}$ we use the rational map
\[
\xymatrix@1@C+55pt{Y^{p} \ar[r]^-{\nu} &
X_{\varphi}^{p}\times_{U^{p}} X_{\beta}^{p}
\ar@{-->}[r]^-{t_{-s_{\varphi}}\times q_{\beta}^{m}} &
X^{p}\times_{U^{p}} X^{p}.
}
\]
Note that $(t_{-s_{\varphi}}\times q_{\beta}^{m})\circ \nu$ is
a morphism from the complement of the exceptional curve $n^{p}$ for
the small resolution $\nu$ to the complement of the unique
singular point in $X^{p}\times_{U^{p}} X^{p}$. Since the standard
Poincare sheaf $\mycal{P}$ is a rank one torsion free sheaf on
$X\times_{B} X$ which fails to be locally free only at the singular
points of $X\times_{B} X$, it follows that the pullback of $\mycal{P}$
via the map $(t_{-s_{\varphi}}\times q_{\beta}^{m})\circ \nu$
makes sense as a line bundle defined on $Y^{p}\setminus
n^{p}$. Combined with the observation that $Y^{p}$ is smooth and that
$n^{p} \subset Y^{p}$ is of codimension two, it follows that this
pullback extends to a unique line bundle $\mycal{Q}^{\, p}$ on all of
$Y^{p}$. We can now use the seesaw principle in the same way we did in
the proof of Theorem~\ref{theo-lifting-vs-coho} to show that the
isomorphism \eqref{eq-ltriv} exists and satisfies the cocycle
condition.

Similarly, to construct $\psi^{p}$ we note that the section
$s_{\beta}$ gives rise to a relative line bundle
$\Sigma_{\beta}\otimes {\mathcal O}_{X_{\beta}^{p}}(-m\cdot s_{\beta})
\in \mycal{P}ic^{0}(X_{\beta}^{p}/U_{p})$. Under the canonical
identification $\mycal{P}ic^{0}(X_{\beta}^{p}/U^{p}) = \mycal{X}^{p} =
\mycal{P}ic^{0}(X_{\varphi}^{p}/U^{p})$, this relative bundle
corresponds to a relative line bundle $\rho^{p}$ of degree zero along
the fibers of $X_{\varphi}^{p} \to U^{p}$ and hence to a globally
defined line bundle $\widetilde{\psi}^{p}$ on $X_{\varphi}^{p}$ which
restricts to $\rho^{p} \in \mycal{P}ic^{0}(X_{\varphi}^{p}/U^{p})$.
We normalize $\widetilde{\psi}^{p}$ by choosing a trivialization
$s_{\varphi}^{*}\widetilde{\psi}^{p} \cong {\mathcal O}_{U^{p}}$. Now
the argument used in Theorem~\ref{theo-extension-vs-coho} implies that
on $(S^{p}\times_{X_{\alpha}^{p}} S^{p})_{|U_{p}\setminus \{p\}}$ we
can find an isomorphism $\Phi_{\beta} \cong
p_{1}^{*}\widetilde{\psi}^{p}\otimes
p_{2}^{*}(\widetilde{\psi}^{p})^{-1}$ which (after possibly rescaling
the trivialization $s_{\varphi}^{*}\widetilde{\psi}^{p} \cong
{\mathcal O}_{U^{p}}$) will also satisfy the cocycle condition. Now
since $S$ is an atlas for the gerbe $\Xe{\beta}{\alpha}$ we conclude
that $\widetilde{\psi}^{p}$ on $(S^{p}\times_{X_{\alpha}^{p}}
S^{p})_{|U_{p}\setminus \{p\}}$ extends to a unique line bundle
$\psi^{p}$ on $S^{p}\times_{X_{\alpha}^{p}} S^{p}$ equipped with an
isomorphism \eqref{eq-etriv} satisfying the cocycle condition.

We now define the coherent sheaf $\Pi^{p} \in
\op{Coh}(Y^{p}\times_{X_{\varphi}^{p}} S^{p})$ by setting
\[
\Pi^{p} := \mycal{P}^{p}\otimes p_{Y}^{*}(\mycal{Q}^{\,
p})^{-1}\otimes p_{S}^{*}(\psi^{p})^{-1},
\]
where $p_{Y} : Y \times S \to Y$ and $p_{S} : Y \times S \to  S$
denote the natural projections.

\

\medskip

\noindent
With this notation we now have

\begin{lem} \label{lem-glue} For any $p \in B\setminus B^{o}$ there is
an isomorphism
\[
\Pi^{o}_{|(Y^{o}\times_{X_{\varphi}^{o}} S^{o})\cap
(Y^{p}\times_{X_{\varphi}^{p}} S^{p})} \cong
\Pi^{p}_{|(Y^{o}\times_{X_{\varphi}^{o}} S^{o})\cap
(Y^{p}\times_{X_{\varphi}^{p}} S^{p})}.
\]
In particular the sheaves $\{\Pi^{p}\}_{p \in B\setminus B^{o}}$ glue
to the sheaf $\Pi^{o}$ to yield a globally defined analytic coherent
sheaf $\Pi$ on $Y\times S$.
\end{lem}
{\bf Proof.} We have to show that $\Pi^{p}$ is naturally isomorphic to
the trivial line bundle on $(Y^{o}\times_{X_{\varphi}^{o}} S^{o})\cap
(Y^{p}\times_{X_{\varphi}^{p}} S^{p})$.

Write $U^{po}$ for the punctured disc $U^{p}\setminus \{p\}$ and for
any space or stack $Z \to B$ write $Z^{po}$ for the fiber product
$Z\times_{B} U^{po}$.  Using the isomorphisms $t_{-s_{\varphi}}$ and
$t_{-s_{\beta}}$ we can now identify $\pi_{\varphi} : X_{\varphi}^{po}
\to U^{po}$ and $\pi_{\beta} : X_{\beta}^{po} \to U^{po}$ with the
smooth elliptic fibration $\pi : X^{po} \to U^{po}$. Under these
identifications the intersection $(Y^{o}\times_{X_{\varphi}^{o}}
S^{o})\cap (Y^{p}\times_{X_{\varphi}^{p}} S^{p})$ gets identified with
the fiber product $X^{po}\times_{U^{po}} X^{po}$. Using the same
trivializations to recast $\mycal{P}^{p}$, $\mycal{Q}^{\, p}$ and
$\psi^{p}$ as sheaves on  $X^{po}\times_{U^{po}}
X^{po}$ we get that $\mycal{P}^{p}$ becomes
$p_{1,m\cdot 2}^{*}\mycal{P}$, $\mycal{Q}^{\, p}$ becomes
$p_{m\cdot 1,2}^{*}\mycal{P}$ and $\psi^{p}$ becomes
${\mathcal O}$. In particular we get that $\Pi^{p}$
corresponds to $p_{1,m\cdot 2}\mycal{P}\otimes p_{m\cdot
1,2}^{*}\mycal{P}^{-1}$ and so it suffices to check that
\[
p_{1,m\cdot 2}\mycal{P}\otimes p_{m\cdot
1,2}^{*}\mycal{P}^{-1} \cong {\mathcal O}_{X^{po}\times_{U^{po}} X^{po}}.
\]
This however follows immediately from the universal property of
$\mycal{P}$ on $X\times_{B} X$ and the seesaw theorem. \ \hfill $\Box$

\

The sheaf $\Pi$ gives rise to a well defined integral transform
\begin{equation} \label{eq-F}
\xymatrix@R-25pt{
D^{b}(Y) \ar[r]^{\bF} & D^{b}(S) \\
L \ar[r] & p_{S*}(p_{Y}^{*}L\otimes \Pi)
}
\end{equation}
between the derived categories of analytic coherent sheaves on $Y$ and
$S$ respectively. If in addition $\alpha \in \TSh(X) \subset
\TSh_{an}(X)$ is also algebraic, then the atlases $Y$ and $S$ of the
lifting and extension gerbes are proper separated algebraic spaces and
so we can invoke the GAGA theorem \cite[Corollary~7.15]{artin-a2} to
conclude that $\Pi$ is an algebraic coherent sheaf on $Y\times
S$. This implies that in the algebraic context $\bF$ makes sense as an
integral transform between algebraic coherent sheaves.

We now have the following

\begin{prop} \label{prop-rel} The functor
$\bF : D^{b}(Y) \to D^{b}(S)$ defined by \eqref{eq-F} maps the descent
 data for the lifting presentation to the descent data for the
 extension presentation and thus defines a functor $\FM :
 D^{b}_{1}(\Xc{\alpha}{\beta}) \to D^{b}_{1}(\Xc{\beta}{\alpha})$.
\end{prop}
{\bf Proof.} Let ${\mathcal L}$ be an object in
$D^{b}_{1}(\Xc{\alpha}{\beta})$ represented by descent datum
$(L,\boldsymbol{f})$ for the presentation $\Xl{\alpha}{\beta}$. In
other words, $L$ is an object in $D^{b}(Y)$ and $\boldsymbol{f} :
p_{1}^{*}L \to  p_{2}^{*}L\otimes \mycal{P}_{\ot}$ is a
quasi-isomorphism on $Y\times_{X_{\beta}} Y$ satisfying the cocycle
condition on $Y\times_{X_{\beta}} Y\times_{X_{\beta}} Y$.

Consider the object $\bF L \in D^{b}(S)$. To prove the proposition we
need to construct a quasi-isomorphism $\boldsymbol{g} : p_{1}^{*}\bF L
\to p_{2}^{*}\bF L\otimes \Phi_{\beta}^{-1}$ on $S\times_{X_{\alpha}} S$
which depends functorially on $\boldsymbol{f}$ and satisfies the
cocycle condition on $S\times_{X_{\alpha}} S\times_{X_{\alpha}} S$.

Let $\Gamma := Y\times_{X_{\varphi}} S$ and let $p_{Y} : \Gamma \to Y$
and $p_{S} : \Gamma \to S$ denote the natural projections. Then $\bF L
= p_{S*}(p_{Y}^{*} L\otimes \Pi)$, and so our problem boils down to
finding a quasi-isomorphism
\[
\boldsymbol{g} : p_{1}^{*}p_{S*}(p_{Y}^{*}L\otimes \Pi)
\stackrel{\cong}{\to} p_{2}^{*}p_{S*}(p_{Y}^{*}L\otimes \Pi)\otimes
\Phi_{\beta}^{-1}
\]
in $D^{b}(S\times_{X_{\alpha}} S)$, which depends functorially on
$\boldsymbol{f}$. In other words we want to compare the objects
$p_{1}^{*}p_{S*}(p_{Y}^{*}L\otimes \Pi)$ and
$p_{2}^{*}p_{S*}(p_{Y}^{*}L\otimes \Pi)$ on $S\times_{X_{\alpha}} S$.

Since we have an obvious commutative diagram
\begin{equation} \label{eq-diagGamma}
\xymatrix@M+5pt{ Y\times_{X_{\beta}} Y
\ar@<.5ex>[r]^-{p_{1}} \ar@<-.5ex>[r]_-{p_{2}} & Y \\
\Gamma\times_{X_{\beta}} \Gamma \ar[u]^-{p_{Y}\times p_{Y}}
\ar@<.5ex>[r]^-{\op{p}_{1}} \ar@<-.5ex>[r]_-{\op{p}_{2}}
\ar[d]_-{p_{S}\times p_{S}} &
\Gamma \ar[d]^-{p_{S}} \ar[u]_{p_{Y}} \\
S\times_{B} S
\ar@<.5ex>[r]^-{\op{pr}_{1}}
\ar@<-.5ex>[r]_-{\op{pr}_{2}} &
S  \\
S\times_{X_{\alpha}} S \ar@{^{(}->}[u]^-{\iota}  \ar@<.5ex>[r]^-{p_{1}}
\ar@<-.5ex>[r]_-{p_{2}} & S \ar@{=}[u]
}
\end{equation}
we see that equivalently we can compare the objects
$\iota^{*}\op{pr}_{1}^{*}p_{S*}(p_{Y}^{*}L\otimes \Pi)$ and \linebreak
$\iota^{*}\op{pr}_{2}^{*}p_{S*}(p_{Y}^{*}L\otimes \Pi)$. To compute
these objects, we would like to perform a base change in the
commutative squares
\begin{equation} \label{eq-squares}
\xymatrix{\Gamma_{X_{\beta}}\Gamma \ar[r]^{\op{p}_{1}} \ar[d]_-{p_{S}\times
p_{S}} &  \Gamma \ar[d]^-{p_{S}} \\
S\times_{B} S \ar[r]_-{\op{pr}_{1}} & S} \qquad
\xymatrix{\Gamma_{X_{\beta}}\Gamma \ar[r]^{\op{p}_{2}} \ar[d]_-{p_{S}\times
p_{S}} &  \Gamma \ar[d]^-{p_{S}} \\
S\times_{B} S \ar[r]_-{\op{pr}_{2}} & S}.
\end{equation}
Unfortunately these squares are not cartesian and so we do not have
the base change property on the nose. However we have the following useful

\begin{rem} \label{rem-cohofp} Let $X$, $Y$, $Z$ and $T$ be analytic (or
algebraic) spaces and let
\[
\xymatrix{ Z \ar[r]^-{p} \ar[d]_-{q} & X \ar[d]^-{f} \\
Y \ar[r]_-{g} & T
}
\]
be a commutative square of proper maps. Suppose further that the
natural map \linebreak $u : Z  \to Y\times_{T} X$ satisfies
$u_{*}{\mathcal O}_{Z} = {\mathcal O}_{Y\times_{T} X}$. Then
for every $F \in D^{b}(X)$ we have a base change identification
$q_{*}p^{*} F = g^{*}f_{*} F$ in $D^{b}(Y)$. Indeed, we can complete
the above square to a commutative diagram
\[
\xymatrix{
Z \ar[drr]^-{p} \ar[ddr]_-{q} \ar[dr]_-{u} & & \\
& Y\times_{T} X \ar[r]_-{\bar{p}} \ar[d]^-{\bar{q}} & X \ar[d]^-{f} \\
& Y \ar[r]_{g} & T
}
\]
Since $\bar{q}_{*}\bar{p}^{*} F = g^{*}f_{*} F$ we get
\begin{align*}
q_{*}p^{*} F & = \bar{q}_{*}u_{*} u^{*}
\bar{p}^{*} F  \tag{transitivity} \\
& = \bar{q}_{*}(\bar{p}^{*} F \otimes u_{*} {\mathcal
O}_{Z}) \tag{projection formula} \\
& = \bar{q}_{*}\bar{p}^{*} F \tag{assumption} \\
& = g^{*}f_{*} F \tag{base change}.
\end{align*}
\end{rem}

\

\noindent
In view of the previous remark we will be able to treat the squares
\eqref{eq-squares} as base change squares if we can show that for
the maps
\[
u_{1} : \Gamma\times_{X_{\beta}} \Gamma \to (S\times_{B} S)
\times_{p_{1}, S, p_{S}} \Gamma \qquad \text{and} \qquad
u_{2} : \Gamma\times_{X_{\beta}} \Gamma \to (S\times_{B} S)
\times_{p_{2}, S, p_{S}} \Gamma
\]
we have $u_{1*} {\mathcal O} = {\mathcal O}$ and
$u_{2*} {\mathcal O} = {\mathcal O}$.

To check this, note that
\[
\begin{split}
\Gamma\times_{X_{\beta}} \Gamma & = (S\times_{B}
S)\times_{X_{\varphi}\times_{B} X_{\varphi}} (Y\times_{X_{\beta}} Y)
\\
(S\times_{B} S) \times_{p_{1}, S, p_{S}} \Gamma & = (S\times_{B}
S)\times_{\varepsilon\circ p_{1}, X_{\varphi}, \nu_{1}} Y \\
(S\times_{B} S) \times_{p_{2}, S, p_{S}} \Gamma & = (S\times_{B}
S)\times_{\varepsilon\circ p_{2}, X_{\varphi}, \nu_{2}} Y,
\end{split}
\]
and so $u_{i}$, $i = 1,2$ are given explicitly by
\[
\xymatrix@R-25pt{
u_{i} : & (S\times_{B}
S)\times_{X_{\varphi}\times_{B} X_{\varphi}} (Y\times_{X_{\beta}} Y)
\ar[r] & (S\times_{B}
S)\times_{\varepsilon\circ p_{i}, X_{\varphi}, \nu_{1}} Y \\
& ((a_{1},a_{2}),(b_{1},b_{2})) \ar[r] & ((a_{1},a_{2}),b_{i}).
}
\]
This implies in particular that the maps $u_{i}$
fit in the cartesian squares
\[
\xymatrix@C+50pt{ (S\times_{B} S)\underset{X_{\varphi}\times_{B}
X_{\varphi}}{\times} (Y\times_{X_{\beta}} Y) \ar[d]^-{u_{1}}
\ar[r]^-{p_{2}\circ p_{Y\times_{X_{\beta}}}} & Y \ar[d]^-{\nu} \\
(S\times_{B} S)\times_{\varepsilon\circ p_{1}, X_{\varphi}, \nu_{1}} Y
\ar[r]_-{(\varepsilon\circ p_{2}\circ \op{pr}_{S\times_{B} S})\times
(\nu_{2}\circ \op{pr}_{Y})} & X_{\varphi}\times_{B} X_{\beta} }
\]
and
\[
\xymatrix@C+50pt{ (S\times_{B}
S)\underset{X_{\varphi}\times_{B} X_{\varphi}}{\times} (Y\times_{X_{\beta}} Y)
\ar[d]^-{u_{2}} \ar[r]^-{p_{1}\circ p_{Y\times_{X_{\beta}}}}
& Y \ar[d]^-{\nu} \\
(S\times_{B}
S)\times_{\varepsilon\circ p_{2}, X_{\varphi}, \nu_{1}} Y
\ar[r]_-{(\varepsilon\circ p_{1}\circ \op{pr}_{S\times_{B} S})\times
(\nu_{2}\circ \op{pr}_{Y})} & X_{\varphi}\times_{B} X_{\beta}
}
\]
where $\nu : Y \to X_{\varphi}\times_{B} X_{\beta}$ is the small
resolution map defining $Y$.

Since the pullback of ${\mathcal O}$ by any morphism is again
${\mathcal O}$, it suffices to check that there is a canonical
isomorphism $\nu_{*}{\mathcal O}_{Y} = {\mathcal
O}_{X_{\varphi}\times_{B} X_{\beta}}$. This is obvious by the cohomology
and base change theorem. Thus $u_{i*}{\mathcal O} =
{\mathcal O}$ and so the squares \eqref{eq-squares} have the base
change property.

In particular,  in $D^{b}(S\times_{B} S)$ we get identifications:
\[
\op{pr}_{i}^{*}p_{S*}(p_{Y}^{*}L\otimes \Pi) = (p_{S}\times
p_{S})_{*}\op{p}_{i}^{*}(p_{Y}^{*}L\otimes \Pi),
\]
for all $L \in D^{b}(Y)$ and  for $i = 1,2$. Furthermore since
\[
(S\times_{X_{\alpha}} S)\underset{S\times_{B} S}{\times}
(\Gamma\times_{X_{\beta}} \Gamma) =
\Gamma\underset{X_{\alpha}\times_{B} X_{\beta}}{\times} \Gamma,
\]
we can identify $\iota^{*}(p_{S}\times
p_{S})_{*}\op{p}_{i}^{*}(p_{Y}^{*}L\otimes \Pi)$ with
$p_{S\times_{X_{\alpha}}S*}p_{i}^{*}(p_{Y}^{*}L\otimes \Pi)$ where
$p_{S\times_{X_{\alpha}}S} : \Gamma\times_{X_{\alpha}\times_{B}
X_{\beta}} \Gamma \to S\times_{X_{\alpha}} S$ is the natural
projection.

Now, using the commutativity of the top double square in
\eqref{eq-diagGamma} we get
\[
\begin{split}
(p_{S\times_{X_{\alpha}} S})_{*}p_{1}^{*}(p_{Y}^{*}L\otimes \Pi) & =
(p_{S\times_{X_{\alpha}} S})_{*}(((p_{Y}\times
p_{Y})^{*}p_{1}^{*}L)\otimes (p_{1}^{*}\Pi)), \text{ and } \\
(p_{S\times_{X_{\alpha}} S})_{*}p_{2}^{*}(p_{Y}^{*}L\otimes \Pi) & =
(p_{S\times_{X_{\alpha}} S})_{*}(((p_{Y}\times
p_{Y})^{*}p_{2}^{*}L)\otimes (p_{2}^{*}\Pi)).
\end{split}
\]
If in addition $L \in D^{b}(Y)$ is part of a descend datum
$(L,\boldsymbol{f})$ defining an object in
$D^{b}_{1}(\Xc{\alpha}{\beta})$, we can use $\boldsymbol{f}$ to obtain
an identification
\[
(p_{S\times_{X_{\alpha}} S})_{*}p_{1}^{*}(p_{Y}^{*}L\otimes \Pi) =
(p_{S\times_{X_{\alpha}} S})_{*}((p_{Y}\times
p_{Y})^{*}p_{2}^{*}L\otimes ((p_{Y}\times p_{Y})^{*}
\mycal{P}_{\ot})\otimes (p_{1}^{*}\Pi)).
\]
Also
\[
((p_{S\times_{X_{\alpha}} S})_{*}p_{2}^{*}(p_{Y}^{*}L\otimes
\Pi))\otimes \Phi_{\beta}^{-1} = (p_{S\times_{X_{\alpha}}
S})_{*}(((p_{Y}\times p_{Y})^{*}p_{2}^{*}L)\otimes
((p_{S\times_{X_{\alpha}} S}^{*} \Phi_{\beta}^{-1}) \otimes
(p_{2}^{*}\Pi)),
\]
and so, in order to get the desired isomorphism $\boldsymbol{g} :
p_{1}^{*}\bF L \to p_{2}^{*}\bF L \otimes \Phi_{\beta}^{-1}$ in
$D^{b}(S\times_{X_{\alpha}} S)$ we only have to construct a canonical
identification
\begin{equation} \label{eq-identify}
p_{1}^{*} \Pi \otimes (p_{Y}\times p_{Y})^{*} \mycal{P}_{\ot} =
p_{2}^{*} \Pi \otimes p_{S\times_{X_{\alpha}} S}^{*} \Phi_{\beta}^{-1}
\end{equation}
of coherent sheaves on $\Gamma\times_{X_{\alpha}\times_{B} X_{\beta}}
\Gamma$.

We will construct the desired map \eqref{eq-identify} by gluing
some locally defined but canonical identifications. Note that at this
point we have completely eliminated the derived category from the
picture. In particular, we are left with a question about sheaves, not
complexes, and so gluing is a relatively simple matter.

\

\medskip

\noindent
{\bf (i)} \quad Over the part of $\Gamma\times_{X_{\alpha}\times_{B}
X_{\beta}} \Gamma$ sitting over $B^{o} \subset B$, the identification
\eqref{eq-identify} is just the one established in
section~\ref{sec-duality}. Indeed, by definition $\Pi_{|\Gamma^{o}} =
{\mathcal O}_{\Gamma^{o}}$ and so constructing \eqref{eq-identify}
over $B^{o}$ becomes equivalent to constructing a canonical
identification \eqref{eq-iso2}. Such a construction was carried out in
the proof of Proposition~\ref{prop-gerbyFM}.

\

\medskip

\noindent
{\bf (ii)} \quad Over the part of $\Gamma\times_{X_{\alpha}\times_{B}
X_{\beta}} \Gamma$ sitting over $U^{p} \subset B$, $p \in B\setminus
B^{o}$, we can use again the line bundles $\mycal{Q}^{p} \to Y^{p}$ and
$\psi^{p} \to S^{p}$ appearing in the construction of $\Pi^{p} =
\Pi_{|\Gamma^{p}}$ to trivialize our gerbes. Recall that
$\mycal{Q}^{p}$ and ${\psi}^{p}$ come equipped with the
natural isomorphisms \eqref{eq-ltriv} and
\eqref{eq-etriv}. Furthermore $\Pi^{p} = \mycal{P}^{p}\otimes
p_{Y}^{*}(\mycal{Q}^{p})^{-1} \otimes p_{S}^{*}(\psi^{p})^{-1}$ and so
establishing \eqref{eq-identify} over $\Gamma^{p}
\times_{X_{\alpha}^{p}\times_{U^{p}} X_{\beta}^{p}} \Gamma^{p}$
reduces to constructing an identification:
\[
\xymatrix@R-15pt@M+5pt{
p_{1}^{*}\mycal{P}^{p}\otimes p_{1}^{*}p_{Y}^{*}(\mycal{Q}^{p})^{-1}\otimes
p_{1}^{*}p_{S}^{*}(\psi^{p})^{-1} \otimes p_{1}^{*}p_{Y}^{*}\mycal{Q}^{p}
\otimes p_{2}^{*}p_{Y}^{*} (\mycal{Q}^{p})^{-1} \ar@{=}[d] \\
p_{2}^{*}\mycal{P}^{p}\otimes p_{2}^{*}p_{Y}^{*}(\mycal{Q}^{p})^{-1}\otimes
p_{2}^{*}p_{S}^{*}(\psi^{p})^{-1} \otimes p_{1}^{*}p_{S}^{*}(\psi^{p})^{-1}
\otimes p_{2}^{*}p_{S}^{*} \psi^{p}
}
\]
or equivalently, after the obvious cancellations, an identification
\[
p_{1}^{*}\mycal{P}^{p} = p_{2}^{*}\mycal{P}^{p}.
\]
However $\mycal{P}^{p} \in \op{Coh}(\Gamma^{p})$ was defined as a
pullback of a sheaf on $X_{\alpha}^{p}\times_{U^{p}} X_{\beta}^{p}$
and so we get a canonical identification $p_{1}^{*}\mycal{P}^{p} =
p_{2}^{*}\mycal{P}^{p}$ on $\Gamma^{p}
\times_{X_{\alpha}^{p}\times_{U^{p}} X_{\beta}^{p}} \Gamma^{p}$. This
yields the desired canonical identification \eqref{eq-identify} over
the part of $\Gamma\times_{X_{\alpha}\times_{B} X_{\beta}} \Gamma$
sitting over $U^{p}$

\

\medskip

Finally it only remains to observe that the isomorphisms chosen in
\eqref{eq-ltriv} and \eqref{eq-etriv} were the ones used in the proof of
Lemma~\ref{lem-glue} to glue $\Pi^{p}$ and $\Pi^{o}$ on the
overlap $\Gamma^{p}\cap \Gamma^{o}$. Therefore the identifications in
items {\bf (i)} and {\bf (ii)} above glue on the overlaps
$(\Gamma\times_{X_{\alpha}\times_{B} X_{\beta}} \Gamma)\times_{B}
U^{po}$ and so we have found our global identification
\eqref{eq-identify}. This finishes the proof of the \linebreak
proposition.
\hfill $\Box$

\

\bigskip

\noindent
We are now ready to complete the

\medskip

\noindent
{\bf Proof of Theorem~\ref{Main-K3}.} The only thing left to show is
that the gerby Fourier-Mukai transform $\FM :
D^{b}_{1}(\Xc{\alpha}{\beta}) \to D^{b}_{1}(\Xc{\beta}{\alpha})$
constructed in Proposition~\ref{prop-rel} is an equivalence of
categories. We will again use the criterion of Bondal-Orlov and
Bridgeland applied to the spanning class $\Omega$ of gerby points in
$\Xc{\alpha}{\beta}$ described in Claim~\ref{claim-spanning}. As
before we need to show that $\FM$ intertwines the Serre functors on
the sheaves ${\mathcal O}_{x} \in \Omega$ and that $\FM$ satisfies the
orthogonality property
\[
\FM : \op{Hom}_{D^{b}_{1}(\Xc{\alpha}{\beta})}^{i}
({{\mathcal O}}_{x_1},{{\mathcal O}}_{x_2}) \widetilde{\to}
\op{Hom}_{D^{b}_{1}(\Xc{\beta}{\alpha})}^{i}
(\bF{{\mathcal O}}_{x_1},\bF{{\mathcal O}}_{x_2}), \quad \text{for
all } i \in {\mathbb Z}, x_1, x_2 \in X_{\beta}.
\]
Since by definition the Fourier-Mukai image $\FM {\mathcal O}_{x}$ is
supported on the fiber $(\Xc{\beta}{\alpha})_{b}$ of
$\Xc{\beta}{\alpha}$ over the point $b = \pi_{\beta}(x) \in B$, it
suffices to check the intertwining and orthogonality properties of
$\FM$ locally in the base $B$.

Over $B^{o}$ these properties were established in
Claims~\ref{claim-ortho} and \ref{claim-intertwine}.  To check the
properties for the parts of our gerbes sitting over $U^{p} \subset B$
we note that the proof of Lemma~\ref{lem-glue} shows that over $U^{p}$
the functor $\FM$ fits in a commutative diagram of functors
\[
\xymatrix@C+25pt{
D^{b}_{1}(\Xc{\alpha}{\beta}^{p}) \ar[r]^{\FM} &
D^{b}_{-1}(\Xc{\beta}{\alpha}^{p}) \\
D^{b}(X^{p}_{\beta}) \ar[u]^-{\nu_{1}^{*}(\bullet)\otimes
\mycal{Q}^{p}} & D^{b}(X_{\alpha}^{p})
\ar[u]_-{\nu_{2}^{*}(\bullet)\otimes (\psi^{p})^{-1}} \\
D^{b}(X^{p}) \ar[u]^-{t_{-s_{\beta}}^{*}}
\ar[r]_-{p_{2*}(p_{1}^{*}(\bullet)\otimes \mycal{P})} & D^{b}(X^{p})
\ar[u]_-{t_{-s_{\alpha}}^{*}}
}
\]
where the vertical arrows are equivalences. However the bottom arrow
is the usual integral transform with respect to the Poincare
sheaf on an elliptic surface having at most $I_{1}$ fibers. Such
a transform  is an equivalence, e.g. by
\cite{bridgeland-maciocia}. Finally, the functor
$(\nu_{1}^{*}(\bullet)\otimes \mycal{Q}^{p})\circ t_{-s_{\beta}}^{*}$
transforms a structure sheaf of a point $x \in X^{p}$ into a
sheaf in the spanning class $\Omega^{p}$ for $\Xc{\alpha}{\beta}^{p}$
and clearly every sheaf in $\Omega^{p}$ is obtained this way. This
implies that $\FM$ has the orthogonality and intertwining properties
for sheaves in $\Omega^{p}$. The theorem is proven, \ \hfill $\Box$

\

\bigskip

\noindent
>From the statement of Theorem~\ref{Main-K3} one can derive a whole
sequence of new cases of C\u{a}ld\u{a}raru's conjecture. Indeed,
suppose $X$ is an elliptic K3 surface whose singular fibers are of
type $I_{1}$ only. Note that for any element $\alpha \in \TSh(X) =
\op{Hom}(\Tr_{X},{\mathbb Q}/{\mathbb Z})$ we have a natural Hodge
isometry $\Tr_{X_{\alpha}} \cong \ker(\alpha)$ induced by the isogeny
of $X_{\alpha}$ and $X$. In terms of the identifications $\TSh(X) =
Br'(X) = \op{Hom}(\Tr_{X},{\mathbb Q}/{\mathbb Z})$ and
$Br'(X_{\alpha}) = \op{Hom}(\ker(\alpha),{\mathbb Q}/{\mathbb Z})$,
the surjective map $T_{alpha} : \TSh(X) \to Br'(X_{\alpha})$ sends a
homomorphism $a : \Tr_{X} \to {\mathbb Q}/{\mathbb Z}$ to its
restriction $a_{|\ker(\alpha)} : \ker(\alpha) \to {\mathbb Q}/{\mathbb
Z}$. Now we have:

\

\medskip

\begin{cor} \label{cor-andrei} Let $X$ be an elliptic K3 surface
whose singular fibers are of type $I_{1}$ only. Let
$\alpha, a \in \TSh(X) = \op{Hom}(\Tr_{X},{\mathbb Q}/{\mathbb Z})$
and let $(b,\beta) \in \TSh(X)^{\times 2}$ be in the $SL(2,{\mathbb
Z})$ orbit of $(\alpha,a)$. Then
\begin{itemize}
\item[(a)] $\ker\left[\ker(\alpha)\stackrel{a}{\to} {\mathbb
Q}/{\mathbb Z}\right]$ and $\ker\left[\ker(\beta)\stackrel{b}{\to}
{\mathbb Q}/{\mathbb Z}\right]$ are Hodge isometric;
\item[(b)] $D^{b}_{1}(\Xc{a}{\alpha})$ and $D^{b}_{1}(\Xc{b}{\beta})$
are equivalent.
\end{itemize}
\end{cor}
{\bf Proof.} Part (b) follows from the fact that for every
$(a,\alpha)$ we have equivalences $D^{b}_{1}(\Xc{a}{\alpha}) \cong
D^{b}_{1}(\Xc{-\alpha}{a})$ (by Theorem~\ref{Main-K3}) and
$D^{b}_{1}(\Xc{a}{\alpha}) \cong D^{b}_{1}(\Xc{a+\alpha}{\alpha})$
(since $T_{\alpha}(\alpha) = 0$).

For part (a) observe that by our identifications we have isometries of
Hodge lattices $\ker\left[\ker(\alpha)\stackrel{a}{\to} {\mathbb
Q}/{\mathbb Z}\right] = \ker(a)\cap\ker(\alpha)$ and
$\ker\left[\ker(\beta)\stackrel{b}{\to} {\mathbb Q}/{\mathbb
Z}\right] = \ker(b)\cap \ker(\beta)$. Since
$(a,\alpha) : \Tr_{X} \to ({\mathbb Q}/{\mathbb Z})^{2}$
and $(b,\beta) : \Tr_{X} \to ({\mathbb Q}/{\mathbb Z})^{2}$ differ by
an element of $\op{SL}(2,{\mathbb Z})$ acting on $({\mathbb
Q}/{\mathbb Z})^{2}$, it follows that $\ker(a)\cap \ker(\alpha) =
\ker(b)\cap \ker(\beta)$ as sublattices in $\Tr_{X}$. \ \hfill $\Box$

\section{Modified $T$-duality and the SYZ conjecture} \label{sec-SYZ}

The celebrated work of Strominger, Yau and Zaslow \cite{SYZ}
interprets mirror symmetry of Calabi-Yaus in terms of special
Lagrangian (SLAG) torus fibrations. If a CY manifold $X$ (with ``large
complex struture") has mirror $X'$, \cite{SYZ} conjecture the
existence of fibrations $\pi: X \to B$ and $\pi': X' \to B$ whose
generic fibers are SLAG tori dual to each other: each parameterizes
$U(1)$ flat connections on the other.  In particular, each of these
fibrations admits a SLAG zero-section, corresponding to the trivial
connection on the dual fibers. The analogy with the situation
considered in the main part of our work is clear: the SLAG torus fibrations
replace the elliptic fibrations, and mirror symmetry (interchanging
D-branes of type B with D-branes of type A) replaces the Fourier-Mukai
transform (which interchanges vector bundles with spectral data).

In this context, the analogue of our gerbes and the Brauer group is
given by the "B-fields" $\alpha \in H^2(X, {\mathbb R}/{\mathbb
Z})$. On the other hand, the SLAG analogue of the Tate-Shafarevich
group of $X'$ is given by $H^1(B,X')$, which over the locus where
$\pi, \pi'$ are smooth can be identified with $H^1(B,
R^1\pi_*({\mathbb R}/{\mathbb Z}))$.  As in the proof of lemma
\ref{lem-Talpha}, $H^2(X, {\mathbb R}/{\mathbb Z})$ is related via a
Leray spectral sequence to the three groups $H^i(B,
R^{2-i}\pi_*({\mathbb R}/{\mathbb Z}))$, for $i=1,2,3$. Now for $i=0$,
the local system $H^{0}(B,R^2\pi_*({\mathbb R}/{\mathbb Z}))$ can be
identified, over the locus where $\pi$ is smooth, with the group of
homotopy classes of sections of $X \to B$. Therefore, if the
fibration $X \to B$ is good in the sense of
\cite{gross-slag1,gross-slag2},  $H^0(B,
R^{2}\pi_*({\mathbb R}/{\mathbb Z}))$ should be thought of as the
analogue of the Mordell-Weil group of $X \to B$. Assume that the SLAG
fibration $X \to B$ is generic, in the sense that the local system
$R^{2}\pi_{*}({\mathbb R}/{\mathbb Z})$ has no global sections.

The Leray spectral sequence therefore gives a Brauer-to-Tate-Shafarevich
map:
$$H^2(X, {\mathbb R}/{\mathbb Z}) \to H^1(B, R^1\pi_*({\mathbb
R}/{\mathbb Z})).$$

We therefore may as well start with a pair of B-fields $\alpha \in
H^2(X, {\mathbb R}/{\mathbb Z}), \beta \in H^2(X', {\mathbb
R}/{\mathbb Z})$ on $X$ and $X'$ respectively. Since mirror symmetry
involves the B-field in an essential way, this suggests that the SYZ
conjecture, which is a SLAG analogue of the Fourier-Mukai transform
for elliptic fibrations, must be supplemented by an analogue of our
Theorem \ref{Main-K3}. Let $\mathcal M$ be the CY moduli space on which
mirror symmetry acts. The emerging picture is that $\mathcal M$ looks, at
least in some approximate sense, like an integrable system. The base
is a real submanifold ${\mathcal M}_{\mathbb R} \subset {\mathcal M}$. The
normal directions are parametrized by the B-fields $\alpha, \beta$ ,
so they form a torus isomorphic to $H^2(X, {\mathbb R}/{\mathbb Z})
\times H^2(X', {\mathbb R}/{\bf Z})$. The original SYZ conjecture
holds along ${\mathcal M}_{\mathbb R}$. As we move in normal directions,
$X$ becomes gerby along $B$-directions on $X$, while along
$B$-directions on $X'$, the SLAG fibration on $X$ loses its SLAG
zero-section. Mirror symmetry interchanges these two behaviors.

Indeed, the moduli space $\mathcal M$ parameterizes pairs $(X,\alpha)$ where
$X$ is a complex manifold together with a Calabi-Yau metric, and $\alpha
\in H^2(X, {\mathbb R}/{\mathbb Z})$ is a B-field on it. Mirror symmetry is an
involution $\MS: {\mathcal M} \to {\mathcal M}$, presumably defined in a
neighborhood of the ``large complex structure" point, taking $(X,\alpha)$
to $(X',\beta)$. We let  ${\mathcal M}_{\mathbb R} \subset {\mathcal M}$ denote the
locus where $\alpha=\beta=0$. It is a component of the fixed locus of the
antilinear involution which reverses the signs of $\alpha$ and $\beta$.
We denote a point of ${\mathcal M}_{\mathbb R}$ by the mirror pair
$X,X'$. Now a B-field $\beta \in H^2(X', {\mathbb R}/{\mathbb Z})$ on $X'$
determines a  point $(X',\beta)$ of $\mathcal M$, hence a mirror point
$X_{\beta} := \MS(X',\beta)$. For small $\beta$, this $X_{\beta}$ is a
deformation of $X$, so the additional B-field  $\alpha \in
H^2(X, {\mathbb R}/{\mathbb Z})$ on $X$
determines a corresponding B-field  $T_{\beta}(\alpha) \in H^2(X_{\beta},
{\mathbb R}/{\mathbb Z})$ on $X_{\beta}$. \newline

\noindent
\begin{con} \label{syzB}
\begin{itemize}
\item The SYZ picture holds (near the large complex
structure/large volume limit) on ${\mathcal M}_{\mathbb R}$.
\item For a B-field $B'$ on $X'$, the deformed Calabi-Yau $X_{B'}$ admits a
SLAG fibration (generally {\em without} a section) whose Jacobian (i.e.
double dual) is the original SLAG fibration (with section) on $X$.
\item Mirror symmetry preserves the integrable system
structure: for any pair $\alpha, \beta$, the mirror of $(X_{\beta},
T_{\beta}(\alpha))$ is $(X'_{\alpha}, T_{\alpha}(\beta))$.
\end{itemize}
\end{con}

We note that this modification of \cite{SYZ} is consistent with recent
interpretations in the literature (see \cite{hitchin-slag} and
references therein) of D-branes on Calabi-Yaus in the presence of a
B-field: A D-brane of type {\sf B} on $X$ is a coherent sheaf on the
gerbe given by the B-filed $\beta$, while a D-brane of type {\sf A} on
$X'$ is, roughly, a flat $U(1)$ connection on the restriction of the
gerbe $\alpha$ to a SLAG submanifold in $X'$. The
third part of this conjecture is the exact SLAG translation of Theorem
\ref{Main-K3}.

\

\bigskip

\bigskip

\

\appendix

\Appendix{(by D.Arinkin) Duality for representations of $1$-motives}
\label{appendix-motives}

In this appendix, we sketch a different approach to the Fourier-Mukai
transform for $\calO^\times$-gerbes over smooth genus one fibrations
(Theorem B). In this approach, Theorem B claims that the dual
commutative group stacks (of a certain type) have equivalent derived
categories of coherent sheaves. Let us review the duality for
commutative group stacks (sometimes called the duality for generalized
1-motives).

Recall that the dual $X^\vee$ of an abelian variety $X$ is the moduli
space of line bundles with zero first Chern class on
$X$. Equivalently, $X^\vee$ parametrizes the extensions of the
algebraic group $X$ by $\GM$. In this form, the definition immediately
generalizes to stacks: for a commutative group stack $\calX$, its dual
$\calX^\vee$ is the moduli stack of extensions of commutative group
stacks
\begin{equation*}
1\to\GM\to G\to\calX\to 0.
\end{equation*}
The sum of extensions defines a group operation on $\calX^\vee$;
actually, $\calX^\vee$ is naturally a commutative group stack.

\REMARK{For technical reasons, we use a slightly different definition
of the dual stack (Definition \ref{defdual}).  This allows to avoid
the discussion of short exact sequences of group stacks; also, the
group structure on $\calX^\vee$ seems somewhat more natural.}

Let $\calP\to\calX^\vee\times\calX$ be the universal
$\calX^\vee$-family of extensions of $\calX$ by $\GM$; in particular,
$\calP$ is a $\GM$-torsor on $\calX^\vee\times\calX$ (in fact, $\calP$
is a biextension of $\calX^\vee\times\calX$ by $\GM$). Notice that we
can also view $\calP$ as a $\calX$-family of extensions of
$\calX^\vee$ by $\GM$; this defines a morphism
$\calX\to(\calX^\vee)^\vee$. The main idea of the Fourier-Mukai
transform for commutative group stacks can be informally stated as
follows:

\begin{equation}
\label{goodstacks}
\txt{\begin{minipage}[c]{5in} For a ``good'' commutative group stack
$\calX$, the morphism $\calX\to(\calX^\vee)^\vee$ is an isomorphism,
and the Fourier-Mukai transform defined by $\calP_\CCC$ is an
equivalence $\FM:D^b(\calX)\to D^b(\calX^\vee)$.  Here $\calP_\CCC$ is
the line bundle on $\calX^\vee\times\calX$ associated to the
$\GM$-torsor $\calP$.
\end{minipage}}
\end{equation}

Now let us explain how Theorem B fits into the framework of the
duality for commutative group stacks. First, we notice that the
$\calO^\times$-gerbe $\Xc{\alpha}{0}$ over $X$ is a group stack. Then
we see that $\Xc{\alpha}{\beta}$ is a torsor over the group stack
$\Xc{\alpha}{0}$; more precisely, the gerbes constructed using the
lifting presentation and the extension presentation (from Section 3.1)
have a natural structure of $\Xc{\alpha}{0}$-torsors.

Torsors over a group stack can be thought of as extensions of $\ZZZ$
by this group stack; in this way, $\Xc{\alpha}{\beta}$ defines a
commutative group stack $\tXc{\alpha}{\beta}$ that fits into an
exact sequence
\begin{equation*}
0\to\Xc\alpha0\to\tXc{\alpha}{\beta}\to\ZZZ\to0.
\end{equation*}
The argument in section \ref{sec-duality} shows
that the constructions of the lifting presentation and
the extension presentation are dual, so $\tXc{\alpha}{\beta}$ and
$\tXc{-\beta}{\alpha}$ are dual commutative group stacks (provided that
we use the lifting presentation for one of the two stacks and the
extension presentation for the other). Moreover, these stacks are
``good'' in the sense of (\ref{goodstacks}), and so the Fourier-Mukai
transform gives an equivalence between $D^b(\tXc{\alpha}{\beta})$ and
$D^b(\tXc{-\beta}{\alpha})$. The Fourier-Mukai transform of Theorem B
is the restriction of this equivalence to direct summands in the
derived categories (see Section \ref{dualtorsors})

In the rest of the appendix, we discuss the notion of the dual of a
group stack (Section \ref{dualstacks}) and the special case when the
group stack is an extension of $\ZZZ$ (Section \ref{dualtorsors}). No
proofs are given, but most statements are almost obvious.

I learned about the duality for commutative group stacks from
A.~Beilinson, and I am deeply grateful to him for the explanation.

\subsection{Duality for commutative group stacks}
\label{dualstacks}

>From now on, the word `stack' means an algebraic stack locally of
finite type over a fixed base scheme $B$.  All results also have an
analytic version.

\DEFINITION{\label{defdual} For a commutative group stack $\calX$, the
dual stack $\calX^\vee$ parametrizes 1-morphisms of commutative group
stacks from $\calX$ to $B\GM$ (the classifying stack of $\GM$). Thus,
for a $B$-scheme $S$, the category $\calX^\vee(S)$ is the category of
1-morphisms of commutative group $S$-stacks
$\calX\times_BS\to B\GM\times S$. Notice that $\calX^\vee$ does not have
to be algebraic.}

\REMARK{For the definition to make sense, we need certain smallness
assumptions. Indeed, if $\calX$ and $\calY$ are stacks on a site $B$,
the 1-morphisms from $\calX$ to $\calY$ form a stack only if $\calX$,
$\calY$, and $B$ are small. However, this problem can be avoided if we
assume that $\calX$ is an algebraic stack which is locally of finite
type and replace the category of finitely presented $B$-schemes by an
equivalent small category.}

\EXAMPLE{If $\calX$ is an abelian scheme over $B$, then $\calX^\vee$
is the dual abelian scheme.}

\EXAMPLE{Let $\calX=G$ be an affine (or ind-affine) abelian group
(over $\CCC$).  Then $\calX^\vee$ is the classifying stack of the
Cartier dual of $G$. In particular, if $\calX=\ZZZ$, we have
$\calX^\vee=B\GM$.}

Another example is provided by the stacks $\tXc{\alpha}{\beta}$
(constructed using either the lifting presentation or the extension
presentation). It is clear from the construction that locally on $B$,
the stack $\tXc{\alpha}{\beta}$ is isomorphic to $X\times
B\GM\times\ZZZ$; globally, it carries a natural filtration
$0\subset{\widetilde{X}}^{(1)}\subset{\widetilde{X}}^{(2)}
\subset\tXc{\alpha}{\beta}$ with $\widetilde{X}^{(1)}=B\GM$,
$\widetilde{X}^{(2)}/\widetilde{X}^{(1)}=X$, and
$\tXc{\alpha}{\beta}/{\widetilde X}^{(2)}=\ZZZ$. This implies the
following statement:

\PROPOSITION{$\tXc{\alpha}{\beta}$ is ``good'' in the sense of
(\ref{goodstacks}).}  \PROOF{The property of being ``good'' is local
on $B$, so it is enough to notice that the stacks $X$, $B\GM$, and
$\ZZZ$ are ``good''.}

In particular, we see that the Fourier-Mukai transform gives an
equivalence between $D^b(\tXc{\alpha}{\beta})$ and
$D^b(\tXc{-\beta}{\alpha})$.

\subsection{Duality for torsors}
\label{dualtorsors}

Now suppose $\calX$ is a commutative group stack which is ``good'',
and let $\calX'$ be a torsor over $\calX$. Denote by $\widetilde{\calX}$
the corresponding extension of $\ZZZ$ by $\calX$: it fits into the
exact sequence
\begin{equation*}
0\to\calX\to\widetilde\calX\to\ZZZ\to0
\end{equation*}
and $\calX'$ is identified with the preimage of $1\in\ZZZ$. Notice
that locally on $B$, the torsor is trivial, so $\widetilde{\calX}$ is
isomorphic to $\calX\times\ZZZ$. Since both $\calX$ and $\ZZZ$ are
``good'', so is $\widetilde{\calX}$. The dual stack
$\widetilde{\calX}^\vee$ is isomorphic to $\calX^\vee\times B\GM$
locally on $B$ (globally, it contains a substack isomorphic to $B\GM$,
and the quotient equals $\calX^\vee$). In particular, if $\calX^\vee$
is actually a space (rather than a stack), then
$\widetilde{\calX}^\vee$ is an $\calO^\times$-gerbe over the
space. The following statement is clear:

\PROPOSITION{The Fourier-Mukai functor $\FM:D^b(\widetilde{\calX})\to
D^b(\widetilde{\calX}^\vee)$ induces an equivalence $D^b(\calX')\to
D^b_1(\widetilde{\calX}^\vee)$, where
$D^b_1(\widetilde{\calX}^\vee)\subset D^b(\widetilde{\calX}^\vee)$ is the
complete subcategory of objects $F$ such that the action of $\GM$ on
$H^i(F)$ is tautological. Here the action is induced by the morphism
$B\GM\to\widetilde\calX^\vee$.}

In the case of duality between $\tXc{\alpha}{\beta}$ and
$\tXc{-\beta}{\alpha}$, both stacks are torsors (over
$\Xc{\alpha}{\beta}$ and $\Xc{-\beta}{\alpha}$, respectively), and so
we get equivalences
\begin{equation*}
D^b(\Xc{\alpha}{\beta})\to D^b_1(\tXc{-\beta}{\alpha})
\end{equation*}
and
\begin{equation*}
D^b_1(\tXc{\alpha}{\beta})\to D^b(\Xc{-\beta}{\alpha}).
\end{equation*}
They induce an equivalence
\begin{equation*}
D^b(\Xc\alpha\beta)
\cap D^b_1(\tXc{\alpha}{\beta})=
D^b_1(\Xc{\alpha}{\beta})\to D^b_1(\Xc{-\beta}{\alpha}).
\end{equation*}
This is exactly what Theorem B claims.

\


\bigskip

\noindent


\begin{thebibliography}{DOPW01b}

\bibitem[AD98]{aspinwall-donagi}
P.~Aspinwall and R.~Donagi.
\newblock The heterotic string, the tangent bundle and derived categories.
\newblock {\em Adv. Theor. Math. Phys.}, 2(5):1041--1074, 1998.

\bibitem[Art70]{artin-a2}
M.~Artin.
\newblock Algebraization of formal moduli {I}{I}. {Existence of modifications}.
\newblock {\em Ann. of Math.}, 91:88--135, 1970.

\bibitem[Art74]{artin-stacks}
M.~Artin.
\newblock Versal deformations and algebraic stacks.
\newblock {\em Invent. math.}, 27:165--189, 1974.

\bibitem[BBRP98]{bbhrmp}
C.~Bartocci, U.~Bruzzo, D.~Hern{\'a}ndez Ruip{\'e}rez, and J.~Mu{\~n}oz Porras.
\newblock Mirror symmetry on {$K3$} surfaces via {F}ourier-{M}ukai transform.
\newblock {\em Comm. Math. Phys.}, 195(1):79--93, 1998.

\bibitem[BJPS97]{bjps}
M.~Bershadsky, A.~Johansen, T.~Pantev, and V.~Sadov.
\newblock On four-dimensional compactifications of {$F$}-theory.
\newblock {\em Nuclear Phys. B}, 505(1-2):165--201, 1997.

\bibitem[BK90]{bondal-kapranov}
A.~Bondal and M.~Kapranov.
\newblock Representable functors, {S}erre functors and mutations.
\newblock {\em Math. USSR. Izv.}, 35:519--541, 1990.
\newblock alg-geom/9506012.

\bibitem[BKR01]{bkr}
T.~Bridgeland, A.~King, and M.~Ried.
\newblock Mukai implies {M}c{K}ay: the {M}c{K}ay correspondence as an
  equivalence of derived categories.
\newblock {\em J. Amer. Math. Soc.}, 14:535--554, 2001.
\newblock math.AG/9908027.

\bibitem[BM02]{bridgeland-maciocia}
T.~Bridgeland and A.~Maciocia.
\newblock Fourier-{M}ukai transforms for {$K3$} and elliptic fibrations.
\newblock {\em J. Algebraic Geom.}, 11(4):629--657, 2002.

\bibitem[BM03a]{ruxandra2}
V.~Brinzanescu and R.~Moraru.
\newblock {Holomorphic rank-2 vector bundles on non-Kahler elliptic surfaces},
  2003, arXiv:math.AG/0306191.

\bibitem[BM03b]{ruxandra1}
V.~Brinzanescu and R.~Moraru.
\newblock {Stable bundles on non-{K}ahler elliptic surfaces}, 2003,
  arXiv:math.AG/0306192.

\bibitem[BO95]{bondal-orlov-flops}
A.~Bondal and D.~Orlov.
\newblock {Semiorthogonal decomposition for algebraic varieties}, 1995,
  arXiv:alg-geom/9506012.

\bibitem[Bre90]{breen-bit}
L.~Breen.
\newblock Bitorseurs et cohomologie non ab\'elienne.
\newblock In {\em The Grothendieck Festschrift, Vol.\ I}, pages 401--476.
  Birkh\"auser Boston, Boston, MA, 1990.

\bibitem[Bre94]{breen-2g}
L.~Breen.
\newblock On the classification of $2$-gerbes and $2$-stacks.
\newblock {\em Ast\'erisque}, 225:160, 1994.

\bibitem[Bri98]{bridgeland-elliptic}
T.~Bridgeland.
\newblock Fourier-{M}ukai transforms for elliptic surfaces.
\newblock {\em J. Reine Angew. Math.}, 498:115--133, 1998.

\bibitem[Bri99]{bridgeland}
T.~Bridgeland.
\newblock Equivalences of triangulated categories and {F}ourier-{M}ukai
  transforms.
\newblock {\em Bull. London Math. Soc.}, 31:25--34, 1999.
\newblock math.AG/9809114.

\bibitem[Bri02]{bridgeland-flops}
T.~Bridgeland.
\newblock Flops and derived categories.
\newblock {\em Invent. Math.}, 147(3):613--632, 2002.

\bibitem[Bry93]{brylinski}
J.-L. Brylinski.
\newblock {\em Loop spaces, characteristic classes and geometric quantization},
  volume 107 of {\em Progress in Mathematics}.
\newblock Birkh\"auser Boston Inc., Boston, MA, 1993.

\bibitem[C{\u{a}}l00]{caldararu-thesis}
A.~C{\u{a}}ld{\u{a}}raru.
\newblock {\em Derived categories of twisted sheaves on {C}alabi-{Y}au
  manifolds}.
\newblock PhD thesis, Cornell University, 2000.

\bibitem[C{\u{a}}l01]{caldararu-k3}
A.~C{\u{a}}ld{\u{a}}raru.
\newblock {Non-Fine Moduli Spaces of Sheaves on $K3$ Surfaces}, 2001,
  arXiv:math.AG/0108180.

\bibitem[C{\u{a}}l02]{caldararu-cy}
A.~C{\u{a}}ld{\u{a}}raru.
\newblock Derived categories of twisted sheaves on elliptic threefolds.
\newblock {\em J. Reine Angew. Math.}, 544:161--179, 2002.

\bibitem[Del74]{deligne-hodge3}
P.~Deligne.
\newblock Th\'eorie de {H}odge. {III}.
\newblock {\em Inst. Hautes \'Etudes Sci. Publ. Math.}, 44:5--77, 1974.

\bibitem[Del90]{deligne-tannaka}
P.~Deligne.
\newblock Cat\'egories tannakiennes.
\newblock In {\em The Grothendieck Festschrift, Vol.\ II}, pages 111--195.
  Birkh\"auser Boston, Boston, MA, 1990.

\bibitem[DG94]{dolgachev-gross}
I.~Dolgachev and M.~Gross.
\newblock Elliptic threefolds. {I}. {O}gg-{S}hafarevich theory.
\newblock {\em J. Algebraic Geom.}, 3(1):39--80, 1994.

\bibitem[DG02]{donagi-gaitsgory}
R.~Donagi and D.~Gaitsgory.
\newblock The gerbe of {H}iggs bundles.
\newblock {\em Transform. Groups}, 7(2):109--153, 2002, arXiv:math.AG/0005132.

\bibitem[Don97]{donagi-principal}
R.~Donagi.
\newblock Principal bundles on elliptic fibrations.
\newblock {\em Asian J. Math.}, 1(2):214--223, 1997.

\bibitem[Don99]{donagi-heteroticF}
R.~Donagi.
\newblock Heterotic {F}-theory duality.
\newblock In {\em XIIth International Congress of Mathematical Physics (ICMP
  '97) (Brisbane)}, pages 206--213. Internat. Press, Cambridge, MA, 1999.

\bibitem[DOPW01a]{smb}
R.~Donagi, B.~Ovrut, T.~Pantev, and D.~Waldram.
\newblock Standard-model bundles.
\newblock {\em Adv. Theor. Math. Phys.}, 5(3):563--615, 2001.

\bibitem[DOPW01b]{sm}
R.~Donagi, B.~Ovrut, T.~Pantev, and D.~Waldram.
\newblock Standard models from heterotic {M}-theory.
\newblock {\em Adv. Theor. Math. Phys.}, 5(1):93--137, 2001.

\bibitem[EHKV01]{ehkv}
D.~Edidin, B.~Hassett, A.~Kresch, and A.~Vistoli.
\newblock Brauer groups and quotient stacks.
\newblock {\em Amer. J. Math.}, 123(4):761--777, 2001.

\bibitem[FMW97]{fmw}
R.~Friedman, J.~Morgan, and E.~Witten.
\newblock Vector bundles and {${\rm F}$}-theory.
\newblock {\em Comm. Math. Phys.}, 187(3):679--743, 1997.

\bibitem[FMW98]{fmw-pb}
R.~Friedman, J.~Morgan, and E.~Witten.
\newblock Principal {$G$}-bundles over elliptic curves.
\newblock {\em Math. Res. Lett.}, 5(1-2):97--118, 1998.

\bibitem[FMW99]{fmw-vb}
R.~Friedman, J.~Morgan, and E.~Witten.
\newblock Vector bundles over elliptic fibrations.
\newblock {\em J. Algebraic Geom.}, 8(2):279--401, 1999.

\bibitem[Gab81]{gabber}
O.~Gabber.
\newblock Some theorems on {A}zumaya algebras.
\newblock In {\em The Brauer group (Sem., Les Plans-sur-Bex, 1980)}, pages
  129--209. Springer, Berlin, 1981.

\bibitem[Gir71]{giraud}
J.~Giraud.
\newblock {\em Cohomologie non ab\'elienne}.
\newblock Springer-Verlag, Berlin, 1971.
\newblock Die Grundlehren der mathematischen Wissenschaften, Band 179.

\bibitem[GMS00]{gms}
O.~Ganor, A.~Mikhailov, and N.~Saulina.
\newblock Constructions of noncommutative instantons on {$T\sp 4$} and {$K\sb
  3$}.
\newblock {\em Nuclear Phys. B}, 591(1-2):547--583, 2000.

\bibitem[Gra91]{antonella}
A.~Grassi.
\newblock On minimal models of elliptic threefolds.
\newblock {\em Math. Ann.}, 290(2):287--301, 1991.

\bibitem[Gro68a]{groth-bg1}
A.~Grothendieck.
\newblock Le groupe de {B}rauer. {I},. {A}lg\`ebres d'{A}zumaya et
  interpr\'etations diverses.
\newblock In {\em Dix Expos\'es sur la Cohomologie des Sch\'emas}, pages
  46--66. North-Holland, Amsterdam, 1968.

\bibitem[Gro68b]{groth-bg2}
A.~Grothendieck.
\newblock Le groupe de {B}rauer. {I}{I}. {T}h\'eorie cohomologique.
\newblock In {\em Dix Expos\'es sur la Cohomologie des Sch\'emas}, pages
  67--87. North-Holland, Amsterdam, 1968.

\bibitem[Gro68c]{groth-bg3}
A.~Grothendieck.
\newblock Le groupe de {B}rauer. {I}{I}{I}. {E}xemples et compl\'ements.
\newblock In {\em Dix Expos\'es sur la Cohomologie des Sch\'emas}, pages
  88--188. North-Holland, Amsterdam, 1968.

\bibitem[Gro98]{gross-slag1}
M.~Gross.
\newblock Special {L}agrangian fibrations. {I}. {T}opology.
\newblock In {\em Integrable systems and algebraic geometry (Kobe/Kyoto,
  1997)}, pages 156--193. World Sci. Publishing, River Edge, NJ, 1998.

\bibitem[Gro99]{gross-slag2}
M.~Gross.
\newblock Special {L}agrangian fibrations. {II}. {G}eometry. {A} survey of
  techniques in the study of special {L}agrangian fibrations.
\newblock In {\em Surveys in differential geometry: differential geometry
  inspired by string theory}, volume~5 of {\em Surv. Differ. Geom.}, pages
  341--403. Int. Press, Boston, MA, 1999.

\bibitem[Har66]{hartshorne-rd}
R.~Hartshorne.
\newblock Residues and duality.
\newblock Lecture Notes in Mathematics, No. 20, Springer-Verlag, 1966.

\bibitem[Hit01]{hitchin-slag}
N.~Hitchin.
\newblock Lectures on special {L}agrangian submanifolds.
\newblock In {\em Winter School on Mirror Symmetry, Vector Bundles and
  Lagrangian Submanifolds (Cambridge, MA, 1999)}, volume~23 of {\em AMS/IP
  Stud. Adv. Math.}, pages 151--182. Amer. Math. Soc., Providence, RI, 2001.

\bibitem[Hoo82]{hoobler}
R.~Hoobler.
\newblock When is ${{\rm {B}r}}({X})={{\rm {B}r}}\sp{\prime} ({X})$?
\newblock In {\em Brauer groups in ring theory and algebraic geometry (Wilrijk,
  1981)}, pages 231--244. Springer, Berlin, 1982.

\bibitem[HS03]{hyubrechts-schroeer}
D.~Huybrechts and S.~Schroeer.
\newblock {The {B}rauer group of analytic $K3$ surfaces}, 2003,
  arXiv:math.AG/0305101.
\newblock to appear in Int. Math. Res. Not.

\bibitem[IN99]{ito-nakamura}
Y.~Ito and I.~Nakamura.
\newblock Hilbert schemes and simple singularities.
\newblock In {\em New trends in algebraic geometry (Warwick, 1996)}, volume 264
  of {\em London Math. Soc. Lecture Note Ser.}, pages 151--233. Cambridge Univ.
  Press, Cambridge, 1999.

\bibitem[Kaw02]{kawamata-flops}
Yujiro Kawamata.
\newblock {D-equivalence and K-equivalence}, 2002, arXiv:math.AG/0205287.

\bibitem[KKO01]{kko}
A.~Kapustin, A.~Kuznetsov, and D.~Orlov.
\newblock Noncommutative instantons and twistor transform.
\newblock {\em Comm. Math. Phys.}, 221(2):385--432, 2001.

\bibitem[Kod63]{kodaira}
K.~Kodaira.
\newblock On compact analytic surfaces. {I}{I}, {I}{I}{I}.
\newblock {\em Ann. of Math. (2) 77 (1963), 563--626; ibid.}, 78:1--40, 1963.

\bibitem[Kon01]{kontsevich}
M.~Kontsevich.
\newblock Deformation quantization of algebraic varieties.
\newblock {\em Lett. Math. Phys.}, 56(3):271--294, 2001.
\newblock EuroConf\'erence Mosh\'e Flato 2000, Part III (Dijon).

\bibitem[LMB00]{laumon-stacks}
G.~Laumon and L.~Moret-Bailly.
\newblock {\em Champs alg\'ebriques}.
\newblock Springer-Verlag, Berlin, 2000.

\bibitem[Mil80]{milne-book}
J.~Milne.
\newblock {\em Etale cohomology}.
\newblock Princeton University Press, Princeton, N.J., 1980.

\bibitem[Mir83]{miranda}
R.~Miranda.
\newblock Smooth models for elliptic threefolds.
\newblock In {\em The birational geometry of degenerations (Cambridge, Mass.,
  1981)}, pages 85--133. Birkh\"auser Boston, Mass., 1983.

\bibitem[Muk81]{mukai}
S.~Mukai.
\newblock Duality between {$D(X)$} and {$D(\hat X)$} with its application to
  {P}icard sheaves.
\newblock {\em Nagoya Math. J.}, 81:153--175, 1981.

\bibitem[Nak01]{nakayama-global}
N.~Nakayama.
\newblock Global structure of an elliptic fibration, 2001,
  www.kurims.kyoto-u.ac.jp/home\_ page/ preprint/PS/RIMS1322.pdf.

\bibitem[NS98]{ns}
N.~Nekrasov and A.~Schwarz.
\newblock Instantons on noncommutative {${\mathbb R}^{4}$}, and {$(2,0)$}
  superconformal six-dimensional theory.
\newblock {\em Comm. Math. Phys.}, 198(3):689--703, 1998.

\bibitem[Orl97]{orlov-k3}
D.~Orlov.
\newblock Equivalences of derived categories and {$K3$} surfaces.
\newblock {\em J. Math. Sci. (New York)}, 84(5):1361--1381, 1997.
\newblock Algebraic geometry, 7.

\bibitem[Orl02]{orlov-abelian}
D.~Orlov.
\newblock Derived categories of coherent sheaves on abelian varieties and
  equivalences between them.
\newblock {\em Izv. Ross. Akad. Nauk Ser. Mat.}, 66(3):131--158, 2002.

\bibitem[Pol96]{polishchuk-biextensions}
A.~Polishchuk.
\newblock Symplectic biextensions and a generalization of the {F}ourier-{M}ukai
  transform.
\newblock {\em Math. Res. Lett.}, 3(6):813--828, 1996.

\bibitem[Pol02]{polishchuk-weil}
A.~Polishchuk.
\newblock Analogue of {W}eil representation for abelian schemes.
\newblock {\em J. Reine Angew. Math.}, 543:1--37, 2002.

\bibitem[Ray70]{raynaud-picard}
M.~Raynaud.
\newblock Sp\'ecialisation du foncteur de {P}icard.
\newblock {\em Inst. Hautes \'Etudes Sci. Publ. Math.}, 38:27--76, 1970.

\bibitem[Sch01]{schroer}
S.~Schr{\"o}er.
\newblock There are enough {A}zumaya algebras on surfaces.
\newblock {\em Math. Ann.}, 321(2):439--454, 2001.

\bibitem[{\mbox{SGA7-I}}]{sga7}
Groupes de monodromie en g\'eom\'etrie alg\'ebrique. {I}.
\newblock Lecture Notes in Mathematics, Vol. 288, Springer, 1972.
\newblock S\'eminaire de G\'eom\'etrie Alg\'ebrique du Bois-Marie 1967--1969
  (SGA 7 I), Dirig\'e par A. Grothendieck. Avec la collaboration de M. Raynaud
  et D. S. Rim.

\bibitem[ST01]{seidel-thomas}
P.~Seidel and R.~Thomas.
\newblock Braid group actions on derived categories of coherent sheaves.
\newblock {\em Duke Math. J.}, 108(1):37--108, 2001.

\bibitem[SYZ96]{SYZ}
A.~Strominger, S.-T. Yau, and E.~Zaslow.
\newblock Mirror symmetry is {$T$}-duality.
\newblock {\em Nuclear Phys. B}, 479(1-2):243--259, 1996.

\end{thebibliography}
\end{document}